\documentclass[pdflatex,sn-mathphys-num]{sn-jnl}

\usepackage[utf8]{inputenc}
\usepackage[T1]{fontenc}
\usepackage{graphicx}%
\usepackage{multirow}%
\usepackage{amsmath}%
\usepackage{mathtools}%
\usepackage{amssymb,amsfonts}%
\usepackage{amsthm}%
\usepackage{mathrsfs}%
\usepackage[title]{appendix}%
\usepackage{xcolor}%
\usepackage{textcomp}%
\usepackage{manyfoot}%
\usepackage{booktabs}%
\usepackage{algorithm}%
\usepackage{algorithmicx}%
\usepackage{algpseudocode}%
\usepackage{listings}%
\usepackage{placeins}%
\usepackage{epsfig}%
\usepackage{tabularx}%
\newcolumntype{Y}{>{\centering\arraybackslash}X}
\usepackage{float}%
\usepackage{hyperref}%

\theoremstyle{thmstyleone}%
\newtheorem{theorem}{Theorem}
\newtheorem{proposition}[theorem]{Proposition}%
\newtheorem{assumption}[theorem]{Assumption}
\newtheorem{lemma}[theorem]{Lemma}
\newtheorem{corollary}[theorem]{Corollary}

\theoremstyle{thmstyletwo}%

\theoremstyle{thmstylethree}%

\begin{document}

\title[Nonlinear Dimensionality Reduction Techniques for Bayesian Optimization]{Nonlinear Dimensionality Reduction Techniques for Bayesian Optimization}


\author[1]{\fnm{Luo} \sur{Long}}\email{luo.long@maths.ox.ac.uk}

\author[1]{\fnm{Coralia} \sur{Cartis}}\email{cartis@maths.ox.ac.uk}

\author[1]{\fnm{Paz} \sur{Fink Shustin}}\email{paz.finkshustin@maths.ox.ac.uk}

\affil[1]{\orgdiv{Mathematical Institute}, \orgname{University of Oxford}, \orgaddress{\street{Radcliffe Observatory, Andrew Wiles Building, Woodstock Rd}, \city{Oxford}, \postcode{OX2 6GG}, \country{United Kingdom}}}

\abstract{Bayesian optimisation (BO) is a standard approach for sample-efficient global optimisation of expensive black-box functions, yet its scalability to high dimensions remains challenging. Here we investigate nonlinear dimensionality reduction techniques that reduce the problem to a sequence of low-dimensional Latent-Space BO (LSBO). While early LSBO methods used (linear) random and supervised embeddings (\citet{Wang2013} and \citet{egorse}), building on \citet{grosnit2021highdimensionalbayesianoptimisationvariational}, we employ Variational Autoencoders (VAEs) for LSBO, focusing on deep metric loss for structured latent manifolds and VAE retraining to adapt the encoder–decoder pair to newly sampled regions.
We then couple LSBO with Sequential Domain Reduction (SDR) directly in latent space (SDR-LSBO), yielding an algorithm that narrows the latent search domains as evidence accumulates. Implemented in a GPU-accelerated \textit{BoTorch} stack with Matérn-5/2 Gaussian-process surrogates, our numerical results show improved optimisation quality across benchmark tasks and that retraining can improve BO performance. Additionally, we compare adaptive supervised linear random embeddings against VAEs as two mechanisms for dimensionality reduction, showing that combining VAEs with BO effectively addresses problems with nonlinear low-dimensional structures. For theory, we analyse BO--VAE with a fixed pretrained representation and decompose its ambient-space simple regret into a latent BO error and a fixed VAE-induced representation gap. Under a PAC--Bayes-certified reconstruction condition and standard fixed-prior assumptions for expected improvement with a Mat\'ern--5/2 kernel, the latent BO error vanishes as the evaluation budget increases, whereas the representation gap remains fixed and may impose a non-vanishing error floor. We complement this analysis with visualisations that empirically examine the accessibility of the ambient optimum through the learned decoder. To the best of our knowledge, this is the first study to combine SDR with VAE-based LSBO, and our analysis clarifies design choices for metric shaping and retraining that are critical for scalable latent-space BO. For reproducibility, our source code is available at \hyperlink{https://github.com/L-Lok/Nonlinear-Dimensionality-Reduction-Techniques-for-Bayesian-Optimization.git}{this link}.}


\maketitle

\section{Introduction}\label{sec1}

Global optimisation aims to find the (approximate) global optimum of a smooth function $f$ within a region of interest, possibly without the use of derivative information and with careful handling of often-costly objective evaluations. In this paper, we consider the derivative-free global optimisation problem
\begin{equation}\label{main problem}
    \operatorname*{minimize}_{\mathbf{x}\in\Omega} f(\mathbf{x}),
    \tag{P}
\end{equation}
where \(\Omega\subseteq\mathbb R^D\) is the problem domain, possibly \(\Omega=\mathbb R^D\) in the unconstrained case, and $f$ is Lipschitz continuous. We write
\[
    f_\Omega^\ast
    \coloneq
    \inf_{\mathbf x\in\Omega}f(\mathbf x).
\]
When this infimum is attained, the set of global minimisers is denoted by
\[
    S_\Omega^\ast
    \coloneq
    \operatorname*{arg\,min}_{\mathbf x\in\Omega}f(\mathbf x).
\]

To address such objective functions, we adopt the Bayesian Optimisation (BO) approach in this work. \\

BO is a state-of-the-art global optimisation framework that typically models \(f\) with a Gaussian-process (GP) surrogate. At each iteration, an acquisition function chooses the next evaluation point \(\mathbf{x}_{n+1}\). After computing \(f(\mathbf{x}_{n+1})\), the pair \((\mathbf{x}_{n+1}, f(\mathbf{x}_{n+1}))\) is added to the data set and the GP posterior is updated, enabling sample-efficient search for the global minimiser via a probabilistic surrogate \cite{Frazier2018}. The performance of BO depends on the acquisition functions that balance exploitation and exploration during the BO search. The former considers the areas with lower posterior mean, while the latter prefers areas with higher posterior variance.

\paragraph{BO with sequential domain reduction.}
Since BO requires acquisition optimisation over a bounded region, the algorithms in this work are applied on a prescribed initial search box
\[
    \mathcal X \coloneq [\mathbf a,\mathbf b]
    = \{\mathbf x\in\mathbb R^D : a_i\le x_i\le b_i,\ i=1,\ldots,D\}
    \subseteq \Omega,
\]
where \(\mathbf a,\mathbf b\in\mathbb R^D\) and \(a_i, b_i \in \mathbb{R}, a_i<b_i, \forall i \in [1, D]\), for Problem \eqref{main problem}. In the unconstrained instances, \(\mathcal X\) is a computational search region chosen from prior knowledge or problem scaling, large enough not to exclude the global minimiser. For the numerical experiments in Section \ref{Results}, it is the standard box specified for each benchmark problem. Thus, throughout this work, the computational domain is selected so that it contains at least one ambient global minimiser:
\begin{equation}
    S_\Omega^\ast\cap\mathcal X\neq\varnothing.
\label{eq:computational-domain-contains-optimum}
\end{equation}
Define
\[
    f_{\mathcal X}^\ast
    \coloneq
    \min_{\mathbf x\in\mathcal X}f(\mathbf x),
    \qquad
    S_{\mathcal X}^\ast
    \coloneq
    \operatorname*{arg\,min}_{\mathbf x\in\mathcal X}f(\mathbf x).
\]
Condition~\eqref{eq:computational-domain-contains-optimum} implies $f_{\mathcal X}^\ast=f_\Omega^\ast.$ For completeness, one could nevertheless define the computational-domain restriction gap by
\[
    \Delta_{\mathcal X}
    \coloneq
    f_{\mathcal X}^\ast-f_\Omega^\ast,
\]
but \(\Delta_{\mathcal X}=0\) under our standing convention. While BO is highly effective for expensive black-box optimisation, its performance is often affected by high ambient dimension and by overly large search domains \cite{hvarfner2024vanillabayesianoptimizationperforms, Shahriari2016}. Such domains, however, exacerbate the exploration burden and the difficulty of acquisition optimisation, often yielding slow convergence and increased computational cost. To address the exploration burden associated with a large search domain, we integrate Sequential Domain Reduction (SDR), a response-surface strategy that maintains and adaptively updates a region of interest. Let $\mathcal{R}_k\subseteq\mathcal{X}$ denote the region used at BO iteration $k$. The SDR update recentres $\mathcal{R}_k$ around the current incumbent and rescales its coordinatewise side lengths according to the observed incumbent movement. Thus, the search is progressively localised around promising solutions while the surrogate--acquisition loop remains unchanged. The implementation used in this work is a panning--zooming procedure rather than a strictly nested-domain method. In particular, although every region remains inside the initial domain,
\begin{equation}
    \mathcal{R}_k\subseteq\mathcal{X}
    \qquad\text{for every }k,
\end{equation}
we do not require $\mathcal{R}_{k+1}\subseteq\mathcal{R}_k$. The region may translate when the incumbent moves and may temporarily expand in some coordinates. Once the incumbent stabilises, the updates typically reduce the scale of the region and concentrate evaluations near the incumbent. We therefore use ``domain reduction'' to mean adaptive localisation of the search, rather than strict set-theoretic nesting. Further details are provided in Appendix~C and \cite{SDR}. As an illustration, Figure~\ref{fig:ackley-sdr} reports a 10-dimensional Ackley benchmark comparing BO with and without SDR. The increasing gap between the two curves indicates an improvement from using SDR.

\begin{figure}[H]
    \centering
    \includegraphics[width=\linewidth]{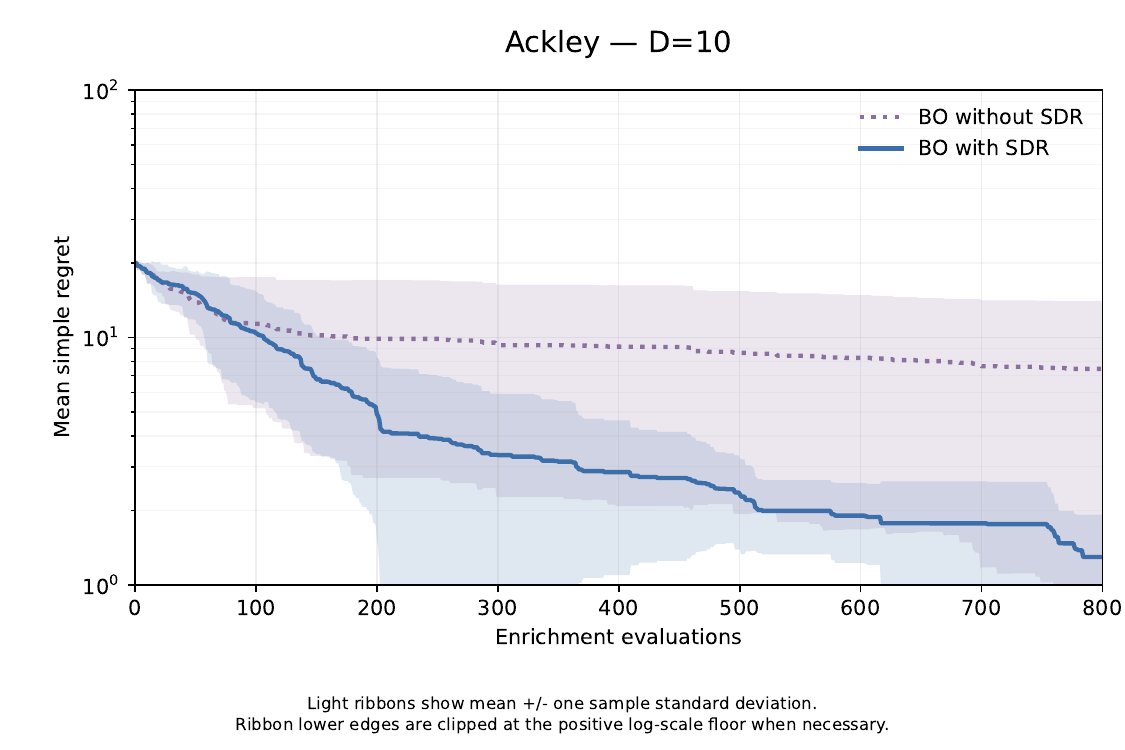}
    \caption{Effect of SDR on BO for the 10-dimensional Ackley function. SDR adaptively recentres and rescales the region of interest as evidence accumulates. $10$ initial design points, $5$ distinct seeds, and $800$ enrichment iterations are used. In this example, it yields improved performance relative to vanilla BO.}
    \label{fig:ackley-sdr}
\end{figure}

However, as noted by \cite{SDR} and our experiments, early BO decisions can be strongly affected by the initial design and by the surrogate uncertainty when only a small number of evaluations have been made. If SDR is applied too frequently, the current region of interest may therefore be contracted around a suboptimal incumbent before the acquisition rule has explored the current region sufficiently. The global minimiser may then lie outside the current region. Under a strictly nested update \(\mathcal X_0 \supseteq \mathcal X_1 \supseteq \cdots\), this exclusion would be irreversible, i.e., if a global minimiser \(\mathbf{x}^\ast\) is excluded at some update, that is, \(\mathbf{x}^\ast\notin \mathcal X_k\), then \(\mathbf{x}^\ast\notin \mathcal X_{k'}\) for every \(k'\ge k\). However, under the panning--zooming implementation adopted here, a later translation or expansion can in principle reintroduce the excluded region, but such recovery is not guaranteed and may become increasingly difficult as the search concentrates around the incumbent. We therefore update SDR only after every $K_{\mathrm{SDR}}$ BO evaluations, allowing additional exploration between successive region updates.

\paragraph{Dimensionality reduction techniques for BO.}
As the dimensionality $D$ of the Problem~\eqref{main problem} grows, the scalability of BO degrades \cite{hvarfner2024vanillabayesianoptimizationperforms}. Although SDR can enhance the robustness of BO, its effectiveness diminishes in high-dimensional settings. A common remedy is through dimensionality reduction techniques \cite{Moriconi2020quantileGP, Rolland2018, Moriconi2020}. The motivation of dimensionality reduction methods is to map the Problem \eqref{main problem} from a high-dimensional space into a lower-dimensional latent subspace so that BO can operate more efficiently, leading to what is commonly termed the Latent Space Bayesian Optimisation (LSBO). Although there are many existing dimensionality reduction methods such as the linear principal component analysis \cite{Wold1987} and the nonlinear t-Distributed Stochastic Neighbor Embedding \cite{Maaten2008}, we suggest thinking of them as an \textit{Encoder-Decoder} framework, in which the \textit{encoder} denotes the process that produces the latent representation given the original high-dimensional data (by feature selection or extraction) and the \textit{decoder} performs the reverse mapping. Thus, dimensionality reduction can be interpreted as a process of data compression where the encoder encodes (compresses) the original data from the ambient space to the latent subspace\footnote{Equivalent names are low-dimensional subspace and encoded space.} and then the decoder decodes (decompresses) it. For example, in \cite{Wang2013, egorse}, random linear mappings are used to reduce the input dimension, and then BO is performed over the latent subspace. However, they are restricted to linear mappings. When dealing with nonlinearities and complex data distributions, deep generative models \cite{Bond_Taylor_2022} are frequently employed for LSBO. For example, Variational AutoEncoders (VAEs) \cite{kingma2022autoencodingvariationalbayes} use neural networks as their encoders and decoders. By incorporating nonlinear activation functions, the encoder acts as a nonlinear mapping, capable of creating general latent data manifolds. For this reason, in this work, we primarily adopt VAEs to learn nonlinear embeddings, and we compare against random embeddings to isolate the effect of the dimensionality reduction choice.\\

SDR still yields notable gains for BO in moderate dimensions,  but its effectiveness can deteriorate as the ambient dimension $D$ increases. This motivates applying SDR after mapping the Problem \eqref{main problem} to a low-dimensional latent space $\mathcal{Z} \subseteq \mathbb{R}^d$, where $d < D,$ to accelerate BO and reduce computational burden. Thus, the region of interest is contracted in the latent space rather than in the original variable space. To reduce the risk of premature contraction, we allow the SDR region of interest to be updated only periodically, after $K$ iterations of BO ($K \geq 1$), rather than necessarily after every evaluation, so that more evaluations are made within the current bounds before contraction. This SDR update period is distinct from the VAE retraining period $q$, used only in the retraining variants and defined with the algorithms in Section \ref{Alg}.

\subsection{Related Work}
BO has been widely studied as a black-box optimisation approach. For general background on derivative-free and black-box optimisation, see \cite{audet2026derivativefree, AudetHare2017} and \cite{Conn2009}. For comprehensive discussions of acquisition functions, including \textit{probability of improvement} \cite{Kushner1964}, \textit{upper confidence bounds} \cite{Cox1997}, and \textit{expected improvement} \cite{Mockus1978}, we refer the reader to more tutorial treatments \cite{Shahriari2016, Frazier2018, Rasmussen2006}. For recent theoretical results, we recommend \cite{srinivas2010ucb, bull2011ego, pmlr-v77-nguyen17a, teckentrup2020convergencegaussianprocessregression}. It is known that BO is hard to scale up to high dimensions. This motivates dimensionality reduction schemes that map Problem \eqref{main problem} to a lower-dimensional search space and exploit structure. Early work assumes specific structure in \(f\): additivity/partial separability \cite{kandasamy2015additive,Rolland2018} or low effective dimension \(d\ll D\) \cite{Wang2013,constantine2015active,fornasier2012learningmutliridge}. Baseline treatments assume an axis-aligned effective subspace (i.e., some variables are inactive) \cite{bensalem2019sequential, chen2012joint}. For the general case of functions constant on an unknown linear subspace, prior work spans BO and related paradigms with three tactics: (i) learn the subspace (e.g., low-rank recovery) then optimise \cite{djolonga2013highdimgpbandits, tyagi2014learning}; (ii) alternate subspace estimation and optimisation \cite{garnett2013activelearninglinearembeddings, cartis2021globaloptimizationusingrandom}; or (iii) skip learning and optimise in randomly sampled low-dimensional subspaces given an estimate of the effective dimension \cite{binois2020choice, Wang2013}.\\

For the latter, \citet{Wang2013} developed the REMBO algorithm, primarily for low-rank functions. It tackles box-constrained BO by drawing a Gaussian random embedding matrix \(A\!\in\!\mathbb{R}^{D\times d}\) and optimising over a low-dimensional set \(\mathcal{Y}\subset\mathbb{R}^d\), with candidates mapped to the feasible set by
\[
\mathcal{Y} \supset y \;\mapsto\; x \;=\; p_{\mathcal X}(A y)\in\mathcal X \subset \mathbb{R}^D,
\]
where $p_{\mathcal{X}}: \mathbb{R}^D \mapsto \mathbb{R}^D.$ The success of REMBO hinges on the size of \(\mathcal{Y}\). When the embedded dimension \(d\) matches the effective dimension and the active subspace is axis-aligned, success probabilities can be quantified (Theorem 3 in \cite{Wang2013}). Under an encoder–decoder view (assuming \(A\) has orthonormal columns), the \emph{encoder} $A$ projects ambient points to the reduced space, the \emph{decoder} \(A^\top\) lifts latent points back to ambient, and \( p_{\mathcal X}\) enforces feasibility.\\

Extending the random-subspace idea of \citet{Wang2013}, \citet{cartis2020dimensionalityreductiontechniqueunconstrained} proposed REGO, a solver-agnostic framework that replaces the original problem with a Gaussian random, low-dimensional, bound-constrained reduced problem. They obtained probabilistic success bounds that depend only on the effective and embedding dimensions (not the ambient dimension) and identify the exact distribution of the reduced minimiser and its Euclidean norm. Later, \citet{cartis2021globaloptimizationusingrandom} extended this line via \emph{X-REGO}, which projects the problem sequentially or simultaneously onto Gaussian low-dimensional subspaces and optimises the reduced problems. Conic-integral-geometry tools give explicit success probabilities and hence global convergence guarantees under mild conditions.  For low effective dimension, they further devise an adaptive variant that increases the embedding dimension until the effective subspace is found, ensuring finite-embedding convergence, corroborated by numerical experiments.\\

More recently, \citet{egorse} proposed EGORSE, which extends purely random-embedding approaches by adaptively combining random linear embeddings, such as Gaussian projections, with supervised linear embeddings learned from the accumulated objective evaluations. In particular, partial least squares and marginal-Gaussian-process methods are used to identify influential linear directions, while random embeddings provide complementary exploration. For each embedding, EGORSE constructs a low-dimensional search box and formulates a feasibility constraint to ensure that candidate points can be mapped back into the original bounded domain. The resulting reduced problem is then solved using constrained BO. The numerical experiments in \cite{egorse} found the combination of partial least squares and Gaussian random embeddings to be particularly effective. Thus, EGORSE uses observed objective values to guide subspace selection more effectively than purely random methods, while retaining a low-dimensional linear-subspace representation.\\

However, the black-box objective functions are not necessarily of low-rank structures that can be well captured by random linear mapping, which is what these works primarily focus on. Although algorithms such as EGORSE and X-REGO can be applied to full-rank functions\footnote{We use the term ``full-rank function" to imply the objective functions do not have any low-rank properties.}, they may require restarts because the latent subspaces generated by random embeddings are linear and may fail to capture the optimum in full-rank settings \cite{Wang2013, cartis2021globaloptimizationusingrandom}. Hence,  we use deep generative models such as VAEs for dimensionality reduction within latent space optimisation, which is the main focus of this paper. The approach was first proposed in \cite{G_mez_Bombarelli_2018} for chemical design, which uses a VAE to embed discrete molecules into a continuous latent space and to decode back to valid structures, augmented by a property predictor. This continuous parametrisation permits latent-space exploration (random sampling, perturbations, interpolation) and gradient-based search for property-optimised molecules. \citet{Moriconi2020}\ later proposed a heuristic BO framework for high-dimensional optimisation by incorporating a nonlinear feature mapping $h: \mathbb{R}^D \to \mathbb{R}^d$ to reduce the dimensionality of the inputs, and a reconstruction mapping $g: \mathbb{R}^d \to \mathbb{R}^D$ based on GPs to evaluate the true objective function. This raises the question of how good the latent spaces should be to avoid invalid decoder outputs. In the context of LSBO, label guidance approaches are suggested by \citet{Siivola2020}, classified into joint and disjoint training.  Joint training means the VAE and the BO surrogate are optimised simultaneously, while disjoint training refers to optimising the two components separately. Joint training treats the VAE and GP models as a single machine by optimising the total loss consisting of the ELBO loss and the labelled data costs. However, \citet{Siivola2020} report that joint training can lead to overfitting and recommend disjoint training instead to yield improved BO performance. Meanwhile, \citet{Tripp2020}\ proposed a weighted retraining latent space optimisation framework based on deep generative models to address potential failures in common latent space optimisation frameworks, such as misalignment between VAE training and optimisation objectives. This framework ensures a well-structured latent space generated by deep generative models, where optimisation algorithms like BO can be used. An alternative method to construct a structured VAE latent space by incorporating Deep Metric Loss (DML) within the VAE training objective was proposed by \citet{grosnit2021highdimensionalbayesianoptimisationvariational}, which also combines weighted retraining techniques within the LSBO framework with applications to tasks such as molecule generation. However, these works provide limited discussion of generic problems with nonlinear low-dimensional structure or of how key parameters, such as the latent dimension, affect optimisation performance.\\

Related work has also explored domain-reduction methods for accelerating BO. For instance, \citet{salgia2021domainshrinkingbasedbayesianoptimization} proposed a domain reduction scheme based on a probabilistic threshold, which reduces the search domain by dividing it into sub-domains and predicting the one most likely to contain the global minimiser. Similarly, the SDR approach adopted here reduces the search region through a panning–zooming scheme. However, to the best of our knowledge, although SDR is an established global-optimisation technique, its application in a GPU-based environment and within VAE-generated latent spaces remains unexplored.\\

Finally, the theoretical foundations of BO with learned generative representations remain comparatively underdeveloped. \citet{grosnit2021highdimensionalbayesianoptimisationvariational} established sublinear regret for VAE-assisted BO, but their Assumption~3 requires that, after sufficiently many retraining epochs, there exist latent inputs from which the decoder recovers the unknown global optimiser with probability approaching one. This is a demanding representation condition because it directly assumes the asymptotic decoder-recovery property needed to transfer convergence from the latent space to the original problem. By contrast, the analyses of REMBO, REGO, and X-REGO in \cite{Wang2013,cartis2021globaloptimizationusingrandom}  concern fixed linear embeddings and therefore do not quantify the approximation error introduced by a learned nonlinear encoder--decoder. We address this gap by analysing BO--VAE with a fixed pretrained representation. Our result separates the ambient-space simple regret into a latent BO error and a VAE-induced representation gap. It shows that increasing the BO budget can eliminate the latent optimisation error while leaving a fixed, potentially nonzero representation floor. This identifies representation improvement as a necessary objective of VAE retraining, without assuming that an arbitrary retraining procedure necessarily achieves it.

\paragraph{Our aims and contributions.} 
Here we investigate scaling BO to high dimensions via dimensionality reduction\footnote{The work in this paper was part of the Master's thesis \cite{oxfordmasterthesis}, which is not published. A subset of this paper was included in a short workshop paper \cite{long2024dimensionalityreductiontechniquesglobal}, namely in the 2024 NeurIPS Workshop \textit{Optimization for ML}.}, focusing on (i) the learned nonlinear embeddings using VAEs and (ii) the adaptive random and supervised linear embeddings with BO, known as EGORSE \cite{egorse}. Building on \cite{grosnit2021highdimensionalbayesianoptimisationvariational}, whose framework is tailored to domain-specific workflows such as molecular design, we reformulate the BO–VAE pipeline and extend the algorithm to incorporate the Matérn-$5/2$ kernel to increase modelling flexibility and robustness for unconstrained global optimisation of generic black-box functions; see Section \ref{Algorithm Implementations}. In parallel, we introduce a variant of SDR that updates the region of interest only every $K$ BO iterations to allow sufficient exploration within the current bounds before contraction, and we deploy it both in the ambient space and in VAE-generated low-dimensional latent spaces within the GPU-based environment. Our main contributions are as follows:

\begin{enumerate}
    \item We propose three BO–VAE algorithms that integrate SDR within VAE latent spaces, exploiting latent regularity to guide contraction and improve optimisation. To our knowledge, this is the first integration of SDR within an LSBO framework. We study how VAE retraining \cite{Tripp2020} and DML \cite{grosnit2021highdimensionalbayesianoptimisationvariational} affect the baseline BO–VAE framework and its optimisation performance. We then conduct controlled comparisons across these algorithms by varying algorithm parameters such as latent dimension $d$ and ambient dimension $D$, and benchmark them against standard BO with SDR. On standard test functions, the BO-VAE approach remains effective in high-dimensional regimes. For fixed $D$, we empirically identify a trade-off between latent-space search difficulty and decoder support as $d$ varies.\\

    \item We conduct a comparative analysis of our BO--VAE algorithms and EGORSE \cite{egorse} on objective functions whose low-dimensional effective structure is induced by nonlinear mappings. This comparison contrasts the nonlinear latent representations learned by VAEs with EGORSE's adaptive use of Gaussian random projections and partial-least-squares-based supervised linear embeddings, thereby evaluating two distinct dimensionality-reduction strategies in terms of optimisation performance. On the nonlinear effective-mapping problems considered here, the VAE-based approaches generally achieve lower regret than EGORSE's adaptive combination of random and supervised linear embeddings.\\
    
    \item We develop a theoretical analysis of BO--VAE with a fixed, pretrained VAE, addressing the limited theoretical understanding of BO with learned generative representations. Our analysis decomposes the ambient-space optimisation error into a latent-space BO error and a VAE-induced representation gap. While the former vanishes as the evaluation budget increases under the stated assumptions, the latter remains fixed without VAE retraining and may impose a non-vanishing error floor relative to the true ambient-space optimum. Complementing this analysis, we provide diagnostic visualisations that empirically reveal these representation gaps and illustrate how they can prevent latent-space optimisation from reaching the ambient-space optimum.

\end{enumerate}

\paragraph{Paper outline.}

Section~\ref{Prelim} reviews BO, VAEs, and deep metric learning. Section~\ref{Alg} formulates decoder-induced latent-space optimisation and presents the three BO--VAE algorithms and the theoretical analysis of fixed-representation BO-VAE. Section~\ref{Results} presents the implementation details, empirical representation-gap findings, comparisons among the BO--VAE variants, and comparisons with EGORSE on problems with nonlinear effective mappings. Section~\ref{Conclusion} discusses the theoretical and empirical findings, limitations, and future directions.

\section{Preliminaries}\label{Prelim}

\subsection{Bayesian Optimisation}
As mentioned in the introduction, BO has two key ingredients: a probabilistic surrogate to model $f$ and acquisition functions to select new points. Let
\[
    \mathcal D_n \coloneq \{(\mathbf x_i,y_i)\}_{i=1}^n \subset (\mathcal{X} \times \mathbb R)^n,
    \qquad y_i \coloneq f(\mathbf x_i),
\]
denote the dataset after \(n\) noiseless evaluations of \(f\). We write \(\mathbf y_n=(y_1,\ldots,y_n)^\top\). A GP surrogate is a random function
\[
    F \sim \mathcal{GP}(m_0,k),
\]
where \(m_0: \mathcal X\to\mathbb R\) denotes the prior mean function and \(k\) is the covariance kernel. A zero prior mean is commonly used in practice. In this work, we use a Matérn-\(5/2\) kernel \cite{Rasmussen2006}. Next, we define \[
    K_n \coloneq [k(\mathbf x_i,\mathbf x_j)]_{i,j=1}^n,
    \qquad
    k_n(\mathbf x) \coloneq (k(\mathbf x_1,\mathbf x),\ldots,k(\mathbf x_n,\mathbf x))^\top,
\]
and
\[
    m_0(\mathbf X_n) \coloneq (m_0(\mathbf x_1),\ldots,m_0(\mathbf x_n))^\top.
\]
Conditioning the GP surrogate on \(\mathcal D_n\), \cite{Bishop2006} gives, for any
\(\mathbf x\in\mathcal X\)
\[
    F(\mathbf x)\mid \mathcal D_n
    \sim
    \mathcal N\bigl(m_n(\mathbf x),s_n^2(\mathbf x)\bigr),
\]
where
\[
    m_n(\mathbf x)
    \coloneq
    m_0(\mathbf x)
    +
    k_n(\mathbf x)^\top K_n^{-1}
    \bigl(\mathbf y_n-m_0(\mathbf X_n)\bigr),
\qquad
    s_n^2(\mathbf x)
    \coloneq
    k(\mathbf x,\mathbf x)
    -
    k_n(\mathbf x)^\top K_n^{-1}k_n(\mathbf x).
\]
Thus \(f\) remains the deterministic objective, while \(F\) denotes the GP surrogate used to model uncertainty about \(f\). For noisy observations \(y_i=f(\mathbf x_i)+\varepsilon_i\), with \(\varepsilon_i\sim\mathcal N(0,\sigma_\varepsilon^2)\), the same formulas hold with \(K_n\) replaced by \(K_n+\sigma_\varepsilon^2 I_n\). At each iteration of BO, we compute the posterior predictive mean $m_n(\mathbf x)$ and variance $s_n^2(\mathbf x)$ for any point $\mathbf{x}$, which will be used in the acquisition function to determine where to sample next.\\

In this work we use the Expected Improvement (EI) acquisition function. For minimisation, let
\[
    \tau_n \coloneq \min_{1\le i\le n} y_i
\]
denote the incumbent objective value after \(n\) evaluations. The improvement obtained by evaluating \(\mathbf x\) is
\[
    I_n(\mathbf x)\coloneq(\tau_n-F(\mathbf x))_+,
\]
where \((a)_+ \coloneq \max\{a,0\}\). Therefore, conditioning on $\mathcal{D}_n$,
\[
    u^{\operatorname{EI}}(\mathbf x)
    \coloneq
    \mathbb E[I_n(\mathbf x)\mid\mathcal D_n].
\]
Since \(F(\mathbf x)\mid\mathcal D_n\sim \mathcal N(m_n(\mathbf x),s_n^2(\mathbf x))\), we obtain
\[
    u^{\operatorname{EI}}(\mathbf x)
    =
    (\tau_n-m_n(\mathbf x))\Phi(\gamma_n(\mathbf x))
    +
    s_n(\mathbf x)\varphi(\gamma_n(\mathbf x)),
    \qquad
    \gamma_n(\mathbf x)
    \coloneq
    \frac{\tau_n-m_n(\mathbf x)}{s_n(\mathbf x)},
\]
with the convention \(u^{\operatorname{EI}}(\mathbf x) \coloneq 0\) when \( s_n(\mathbf{x}) = 0\) \cite{Mockus1978}. The next evaluation point is chosen by approximately maximising \(\operatorname{EI}_n\) over the current search region.\\

To accelerate BO while avoiding premature shrinkage, we propose to implement SDR \cite{SDR} into the standard BO framework so that the search region can be refined to locate the global minimiser more efficiently according to the incumbent solutions found so far. Compared to the traditional SDR implementation that updates the search region at each iteration, we update the search region only every \(K\) evaluations based on the current incumbent(s), rather than at every iteration. This allows sufficient exploration within the present bounds before contraction, reducing the risk of excluding the global minimiser.

\subsection{Variational Autoencoders}\label{subsection VAE}
Dimensionality reduction methods reduce the number of variables in a dataset while preserving essential information \cite{Velliangiri2019}. They can often be framed as an \textit{Encoder-Decoder} process, where the \textit{encoder} maps high-dimensional data to a lower-dimensional latent space, and the \textit{decoder} reconstructs the original data. Unlike autoencoders, which focus solely on minimising reconstruction error, VAEs optimise both reconstruction and latent space regularisation, facilitating more meaningful exploration of the latent space. Consequently, VAEs are useful for LSBO because they can generate novel candidate points while maintaining a structured latent representation, allowing for more efficient and reliable BO. For this reason, we focus on VAEs \cite{kingma2022autoencodingvariationalbayes, doersch2021tutorialvariationalautoencoders}, which provide a probabilistic encoder–decoder framework based on variational inference \cite{Hinton1993, Jordan1998}. VAEs utilise neural networks as encoders and decoders to generate latent manifolds. Let \(Z \coloneq \mathbb R^d\), with \(d < D\), denote the latent space. For a training input \(\mathbf x\in\mathcal X\subseteq\mathbb R^D\), the encoder network \(E_\phi\) outputs two \(d\)-dimensional vectors,
\[
    E_\phi(\mathbf x)
    =
    \bigl(\mathbf m_\phi(\mathbf x),\mathbf s_\phi(\mathbf x)\bigr),
    \qquad \mathbf s_\phi(\mathbf x)\in\mathbb R^d_{>0}.
\]
These outputs parametrise the approximate posterior distribution
\[
    q_\phi(\mathbf z\mid \mathbf x)
    \coloneq
    \mathcal N\!\left(
        \mathbf z;\,
        \mathbf m_\phi(\mathbf x),
        \operatorname{diag}(\mathbf s_\phi^2(\mathbf x))
    \right)
\]
over \(Z\). Thus \(q_\phi(\cdot\mid\mathbf x)\) is a probability distribution on \(Z\). Similarly, the decoder network \(\mathbf r_\theta:Z\to\mathbb R^D\) parametrises the reconstruction distribution \(p_\theta(\mathbf x\mid\mathbf z)\), such as 
\[
    p_\theta(\mathbf x\mid\mathbf z)
    =
    \mathcal N\!\left(
        \mathbf x;\mathbf r_\theta(\mathbf z),\sigma_x^2 I_D
    \right).
\]
The vector \(\mathbf r_\theta(\mathbf z)\) is the decoder mean, or reconstruction, of the latent point \(\mathbf z\). The VAE training objective for a single data point \(\mathbf x\) is the Evidence Lower BOund (ELBO)
\begin{equation}\label{VAE ELBO raw single point}
    \begin{split}
       \mathcal{J}_{\mathrm{VAE}}(\theta, \phi; \mathbf{x}) & \coloneq \ln p_{\theta}(\mathbf{x}) - D_{KL} [q_{\phi}(\mathbf{z|x}) \| p_{\theta}(\mathbf{z|x})] \\
       & = \underbrace{\mathbb{E}_{q_{\phi}(\mathbf{z|x})} [\ln p_{\theta}(\mathbf{x|z})]}_{\mathcal{J}_{\mathrm{recon}}} - \underbrace{D_{KL} [q_{\phi}(\mathbf{z|x}) \| p(\mathbf{z})]}_{\mathcal{J}_{KL}},
    \end{split}
\end{equation}
where $\ln p_{\theta}(\mathbf{x})$ is the marginal log-likelihood, \(p(\mathbf z)=\mathcal N(\mathbf 0,I_d)\) is the latent prior, and \(D_{\mathrm{KL}}\) denotes the Kullback--Leibler divergence between the true and the approximate posteriors. For a training set \(\mathcal D\), we maximise the average objective
\begin{equation}\label{VAE ELBO raw}
    \mathcal L(\theta,\phi;\mathcal D)
    \coloneq
    \frac{1}{|\mathcal D|}
    \sum_{\mathbf x\in\mathcal D}
    \mathcal J_{\mathrm{VAE}}(\theta,\phi;\mathbf x).
\end{equation}
Equivalently, one may minimise the negative ELBO loss \(\mathcal L_{\mathrm{VAE}}(\theta,\phi;\mathbf x)\coloneq -\mathcal J_{\mathrm{VAE}}(\theta,\phi;\mathbf x)\).

To make the optimisation of the ELBO \eqref{VAE ELBO raw single point} tractable, the prior $p(\mathbf{z})$ and posterior $q_{\phi}(\mathbf{z|x})$ are assumed to be parametrised as Gaussians with diagonal covariance matrices. In particular, $p(\mathbf{z})$ is commonly set to be the standard Gaussian $\mathcal{N}(\boldsymbol{0}, \mathbf{I})$ and thus the posterior $q_{\phi}(\mathbf{z|x})$ is Gaussian with mean $\mathbf m_\phi(\mathbf x) = \left[\mathbf{m}_{\phi,1}(\mathbf{x}), \ldots, \mathbf{m}_{\phi, d}(\mathbf{x})\right]$ and the covariance $\operatorname{diag}(\mathbf s_\phi^2(\mathbf x)) = \operatorname{diag}(s_{\phi,1}^2(\mathbf x), \ldots, s_{\phi,d}^2(\mathbf x))$. As shown in \cite{kingma2022autoencodingvariationalbayes}, under the Gaussian choices above, the KL term has the closed form
\[
    D_{\mathrm{KL}}
    \bigl(q_\phi(\mathbf z\mid\mathbf x)\,\|\,p(\mathbf z)\bigr)
    =
    \frac12
    \sum_{j=1}^d
    \left(
        s_{\phi,j}^2(\mathbf x)
        +
        m_{\phi,j}^2(\mathbf x)
        -
        1
        -
        \log s_{\phi,j}^2(\mathbf x)
    \right).
\]
With \(p_\theta(\mathbf x\mid\mathbf z) = \mathcal N(\mathbf x;\mathbf r_\theta(\mathbf z),\sigma_x^2 I_D)\), up to an additive constant independent of \(\theta\) and \(\phi\),
\begin{equation}\label{ELBO}
    \begin{split}
        \mathcal J_{\mathrm{VAE}}(\theta,\phi;\mathbf x)
    =
    -
    \frac{1}{2\sigma_x^2}
    \mathbb E_{q_\phi(\mathbf z\mid\mathbf x)}
    \left[
        \|\mathbf x-\mathbf r_\theta(\mathbf z)\|_2^2
    \right]
    -
    D_{\mathrm{KL}}
    \bigl(q_\phi(\mathbf z\mid\mathbf x)\,\|\,p(\mathbf z)\bigr).
    \end{split}
\end{equation}
These assumptions allow the use of the reparameterisation trick \cite{kingma2022autoencodingvariationalbayes} which enables gradient-based optimisation via Adam \cite{kingma2017adammethodstochasticoptimization}. Sampling from \(q_\phi(\mathbf z\mid\mathbf x)\) is implemented through the
reparameterisation
\[
    \mathbf z
    =
    \mathbf m_\phi(\mathbf x)
    +
    \mathbf s_\phi(\mathbf x)\odot\boldsymbol\varepsilon,
    \qquad
    \boldsymbol\varepsilon\sim\mathcal N(\mathbf 0,I_d),
\]
where \(\odot\) denotes componentwise multiplication. Additionally, as a modification of the classical VAE model, an additional weight hyperparameter $\beta$ is introduced in \eqref{VAE ELBO raw single point} to scale the $\mathcal{J}_{KL}$ term to trade off between reconstruction accuracy and the latent space regularity. It helps mitigate $\mathcal{J}_{KL}$ collapse, in which the latent variables carry little useful information \cite{bowman2016generatingsentencescontinuousspace, burgess2018understandingdisentanglingbetavae}. A practical implementation is to initialise $\beta$ at $0$ and gradually increase it in uniform increments over equal intervals until $\beta$ reaches $1$.

\subsection{Deep Metric Loss}

As described in \cite{grosnit2021highdimensionalbayesianoptimisationvariational, Ishfaq2018}, deep metric learning can be incorporated into VAEs by adding a metric-learning term to the VAE training objective. In this work, we focus on the standard triplet loss \cite{Hoffer2015}, although other metric-learning losses could also be used, see \cite{grosnit2021highdimensionalbayesianoptimisationvariational}. For a triplet
\[
    \mathbf x_{ijk}\coloneq\langle \mathbf x_i,\mathbf x_j,\mathbf x_k\rangle,
\]
the point \(\mathbf x_i\) is the anchor, \( \mathbf x_j\) is a positive point, and \(\mathbf x_k\) is a negative point. Let
\[
    \mathbf z_i\sim q_\phi(\cdot\mid \mathbf x_i), \qquad
    \mathbf z_j\sim q_\phi(\cdot\mid \mathbf x_j), \qquad
    \mathbf z_k\sim q_\phi(\cdot\mid \mathbf x_k),
\]
where $q_\phi(\cdot|\mathbf x)$ is the latent distribution produced by the VAE encoder for input data $\mathbf x$, and write \(\mathbf z_{ijk}\coloneq\langle \mathbf z_i,\mathbf z_j,\mathbf z_k\rangle\). Throughout this section, \(\|\cdot\|\) denotes the Euclidean $\ell_2$ norm. The classical hard triplet loss aims to place the anchor closer to the positive point than to the negative point by a margin \(\rho>0\):
\[
    \mathcal L_{h\text{-trip}}(\mathbf z_{ijk})
    \coloneq
    \max\{0,\|\mathbf z_i- \mathbf z_j\|+\rho-\| \mathbf z_i- \mathbf z_k\|\}.
\]
Minimising $\mathcal{L}_{h\text{-trip}} (\cdot)$ yields a structured embedding space with the positive and negative pairs being separated by a margin $\rho$. In black-box optimisation, class labels are not available. Following \cite{grosnit2021highdimensionalbayesianoptimisationvariational}, we instead construct positives and negatives using objective values by introducing a parameter $\eta$ to quantify differences between objective values. Let
\[
    \mathcal D_L\coloneq\{(\mathbf x_i,f_i)\}_{i=1}^N, \qquad f_i \coloneq f(\mathbf x_i),
\]
be a labelled dataset. For triplet construction, the objective values are rescaled to \([0,1]\). We keep the notation \(f_i\) for these rescaled values in this subsection. Given a threshold \(\eta\in(0,1)\), define
\[
    \mathcal D_p(\mathbf x_i;\eta)
    \coloneq
    \{\mathbf x_j : |f_i-f_j|<\eta\},
\qquad 
    \mathcal D_n(\mathbf x_i;\eta)
    \coloneq
    \{\mathbf x_k : |f_i-f_k|\geq \eta\}.
\]
Thus, \(\mathcal D_p( \mathbf x_i;\eta)\) contains points whose objective values are close to that of the anchor, while \(\mathcal D_n(\mathbf x_i;\eta)\) contains points whose objective values are sufficiently different. We define the valid triplet index set
\[
    \mathcal T(\mathcal D_L;\eta)
    \coloneq
    \{(i,j,k): \mathbf x_j\in \mathcal D_p(\mathbf x_i;\eta),\ \mathbf x_k\in \mathcal D_n(\mathbf x_i;\eta)\}.
\]

\begin{figure}[h]
    \centering
    \includegraphics[width=\linewidth]{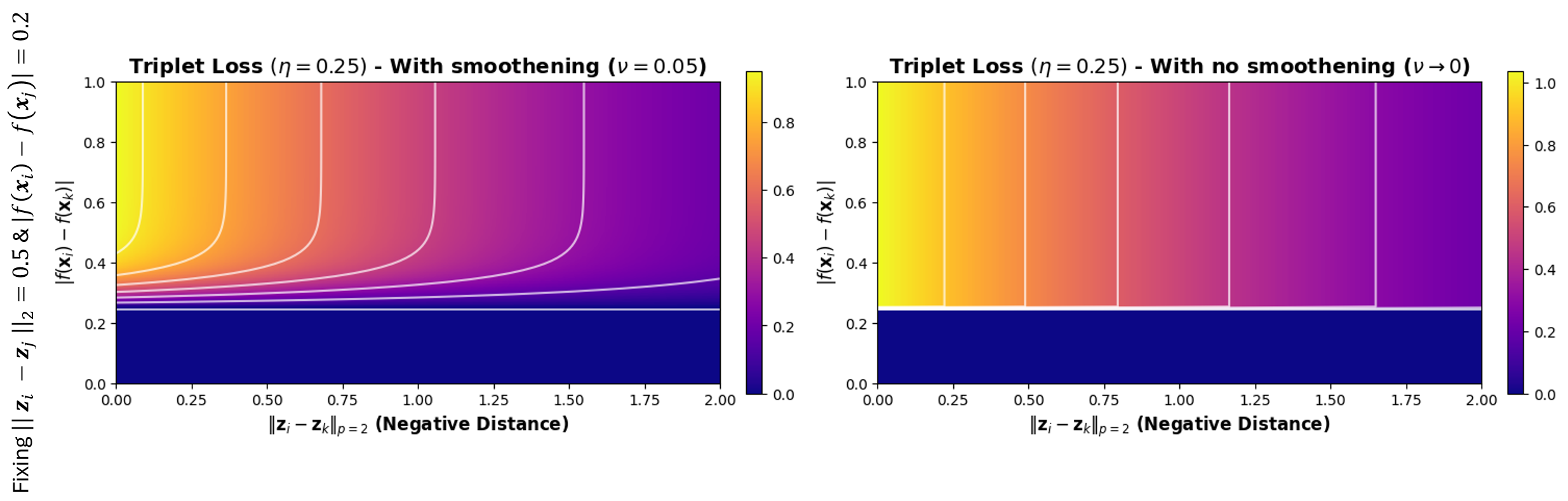}
    \caption{\textbf{Soft Triplet Loss}. Given a latent triplet $\mathbf{z}_{ijk} = \langle \mathbf{z}_i, \mathbf{z}_j, \mathbf{z}_k \rangle$ with the anchor $\mathbf{z}_i $ corresponding to the original anchor point $\mathbf{x}_i$, the right plot illustrates the discontinuous behaviour of $\mathcal{L}_{s-trip}(\cdot)$ when no weights are used. The discontinuity arises at the plane \( |f(\mathbf{x}_i) - f(\mathbf{x}_k)| = \eta = 0.25 \), marking the threshold beyond which \( \mathbf{x}_k \) no longer belongs to the set of negative data points relative to the anchor \( \mathbf{x}_i \). The left plot illustrates a smooth transition when approaching the plane.}
    \label{fig: triplet loss visualisation}
\end{figure}

For a valid triplet \((i,j,k)\in\mathcal T(\mathcal D_L;\eta)\), let
\[
    d_z^+ \coloneq \|\mathbf z_i- \mathbf z_j\|,
    \qquad
    d_z^- \coloneq \|\mathbf z_i- \mathbf z_k\|.
\]
To obtain a differentiable metric-learning objective, we use the soft triplet loss
\begin{equation}\label{L_s_trip}
        \mathcal L_{s\text{-trip}}(\mathbf z_{ijk})
    \coloneq
    \ln\!\left(1+\exp(d_z^+-d_z^-)\right)\omega_{ij}\omega_{ik},
\end{equation}
where
\[
    \omega_{ij}
    \coloneq
    \frac{\psi_\nu(\eta-|f_i-f_j|)}{\psi_\nu(\eta)},
    \qquad
    \omega_{ik}
    \coloneq
    \frac{\psi_\nu(|f_i-f_k|-\eta)}{\psi_\nu(1-\eta)}.
\]
Here
\[
    \psi_\nu(t)\coloneq\tanh\!\left(\frac{t}{2\nu}\right),
    \qquad t\geq 0,\quad \nu>0,
\]
is a smoothing function with $\lim_{\nu \rightarrow 0} \psi_\nu(t) = 1$. The parameter \(\nu\) controls the sharpness of the transition near the threshold \(\eta\), such that smaller values give a sharper transition. The weights $\omega_{ij}$ and $\omega_{ik}$ modulate the attraction and repulsion terms according to the corresponding objective-value differences.
The metric-learning term used in Algorithm \ref{BOVAE-with-DML} is obtained by aggregating the soft triplet losses over the current labelled dataset. Specifically, define
\begin{equation}\label{metric_loss_avg}
    \mathcal L_{\mathrm{metric}}(\phi;\mathcal D_L)
    \coloneq
    \frac{1}{|\mathcal T(\mathcal D_L;\eta)|}
    \sum_{(i,j,k)\in\mathcal T(\mathcal D_L;\eta)}
    \mathbb E_{q_\phi(\mathbf z_{ijk}\mid \mathbf x_{ijk})}
    \left[
        \mathcal L_{s\text{-trip}}(\mathbf z_{ijk})
    \right],
\end{equation}
where
\[
    q_\phi(\mathbf z_{ijk}\mid \mathbf x_{ijk})
    \coloneq
    q_\phi(\mathbf z_i\mid \mathbf x_i)q_\phi(\mathbf z_j\mid \mathbf x_j)q_\phi(\mathbf z_k\mid \mathbf x_k).
\]
If \(\mathcal T(\mathcal D_L;\eta)\) is empty, the metric term is omitted. Intuitively, $\mathcal{L}_{s\text{-trip}}(\cdot)$ encourages points with objective values similar to that of the anchor to lie close together, while pushing the points with dissimilar function values farther away. Consequently, in the latent space, points with similar function values are grouped into small clusters. The weights $\omega_{ij}$ and $\omega_{ik}$ smooth the transition accross $\mathcal{L}_{s-trip}(\cdot)$. Without them, discontinuities occur around the planes $|f(\mathbf{x}_i) - f(\mathbf{x}_j)| = \eta$ and $|f(\mathbf{x}_i) - f(\mathbf{x}_k)| = \eta$ if the weights are not used. Figure \ref{fig: triplet loss visualisation} illustrates this effect\footnote{To guarantee this figure’s reproducibility, we revised and re-implemented the triplet DML routine and a stand-alone script that reliably generates the visualisation. Further details are available in Section \ref{Algorithm Implementations} and our GitHub page.}.
The DML-augmented VAE objective is therefore
\begin{equation}\label{VAE_DML_augmented_objective}
    \mathcal L_{\mathrm{DML}}(\theta,\phi;\mathcal D_L)
    \coloneq
    \mathcal{L}(\theta,\phi;\mathcal D_L)- \lambda_{\mathrm{DML}}\mathcal L_{\mathrm{metric}}(\phi;\mathcal D_L),
\end{equation}
where $\lambda_{\mathrm{DML}} > 0$ indicates the weight for the metric loss term and $\mathcal L(\theta,\phi; \mathcal D_L)$ follows from \eqref{VAE ELBO raw}. Thus, during the DML retraining step in Algorithm \ref{BOVAE-with-DML}, the VAE is trained to retain reconstruction quality and latent regularity while also encouraging points with similar objective values to be close in latent space. This aggregation is computationally manageable in our BO setting because the triplet loss is formed only from the current labelled dataset \(\mathcal D_L\), which consists of the initial labelled points if used and the objective evaluations accumulated so far. This dataset is small compared with the unlabelled pre-training set \(\mathcal D_U\). In the implementations, the triplet term is evaluated over the valid triplets in the current training batch, which avoids constructing triplets from the full unlabelled dataset.

\section{VAE-Based Latent-Space BO: Algorithms and Fixed-Representation Analysis}\label{Alg}
As mentioned above, dimensionality reduction techniques help reduce the optimisation problem's dimensionality. Using a VAE within BO enables the surrogate optimisation to be performed in the lower-dimensional latent space $\mathcal{Z}$, rather than directly in the ambient space. In BO--VAE\footnote{For brevity, we use \textit{BO--VAE} for BO conducted in a VAE-induced latent space.}, a trained VAE decoder \(p_{\theta^\ast}(\cdot\mid \mathbf z)\) induces a latent objective
\[
    \bar f_{\theta^\ast}(\mathbf z)
    \coloneq
    \mathbb E_{\mathbf x\sim p_{\theta^\ast}(\cdot\mid \mathbf z)}
    [f(\mathbf x)],
    \qquad \mathbf z\in Z.
\]
The BO--VAE algorithms do not solve Problem \eqref{main problem} directly. Instead, we apply BO to the decoder-induced latent problem
\begin{equation}\label{latent problem}
    \bar f^\ast_{\theta^\ast}
    \coloneq
    \inf_{\mathbf z\in R_0}\bar f_{\theta^\ast}(\mathbf z),
    \tag{LP}
\end{equation}
where \(R_0\subset Z\) is the initial latent search region used by the algorithm. The latent problem \eqref{latent problem} is an approximation, since its quality depends on whether the decoder can generate points near high-quality solutions of Problem \eqref{main problem}. When \(f^\ast=\min_{\mathbf x\in \mathcal X} f(\mathbf x)\) exists and the decoder distribution is supported on \(\mathcal X\), we have
\[
    f^\ast \leq \bar f_{\theta^\ast}(z)
    \quad\text{for all } z\in R_0,
\]
and hence \(f^\ast \leq \bar f^\ast_{\theta^\ast}\). Thus equality is not guaranteed, since \eqref{latent problem} is close to \eqref{main problem} only when the decoder represents near-optimal regions of \(\mathcal X\) sufficiently well. In practice, \(\bar f_{\theta^\ast}\) is unknown and is queried only through decoded candidate points and evaluations of \(f\). The GP surrogate in the BO loop is therefore used to model the latent objective from the accumulated latent data. To understand BO--VAE, we provide two useful insights about the quality of the latent problem:

\begin{itemize}
    \item \emph{Coverage}: the data-generating distribution gives positive mass in at least one neighbourhood of the global minimiser \(x^{\ast}\) (so \(x^{\ast}\) is representable by the VAE);\\
    \item \emph{Regularity}: a well-trained VAE can provide a regular latent geometry, empirically encouraged by the $\mathcal{L}_{\mathrm{KL}}$ from the ELBO objective, and the upper-bounded Wasserstein–1 generation/regeneration gaps between the data-generated distribution and the VAE's generated/regenerated distributions (c.f. the bounds from Theorems $5.1-5.4$ from \cite{mbacke2023statisticalguaranteesvariationalautoencoders}).
\end{itemize}

The KL term in the ELBO encourages regularity of the latent distribution, while reconstruction training encourages nearby latent representations to decode to points on the learned data manifold, allowing BO to search efficiently in \(\mathcal{Z}\) while the decoder maps promising latent points back to feasible candidates in \(\mathcal{X}\). However, we note that neither the KL term in the VAE objective nor a global distributional reconstruction alone guarantees that the decoder can recover the ambient optimiser. A theoretical analysis of a simplified fixed-representation setting is provided in Section \ref{sec:fixed-representation-theory}. There, the VAE and latent search region are held fixed, and the decoder mean is treated as the reconstructed point. This permits the optimisation error within the latent space to be separated from the approximation error introduced by the learned representation.\\

Below, we introduce our three BO-VAE algorithms. We first present the baseline BO-VAE algorithms, followed by the variants incorporating VAE retraining, and soft-triplet DML. For the standard BO algorithm with SDR, we include it as Algorithm~\ref{BOSDR} in Appendix~\ref{appendix: BO-SDR} to keep our focus on VAE-assisted BO frameworks. We distinguish two update schedules. The BO iteration index \(k\) counts objective evaluations. The SDR update period \(K_\mathrm{SDR}\) controls how often the region is updated, where \(K_\mathrm{SDR}=1\) corresponds to updating the region after every BO evaluation. Separately, the VAE retraining period \(q\), used only in Algorithms \ref{BOVAE-retraining} and \ref{BOVAE-with-DML}, controls how often the encoder--decoder is updated using the current labelled data. Algorithm \ref{BOVAE-SDR} trains the VAE once and keeps it fixed throughout the BO loop.

\subsection{BO--VAE with SDR and a Fixed Pretrained VAE}
Algorithm~\ref{BOVAE-SDR} is the baseline VAE-based latent-space BO method. The VAE is trained once before the BO loop and is then kept fixed, i.e., no VAE retraining and no deep metric learning are used in this algorithm. This baseline allows us to isolate the effect of performing BO and SDR in a fixed VAE-induced latent space before introducing the retraining variants in Algorithms~\ref{BOVAE-retraining} and~\ref{BOVAE-with-DML}. We distinguish three datasets. The unlabelled dataset
\[
    \mathcal D_U \coloneq \{\mathbf{x}_m^U\}_{m=1}^M \subset \mathcal X
\]
is used only for VAE pre-training and does not require objective evaluations. The initial labelled dataset is
\[
    \mathcal D_L^0
    \coloneq
    \{(\mathbf{x}_i^L,y_i)\}_{i=1}^N,
    \qquad
    y_i \coloneq f(\mathbf{x}_i^L),
\]
and contains the initial points at which \(f\) has been evaluated. After training the VAE, the labelled points are encoded to form the initial latent GP dataset
\[
    \mathcal D_Z^0
    \coloneq
    \{(\mathbf{z}_i^L,y_i)\}_{i=1}^N,
    \qquad
    \mathbf{z}_i^L
    \coloneq
    \mathbb E_{q_{\phi^\ast}(\cdot\mid \mathbf{x}_i^L)}[\mathbf z].
\]
In our implementation, VAE pre-training uses the \(\beta\)-annealed ELBO described in Section~\ref{subsection VAE}, which affects only line~1 of Algorithm~\ref{BOVAE-SDR} and does not change the BO loop.
\begin{algorithm}[h]
\caption{BO--VAE with SDR and a Fixed Pretrained VAE}
\label{BOVAE-SDR}
\begin{algorithmic}[1]
\small
\Require Unlabelled dataset $\mathcal D_U \coloneq \{\mathbf{x}_m^U\}_{m=1}^M$; initial labelled dataset $\mathcal D_L^0 \coloneq \{(\mathbf{x}_i^L,y_i)\}_{i=1}^N$, where $y_i \coloneq f(\mathbf{x}_i^L)$; budget $B$; initial latent region $\mathcal R_0\subseteq\mathcal Z$; SDR update period $K_{\mathrm{SDR}}$; acquisition rule $u(\cdot|\cdot)$ \textrm{(EI)}; encoder $q_\phi(\mathbf z\mid\mathbf x)$; decoder $p_\theta(\mathbf x\mid\mathbf z)$.

\Ensure Best objective value $f_{\min}$ found.

\State Train VAE: $(\theta^\ast,\phi^\ast) \leftarrow \operatorname*{arg\,max}_{\theta,\phi} \mathcal J_{\mathrm{VAE}}(\theta,\phi;\mathcal D_U)$.

\State Encode the initial labelled data: $\mathcal D_Z^0 \leftarrow \left\{(\mathbf z_i^L,y_i): \mathbf z_i^L \coloneq \mathbb E_{q_{\phi^\ast}(\cdot\mid\mathbf x_i^L)}[\mathbf z], \;(\mathbf x_i^L,y_i)\in\mathcal D_L^0 \right\}$.

\For{$k=0,1,\ldots,B-1$}
    \State Fit GP surrogate $h_k:\mathcal Z\to\mathbb R$ on $\mathcal D_Z^k$.
    \State $u_k(\mathbf z)\coloneq u(\mathbf z\mid\mathcal D_Z^k)$. \Comment{The current acquisition function.}
    \State $\hat{\mathbf z}_{k+1} \leftarrow \operatorname*{arg\,max}_{\mathbf z\in\mathcal R_k} u_k(\mathbf z)$.
    \State $\hat{\mathbf x}_{k+1} \sim p_{\theta^\ast}(\cdot\mid\hat{\mathbf z}_{k+1})$. \Comment{Decoding.}
    \State $y_{k+1}\leftarrow f(\hat{\mathbf x}_{k+1})$. \Comment{Evaluation of $f.$}
    \State $\mathcal D_L^{k+1} \leftarrow \mathcal D_L^k \cup \{(\hat{\mathbf x}_{k+1},y_{k+1})\}$. 
    \State $\mathcal D_Z^{k+1} \leftarrow \mathcal D_Z^k \cup \{(\hat{\mathbf z}_{k+1},y_{k+1})\}$. 
    \If{$(k+1)\bmod K_{\mathrm{SDR}}=0$} \Comment{Update search region.}
        \State $\mathcal R_{k+1} \leftarrow \mathrm{SDR\_Update}(\mathcal R_k,\mathcal D_Z^{k+1})$.
    \Else
        \State $\mathcal R_{k+1}\leftarrow \mathcal R_k$.
    \EndIf
\EndFor

\State \Return $f_{\min} \leftarrow \min\{y:(\mathbf x,y)\in\mathcal D_L^B\}$.
\end{algorithmic}
\end{algorithm}
\FloatBarrier
\normalsize
At BO iteration \(k\), a GP surrogate is fitted to \(\mathcal D_Z^k\), and the corresponding acquisition function is maximised over the current latent region \(\mathcal R_k\). The selected latent point is decoded to an ambient point, where \(f\) is evaluated. The new ambient evaluation is added to \(\mathcal D_L^{k+1}\), while the pair consisting of the latent query and its objective value is added to \(\mathcal D_Z^{k+1}\). The SDR update then recentres and rescales the latent region of interest using the current incumbent and the observed incumbent movement. The detailed panning--zooming update is given in Appendix~\ref{appendix: BO-SDR}.

\subsection{A Fixed-Representation Convergence Analysis}
\label{sec:fixed-representation-theory}

We now isolate the two sources of error in BO--VAE, that is, the optimisation error within the latent problem and the approximation error caused by the learned representation.  The analysis establishes a guarantee for an idealised fixed-pretrained-VAE variant in which the SDR is not enabled, the decoder mean is evaluated deterministically, the GP hyperparameters are fixed, and the EI acquisition function is globally optimised exactly with the fixed prior used in the convergence theorem of \cite{bull2011ego}.  Thus, set $\mathcal{R}_k=\mathcal{R}_0$ and $\mathbf{x}_i=r_{\theta^\star}(\mathbf{z}_i)$ throughout this subsection.  We note that these are idealised changes to the implementations of Algorithm \ref{BOVAE-SDR} for theoretical analysis. 

Write $\mathcal{Z} \coloneq \mathbb{R}^d$.  Define the mean-decoded latent objective and its optimum by
\begin{equation}
  g(\mathbf{z}) \coloneq f\bigl(r_{\theta^\star}(\mathbf{z})\bigr),
  \qquad
  g^\star \coloneq \min_{\mathbf{z}\in\mathcal{R}_0}g(\mathbf{z}),
  \label{eq:fixed-latent-objective}
\end{equation}
and recall $f_{\mathcal{X}}^\star \coloneq \min_{\mathbf{x}\in\mathcal{X}}f(\mathbf{x})$ and $\Delta_{\mathcal{X}} \coloneq f_{\mathcal{X}}^\star-f_\Omega^\star = 0$. Since $r_{\theta^\star}(\mathcal{R}_0)\subseteq\mathcal{X}$, the VAE-induced representation gap is
\begin{equation}
  \Delta_{\mathrm{rep}}
   \coloneq g^\star-f_{\mathcal{X}}^\star\geq 0.
  \label{eq:fixed-representation-gap}
\end{equation}
For any latent evaluations $\mathbf{z}_1,\ldots,\mathbf{z}_n$, the ambient simple regret has the exact decomposition
\begin{equation}
  \min_{1\leq i\leq n}g(\mathbf{z}_i)-f_\Omega^\star
  =
  \underbrace{\min_{1\leq i\leq n}g(\mathbf{z}_i)-g^\star}_{\text{latent optimization error}}
  +\Delta_{\mathrm{rep}}.
  \label{eq:fixed-regret-decomposition}
\end{equation}

\begin{assumption}[Reconstruction certificate]
\label{ass:fixed-representation-certificate}

The following conditions hold. The first two concern the ambient problem and the fixed VAE, while the remaining conditions concern the reconstruction certificate.

\begin{enumerate}
  \item $\mathcal{X}\subseteq\Omega\subseteq\mathbb{R}^D$ is nonempty and compact, $f_\Omega^\star> -\infty$,   $f|_{\mathcal{X}}$ is $L_f$-Lipschitz, and $\mathcal{S}_{\mathcal{X}}^\star  \coloneq \arg\min_{\mathbf{x}\in\mathcal{X}}f(\mathbf{x})$ is nonempty.
  
  \item The fixed decoder mean $r_{\theta^\star}:\mathcal{Z}\to\mathcal{X}$ is $K_\theta$-Lipschitz.  The fixed encoder is
  \[
    q_{\phi^\star}(\,\cdot\mid\mathbf{x})
    =\mathcal{N}\!\left(\mathbf{m}_{\phi^\star}(\mathbf{x}),
      \operatorname{diag}\bigl(\mathbf{s}_{\phi^\star}^2(\mathbf{x})\bigr)\right),
  \]
  with strictly positive standard deviations, and the parameter map
  $\mathbf{x}\mapsto
  (\mathbf{m}_{\phi^\star}(\mathbf{x}),
   \mathbf{s}_{\phi^\star}(\mathbf{x}))$ 
   is $K_\phi$-Lipschitz.  The PAC--Bayes prior is $p=\mathcal{N}(\mathbf{0},\mathbf{I}_d)$.

  \item The population measure $\mu$ is supported on a nonempty compact set $\mathcal{M}\subseteq\mathcal{X}$.  For some $c_\mu>0$, $m>0$, and
  $\rho_0>0$,
  \begin{equation}
    \mu\bigl(B(\mathbf{x},\rho)\cap\mathcal{M}\bigr)
    \geq c_\mu\rho^m,
    \qquad
    \mathbf{x}\in\mathcal{M},\quad 0<\rho\leq\rho_0.
    \label{eq:local-mass-condition}
  \end{equation}
  Write $\delta_{\mathcal{M}}  \coloneq \operatorname{dist}(\mathcal{M},\mathcal{S}_{\mathcal{X}}^\star)$.

  \item The latent region $\mathcal{R}_0\subset\mathcal{Z}$ is compact with nonempty interior and
  \begin{equation}
    \beta_0 \coloneq \inf_{\mathbf{x}\in\mathcal{M}}
      q_{\phi^\star}(\mathcal{R}_0\mid\mathbf{x})>0.
    \label{eq:latent-region-mass}
  \end{equation}

  \item An independent certificate sample
  $\mathcal{D}_{\mathrm{cert}}
  =(\mathbf{X}_1,\ldots,\mathbf{X}_N)\sim\mu^{\otimes N}$ 
  is drawn after the VAE has been fixed.  The confidence level $\eta\in(0,1)$ and $\lambda>0$ are chosen independently of this sample.
  
\end{enumerate}
\end{assumption}

Let $D_{\mathcal{X}} \coloneq \operatorname{diam}(\mathcal{X})$ and define the empirical reconstruction term
\begin{equation}
  \widehat{\mathcal{R}}_N
   \coloneq \frac{1}{N}\sum_{i=1}^N
    \mathbb{E}_{\mathbf{Z}\sim q_{\phi^\star}(\cdot\mid\mathbf{X}_i)}
      \bigl[\|\mathbf{X}_i-r_{\theta^\star}(\mathbf{Z})\|\bigr].
  \label{eq:empirical-reconstruction-certificate}
\end{equation}
The associated PAC--Bayes quantity is
\begin{align}
  \varepsilon_N^{\mathrm{PB}}
   \coloneq \widehat{\mathcal{R}}_N
  +\frac{1}{\lambda}\biggl[
    &\sum_{i=1}^N
      D_{\mathrm{KL}}\!\left(q_{\phi^\star}(\cdot\mid\mathbf{X}_i)\|p\right)
      +\lambda K_\phi K_\theta D_{\mathcal{X}}
      +\log\frac{1}{\eta}
      +\frac{\lambda^2D_{\mathcal{X}}^2}{8N}
  \biggr],
  \label{eq:pac-bayes-certificate}
\end{align}
and put $\overline\varepsilon_N^{\mathrm{PB}} \coloneq \min\{D_{\mathcal{X}},\varepsilon_N^{\mathrm{PB}}\}$.  Now, set
\begin{equation}
  \varepsilon_N^{\mathrm{cert}}
   \coloneq \min\left\{
    D_{\mathcal{X}},
    \frac{1}{\beta_0}
    \inf_{0<\rho\leq\rho_0}
    \left[
      (1+K_\phi K_\theta)\rho
      +\frac{\overline\varepsilon_N^{\mathrm{PB}}}
             {c_\mu\rho^m}
    \right]
  \right\}.
  \label{eq:uniform-reconstruction-certificate}
\end{equation}
The lower-mass condition in \eqref{eq:local-mass-condition} is what permits an average reconstruction certificate to control reconstruction near an unknown optimizer. Without such a condition, an arbitrarily poorly reconstructed region could have negligible $\mu$-mass.

\begin{assumption}[The prior specification required by \cite{bull2011ego}]
\label{ass:fixed-prior-ei}
The function $g$ in \eqref{eq:fixed-latent-objective} belongs to the RKHS $\mathcal{H}_{K_{5/2,\ell}}(\mathcal{R}_0)$ associated with a Mat\'ern-$5/2$ correlation kernel with fixed positive length-scales $\ell$, and $\|g\|_{\mathcal{H}_{K_{5/2,\ell}}(\mathcal{R}_0)}\leq B_g<\infty$.  All observations are the exact noiseless values $g(\mathbf{z}_i)$. The GP scale and length-scales are fixed throughout the EI run. A fixed number of distinct initial points is selected by a rule independent of $g$, and every subsequent $\mathbf{z}_i$ is a global maximizer of the exact EI acquisition over $\mathcal{R}_0$.
\end{assumption}

We now present the main theorem. We note that this is a minimal fixed-representation result, not a guarantee for the implemented SDR and the retraining algorithms.

\begin{theorem}[Fixed-representation BO--VAE]
\label{thm:fixed-representation-bovae}
Under the Assumptions~\ref{ass:fixed-representation-certificate} and \ref{ass:fixed-prior-ei}, with probability at least $1-\eta$ over $\mathcal{D}_{\mathrm{cert}}$, there exist constants $C_{\mathrm{EI}}<\infty$ and $n_0\in\mathbb{N}$, independent of the BO budget $n$, such that, for every $n\geq n_0$,
\begin{align}
  \mathbb{E}_{\mathrm{EI}}\!\left[
    \min_{1\leq i\leq n}f\bigl(r_{\theta^\star}(\mathbf{z}_i)\bigr)
    -f_\Omega^\star
  \right]
  &\leq C_{\mathrm{EI}}n^{-1/d}
       +\Delta_{\mathrm{rep}}
  \notag\\
  &\leq C_{\mathrm{EI}}n^{-1/d}
       +L_f\bigl(\delta_{\mathcal{M}}
                  +\varepsilon_N^{\mathrm{cert}}\bigr).
  \label{eq:fixed-representation-main-bound}
\end{align}
Here $\mathbb{E}_{\mathrm{EI}}$ is only over randomness in the initial design or the EI acquisition function, while the objective, VAE, and certificate sample are held fixed.
\end{theorem}

The proof is given in Appendix~\ref{app:fixed-representation-proofs}. The first line of \eqref{eq:fixed-representation-main-bound} is the exact
optimisation--representation separation.  The second line is an explicit, data-dependent high-probability upper certificate for the otherwise unknown $\Delta_{\mathrm{rep}}$.  It combines the reconstruction theorem of \cite{mbacke2023statisticalguaranteesvariationalautoencoders} with the fixed-prior EI rate of \cite{bull2011ego}. 

\begin{corollary}[Fixed-representation asymptote]
\label{cor:fixed-representation-floor}
On the event of Theorem~\ref{thm:fixed-representation-bovae},
\begin{equation}
  \lim_{n\to\infty}
  \mathbb{E}_{\mathrm{EI}}\!\left[
    \min_{1\leq i\leq n}f\bigl(r_{\theta^\star}(\mathbf{z}_i)\bigr)
    -f_\Omega^\star
  \right]
  =\Delta_{\mathrm{rep}}.
  \label{eq:fixed-representation-asymptote}
\end{equation}
Consequently, increasing only the BO budget removes the latent optimisation error, but it cannot remove a nonzero representation gap.
\end{corollary}

\begin{corollary}[Exact representation]
\label{cor:fixed-exact-representation}
If $\Delta_{\mathcal{X}}=0$ and $\Delta_{\mathrm{rep}}=0$, then the expected ambient simple regret is $O(n^{-1/d})$.  A sufficient, but not necessary, condition for $\Delta_{\mathrm{rep}}=0$ in the certificate is $\delta_{\mathcal{M}}=\varepsilon_N^{\mathrm{cert}}=0$.
\end{corollary}

The rate depends on the latent dimension $d$, rather than the ambient dimension $D$, while the representation floor can depend on both the VAE and
the ambient problem.  For one fixed VAE and one fixed certificate sample, $\varepsilon_N^{\mathrm{cert}}$ does not change with the BO budget $n$.

Finally, we note that adaptive SDR can remove a latent minimiser, decoder sampling changes the noiseless objective, numerical acquisition optimisation need not find the global EI maximizer, and VAE retraining makes both $g$ and its coordinates stage-dependent.  Extending the
guarantee to any of these mechanisms requires additional assumptions.  The fixed result nevertheless highlights that successful retraining should reduce the representation gap without introducing a larger optimisation error. An extension of the analysis to BO–VAE with retraining is currently under development.

\subsection{BO--VAE with SDR and Retraining}

Retraining under the latent-space optimisation framework is used to mitigate common failure modes of latent-space optimisation \cite{Tripp2020}. In the context of LSBO, it helps propagate new information associated with newly evaluated points from the BO routine into the VAE model. Without updates, the generative model remains static and may miss high-performing areas revealed during BO. Retraining incorporates newly evaluated points, adapting and extending the latent space toward promising regions by augmenting the standard generative model losses.

Algorithm~\ref{BOVAE-retraining} outlines our BO--VAE approach with SDR and periodic VAE retraining, as a variant of Algorithm~\ref{BOVAE-SDR}. The VAE is pre-trained once on the unlabelled dataset \(\mathcal D_U\) (line~\ref{alg:pretrain}). The algorithm then alternates between outer retraining stages and inner BO--SDR steps. At outer stage \(l\), the VAE is warm-started from the previous checkpoint and retrained on the current labelled dataset \(\mathcal D_L^{(l)}\) (line~\ref{alg:retrain}). The updated encoder--decoder \((\theta_l^\ast,\phi_l^\ast)\) is then used to re-encode \(\mathcal D_L^{(l)}\), forming the latent GP dataset \(\mathcal D_Z^{(l)}\) (line~\ref{alg:encode}), which seeds the subsequent inner BO routine. Between retraining steps, BO with SDR runs for at most \(q\) inner iterations in the latent space \(\mathcal Z\) (lines~\ref{alg:inner-start}--\ref{alg:inner-end}). At inner step \(k\), the algorithm fits a GP surrogate on \(\mathcal D_Z^{(l;k)}\), selects a latent query \(\hat{\mathbf z}_{l;k+1}\) by maximising EI over the current latent region \(\mathcal R_{(l;k)}\), decodes it to an ambient point \(\hat{\mathbf x}_{l;k+1}\), evaluates \(f\), and augments both the labelled ambient dataset \(\mathcal D_L^{(l;k+1)}\) and the latent GP dataset \(\mathcal D_Z^{(l;k+1)}\). The SDR region is updated according to the SDR update period \(K_{\mathrm{SDR}}\). After \(k\) inner steps, the datasets contain \(N+(l-1)q+(k+1)\) labelled evaluations. Thus, after \(q\) inner steps, they contain \(N+lq\) labelled evaluations and are passed to the next outer retraining stage.

\begin{algorithm}
\caption{BO--VAE with SDR and Periodic VAE Retraining}
\label{BOVAE-retraining}
\begin{algorithmic}[1]
\Require Initial labelled dataset
\(\mathcal D_L^{(1)}\coloneq \{(\mathbf x_i,y_i)\}_{i=1}^{N}\), where
\(y_i\coloneq f(\mathbf x_i)\);
unlabelled dataset
\(\mathcal D_U\coloneq\{\mathbf x_m^U\}_{m=1}^{M}\);
budget \(B\);
VAE retraining period \(q\);
SDR update period \(K_{\mathrm{SDR}}\);
initial latent region \(\mathcal R_0\subseteq\mathcal Z\);
EI acquisition \(u(\cdot\mid\cdot)\);
encoder \(q_\phi(\mathbf z\mid\mathbf x)\);
decoder \(p_\theta(\mathbf x\mid\mathbf z)\).
\Ensure Best objective value \(f_{\min}\) found.

\State
\((\theta_0^\ast,\phi_0^\ast)
\gets
\operatorname*{arg\,max}_{\theta,\phi}
\mathcal J_{\mathrm{VAE}}(\theta,\phi;\mathcal D_U)\).
\label{alg:pretrain}
\Comment{Pre-train VAE on \(\mathcal D_U\).}

\State
\((\theta_1^\ast,\phi_1^\ast)\gets(\theta_0^\ast,\phi_0^\ast)\);
\quad \(L\gets\lceil B/q\rceil\).

\For{\(l=1\) \textbf{to} \(L\)}
\label{alg:outer-start}

  \State
  \((\theta_l^\ast,\phi_l^\ast)
  \gets
  \operatorname*{arg\,max}_{\theta,\phi}
  \mathcal J_{\mathrm{VAE}}(\theta,\phi;\mathcal D_L^{(l)})\)
  \label{alg:retrain}
  \Comment{Retrain on current labelled set.}

  \State
  \(\mathcal D_Z^{(l)}
  \gets
  \left\{
  (\mathbf z_i,y_i):
  \mathbf z_i
  \coloneq
  \mathbb E_{q_{\phi_l^\ast}(\cdot\mid\mathbf x_i)}[\mathbf z],
  \;
  (\mathbf x_i,y_i)\in\mathcal D_L^{(l)}
  \right\}\)
  \label{alg:encode}

  \State
  \(\mathcal D_L^{(l;0)}\gets \mathcal D_L^{(l)}\);
  \quad
  \(\mathcal D_Z^{(l;0)}\gets \mathcal D_Z^{(l)}\);
  \quad
  \(\mathcal R_{(l;0)}\gets \mathcal R_0\).
  \Comment{Initialise inner loop and SDR region.}

  \State \(q_l\gets \min\{q,\,B-(l-1)q\}\).

  \For{\(k=0\) \textbf{to} \(q_l-1\)}
  \label{alg:inner-start}

    \State Fit GP surrogate \(h_{l;k}:\mathcal Z\to\mathbb R\) on
    \(\mathcal D_Z^{(l;k)}\).

    \State
    \(u_{l;k}(\mathbf z)
    \coloneq
    u(\mathbf z\mid\mathcal D_Z^{(l;k)})\).
    \Comment{Current EI acquisition.}

    \State
    \(\hat{\mathbf z}_{l;k+1}
    \gets
    \operatorname*{arg\,max}_{\mathbf z\in\mathcal R^{(l;k)}}
    u_{l;k}(\mathbf z)\).

    \State
    \(\hat{\mathbf x}_{l;k+1}
    \sim
    p_{\theta_l^\ast}(\cdot\mid\hat{\mathbf z}_{l;k+1})\);
    \quad
    \(y_{l;k+1}\gets f(\hat{\mathbf x}_{l;k+1})\).

    \State
    \(\mathcal D_L^{(l;k+1)}
    \gets
    \mathcal D_L^{(l;k)}
    \cup
    \{(\hat{\mathbf x}_{l;k+1},y_{l;k+1})\}\).

    \State
    \(\mathcal D_Z^{(l;k+1)}
    \gets
    \mathcal D_Z^{(l;k)}
    \cup
    \{(\hat{\mathbf z}_{l;k+1},y_{l;k+1})\}\).

    \State \(t_{l;k+1}\coloneq (l-1)q+k+1\).

    \If{\(t_{l;k+1}\bmod K_{\mathrm{SDR}}=0\)}
      \State
      \(\mathcal R_{(l;k+1)}
      \gets
      \mathrm{SDR\_Update}
      \bigl(\mathcal R_{(l;k)},\mathcal D_Z^{(l;k+1)}\bigr)\).
    \Else
      \State
      \(\mathcal R_{(l;k+1)}
      \gets
      \mathcal R_{(l;k)}\).
    \EndIf

  \EndFor
  \label{alg:inner-end}

  \State
  \(\mathcal D_L^{(l+1)}
  \gets
  \mathcal D_L^{(l;q_l)}\);
  \quad
  \(\mathcal D_Z^{(l+1)}
  \gets
  \mathcal D_Z^{(l;q_l)}\).

\EndFor
\label{alg:outer-end}

\State
\Return
\(f_{\min}
\gets
\min\{\,y:(\mathbf x,y)\in\mathcal D_L^{(L+1)}\,\}\).

\end{algorithmic}
\end{algorithm}
The retraining cost is controlled by the retraining period \(q\), warm-starting from the previous VAE checkpoint, and the size of the current labelled dataset \(\mathcal D_L^{(l)}\). In practice, \(\mathcal D_Z^{(l+1)}\) is recomputed after the next VAE retraining step using the freshly updated encoder, so carrying the latent dataset forward is only a bookkeeping device. To quantify this cost, Section \ref{Results} reports optimisation performance against total computation time.

\subsection{BO--VAE with Periodic DML Retraining}
We follow \cite{grosnit2021highdimensionalbayesianoptimisationvariational} and use DML to generate well-structured VAE-generated latent spaces. Algorithm~\ref{BOVAE-with-DML} is the DML-based retraining variant of BO--VAE. As in Algorithm~\ref{BOVAE-retraining}, Algorithm~\ref{BOVAE-with-DML} first pre-trains a VAE on the unlabelled dataset \(\mathcal D_U\) using the standard VAE objective. The DML term is not used during this pre-training stage. This separation is useful in the BO setting, where labelled objective evaluations are scarce early in the optimisation process.

\begin{algorithm}[h]
\caption{BO--VAE with SDR and Periodic VAE Retraining with DML}
\label{BOVAE-with-DML}
\begin{algorithmic}[1]
\small
\Require Initial labelled dataset \(\mathcal D_L^{(1)}\coloneq \{(\mathbf x_i,y_i)\}_{i=1}^{N}\), where \(y_i\coloneq f(\mathbf x_i)\); unlabelled dataset \(\mathcal D_U\coloneq\{\mathbf x_m^U\}_{m=1}^{M}\); budget \(B\); VAE retraining period \(q\); SDR update period \(K_{\mathrm{SDR}}\); initial latent region \(\mathcal R_0\subseteq\mathcal Z\); EI acquisition \(u(\cdot\mid\cdot)\); encoder \(q_\phi(\mathbf z\mid\mathbf x)\); decoder \(p_\theta(\mathbf x\mid\mathbf z)\).

\Ensure Best objective value \(f_{\min}\) found.

\State \((\theta_0^\ast,\phi_0^\ast)
\gets
\operatorname*{arg\,max}_{\theta,\phi}
\mathcal J_{\mathrm{VAE}}(\theta,\phi;\mathcal D_U)\).
\label{alg:dml-pretrain}
\Comment{Pre-train VAE on $\mathcal D_U$.}

\State $(\theta_1^\ast,\phi_1^\ast) \gets (\theta_0^\ast,\phi_0^\ast);\quad L\gets\lceil B/q\rceil$.

\For{$l=1$ \textbf{to} $L$}
\label{alg:dml-outer-start}

  \State $(\theta_l^\ast,\phi_l^\ast) \gets \operatorname*{arg\,max}_{\theta,\phi} \mathcal L_{\mathrm{DML}}(\theta,\phi;\mathcal D_L^{(l)})$.
  \label{alg:dml-retrain}
  \Comment{Retrain with DML on the current labelled set.}

  \State $\mathcal D_Z^{(l)} \gets \left\{(\mathbf z_i,y_i): \mathbf z_i  \coloneq  \mathbb E_{q_{\phi_l^\ast}(\cdot\mid\mathbf x_i)}[\mathbf z],\; (\mathbf x_i,y_i)\in\mathcal D_L^{(l)} \right\}$.
  \label{alg:dml-encode}

  \State
  \(\mathcal D_L^{(l;0)}\gets \mathcal D_L^{(l)}\);
  \quad
  \(\mathcal D_Z^{(l;0)}\gets \mathcal D_Z^{(l)}\);
  \quad
  \(\mathcal R_{(l;0)}\gets \mathcal R_0\).
  \Comment{Initialise inner loop and SDR region.}

  \State $q_l\gets \min\{q,\,B-(l-1)q\}$.

  \For{$k=0$ \textbf{to} $q_l-1$}
  \label{alg:dml-inner-start}

    \State Fit GP surrogate \(h_{l;k}:\mathcal Z\to\mathbb R\) on
    \(\mathcal D_Z^{(l;k)}\).

    \State
    \(u_{l;k}(\mathbf z)
    \coloneq
    u(\mathbf z\mid\mathcal D_Z^{(l;k)})\).
    \Comment{Current EI acquisition.}

    \State
    \(\hat{\mathbf z}_{l;k+1}
    \gets
    \operatorname*{arg\,max}_{\mathbf z\in\mathcal R^{(l;k)}}
    u_{l;k}(\mathbf z)\).

    \State
    \(\hat{\mathbf x}_{l;k+1}
    \sim
    p_{\theta_l^\ast}(\cdot\mid\hat{\mathbf z}_{l;k+1})\);
    \quad
    \(y_{l;k+1}\gets f(\hat{\mathbf x}_{l;k+1})\).

    \State
    \(\mathcal D_L^{(l;k+1)}
    \gets
    \mathcal D_L^{(l;k)}
    \cup
    \{(\hat{\mathbf x}_{l;k+1},y_{l;k+1})\}\).

    \State
    \(\mathcal D_Z^{(l;k+1)}
    \gets
    \mathcal D_Z^{(l;k)}
    \cup
    \{(\hat{\mathbf z}_{l;k+1},y_{l;k+1})\}\).

    \State \(t_{l;k+1}\coloneq (l-1)q+k+1\).

    \If{\(t_{l;k+1}\bmod K_{\mathrm{SDR}}=0\)}
      \State
      \(\mathcal R_{(l;k+1)}
      \gets
      \mathrm{SDR\_Update}
      \bigl(\mathcal R_{(l;k)},\mathcal D_Z^{(l;k+1)}\bigr)\).
    \Else
      \State
      \(\mathcal R_{(l;k+1)}
      \gets
      \mathcal R_{(l;k)}\).
    \EndIf

  \EndFor
  \label{alg:dml-inner-end}

  \State $\mathcal D_L^{(l+1)} \gets \mathcal D_L^{(l;q_l)};\quad \mathcal D_Z^{(l+1)} \gets \mathcal D_Z^{(l;q_l)}$.
\EndFor

\State \Return $f_{\min} \gets \min\{\,y:(\mathbf x,y)\in\mathcal D_L^{(L+1)}\,\}$.

\end{algorithmic}
\end{algorithm}
\normalsize

At outer retraining stage \(l\), the current labelled ambient dataset is $\mathcal D_L^{(l)}=\{(\mathbf x_i,y_i)\}_{i=1}^{N_l}$, where $ y_i \coloneq f(\mathbf x_i)$. Using the notation of Section~\ref{Prelim}, the VAE is warm-started from the previous checkpoint and retrained with the DML-augmented ELBO objective \eqref{VAE_DML_augmented_objective}. The metric term $\mathcal{L}_{\mathrm{metric}}(\phi;\mathcal D_L^{(l)})$ is formed from valid anchor-positive-negative triplets constructed from \(\mathcal D_L^{(l)}\), after normalising the observed objective values to $[0,1]$. Thus, every labelled point with at least one positive and one negative partner may serve as an anchor\footnote{We note that the current incumbent is not the only anchor.}. The triplet term is not constructed from the large unlabelled pre-training dataset \(\mathcal D_U\), which keeps the additional computational burden manageable. In implementation, the triplet contribution is evaluated over the valid triplets in the current training batch, giving a minibatch estimate of \(\mathcal L_{\mathrm{metric}}(\phi;\mathcal D_L^{(l)})\). After retraining, the labelled dataset \(\mathcal D_L^{(l)}\) is re-encoded using the updated encoder \(q_{\phi_l^\ast}\), giving the latent GP dataset \(\mathcal D_Z^{(l)}\). BO then proceeds for at most \(q\) inner iterations in the latent space. At each inner step, a GP surrogate is fitted on the current latent data, EI is maximised over the fixed latent search region \(\mathcal R_0\), the selected latent point is decoded to an ambient point, and \(f\) is evaluated there.

While Algorithms~\ref{BOVAE-retraining} and~\ref{BOVAE-with-DML} adopt a fixed-period policy that retrains the VAE and accepts every update at predetermined intervals, this formulation serves as the algorithmic definition. Section~\ref{subsec:latent-dimension-representation-gap} uses it to reveal how the retraining policy affects the representation gap, whereas Section~\ref{sec:bovae-benchmark-performance} introduces an adaptive proposal-and-acceptance policy that varies the retraining interval and adopts an update only when it passes validation reconstruction error checks.

\section{Numerical Experiments}\label{Results}

We conduct numerical experiments and studies to evaluate the three BO--VAE algorithms (Algorithms~\ref{BOVAE-SDR}, \ref{BOVAE-retraining}, \ref{BOVAE-with-DML}) and to examine the main design choices that affect their performance. The results are organised in five parts. We first illustrate the effects of SDR within the VAE-generated latent spaces using Ackley and Rosenbrock to isolate the contribution of domain reduction in the latent space. We then study the effects of the latent dimension $d$ on the optimisation quality with Ackley, Rosenbrock, and Rastrigin. In the study, we further examine the VAE-induced representation gaps through empirical diagnostics. Next, we visualise the optimisation traces of the algorithms to compare optimisation performance as a function of the number of evaluations and total computation time. We subsequently compare the three algorithms through performance and data profiles and benchmark them against the standard BO-SDR (Alg. \ref{BOSDR}), using the full-rank test functions in Table~\ref{tab:full-rank-test-problems}. Finally, we compare our BO--VAE algorithms against EGORSE \cite{egorse} to assess their performance over test problems with nonlinear effective mappings. This comparison examines the ability of VAEs to capture nonlinear low-dimensional structure. 

Throughout the experiments, except the comparisons against EGORSE, we adopt the VAEs with the architectures shown in Table~\ref{tab:CN-reference VAE family}. For each $d$, the VAE training domain is fixed to $[-5, 5]^D$. We rescale between the problem domain and the VAE training domain so that the VAE architecture is independent of the benchmark domains. We use the initial design with $D$ points for a D-dimensional test problem, and set the budget of enrichment evaluations $B \coloneq 800$. Performance is reported via the \emph{simple regret}, defined as the difference between the minimum function value found by the algorithm and the known optimal function value. Further configuration details are provided in Appendix~\ref{Appendix Numerical Experiments Details}.

\begin{table}[h]
\caption{The series of VAEs used in the numerical experiments. Bracketed lists give layer widths (encoder: input$\to$latent; decoder: latent$\to$output). $[10,5,2]$ indicates a three-layer feedforward neural network with an input layer of width $10$, a hidden layer of width $5$, and a latent layer of width $2$.}
\label{tab:CN-reference VAE family}
\centering
\setlength{\tabcolsep}{3pt}
\renewcommand{\arraystretch}{1.1}

\begin{tabularx}{\linewidth}{@{}l c Y Y Y c c@{}}
\toprule
VAE ID & $D$ & $d$ & Encoder & Decoder
& \shortstack{Hidden\\activation}
& \shortstack{Decoder\\output} \\
\midrule

VAE--D10
& 10
& $[1,\allowbreak 2,\allowbreak 3,\allowbreak 4,\allowbreak
   5,\allowbreak 6,\allowbreak 8,\allowbreak 10]$
& $[10,\allowbreak 64,\allowbreak 32,\allowbreak d]$
& $[d,\allowbreak 32,\allowbreak 64,\allowbreak 10]$
& SiLU
& Linear \\

VAE--D100
& 100
& $[2,\allowbreak 10,\allowbreak 25,\allowbreak 50,\allowbreak 75]$
& $[100,\allowbreak 512,\allowbreak 256,\allowbreak d]$
& $[d,\allowbreak 256,\allowbreak 512,\allowbreak 100]$
& SiLU
& Linear \\

\botrule
\end{tabularx}
\end{table}

For the experiments, we adopt the \textit{BoTorch} solver for its availability in the GPU-based environment, and implement a GPU-compatible SDR update within this framework. Prior public SDR implementations targeted \emph{BayesOpt} \cite{Bayes_Opt} on CPU. For completeness, we also compare three BO software packages, \textit{GPyOpt} \cite{GPyOpt2016}, \textit{BoTorch} \citet{balandat2020botorch}, and \textit{BayesOpt} \cite{Bayes_Opt}, for noisy and smooth $f$ using the data and performance profiles introduced in \cite{More2009}. The test functions are listed in Table~\ref{tab:solver-test}. In this comparison, \textit{BoTorch} outperforms the other two, and the details are included in Appendix~\ref{appendix: solver results}.

\subsection{Algorithm Implementations}\label{Algorithm Implementations}

Our implementation is inspired by \cite{grosnit2021highdimensionalbayesianoptimisationvariational} but departs from their codebase in several practical ways relevant to general-purpose, black-box mathematical optimisation:

\begin{enumerate}
    \item \textbf{Codebase and modularity}: The official repository in \cite{grosnit2021highdimensionalbayesianoptimisationvariational} is tightly coupled to domain-specific pipelines (e.g. molecular design and gene-expression reconstruction). This coupling makes it difficult to isolate components for reuse in general function optimisation. We therefore re-implemented the full stack from scratch with an explicitly modular design: data handling, VAE pre-training, BO in the latent space, and DML fine-tuning are exposed as independent modules. This separation allows researchers to swap components and reuse only what is needed for their tasks.\\
    
    \item \textbf{Pre-training strategy aligned with BO’s data regime}: While Algorithm~\ref{BOVAE-with-DML} in \cite{grosnit2021highdimensionalbayesianoptimisationvariational} uses a DML objective $\mathcal{L}_{DML}$ during pre-training to enforce a strongly structured latent space, we found this less suitable in the BO setting where labelled evaluations are scarce and progressively acquired. VAE pretraining typically benefits from abundant data; incorporating objective values at this stage would couple representation learning to the limited supervision available to BO... and can undermine BO’s sample-efficiency rationale. In our pipeline, the VAE is first pre-trained with the standard ELBO objective (unsupervised) to ensure robust reconstruction under limited supervision. Only after BO begins to accumulate informative query points do we fine-tune the latent geometry with DML. This sequencing preserves the separation between representation learning and the limited supervision available to BO, while remaining compatible with the DML objective used in \cite{grosnit2021highdimensionalbayesianoptimisationvariational} when sufficient task-specific signal is present.
\end{enumerate}

We emphasise that these are pragmatic design choices aimed at broadening applicability and facilitating fair component-wise comparisons. 

To initialise and update SDR consistently in the latent space, we avoid keeping $\mathcal{R}_0$ fixed as the latent search box throughout the optimisation. VAE retraining, particularly when DML is used, may alter the location, scale, and correlation structure of the latent representation. Consequently, an SDR region defined using the raw latent coordinates before retraining may no longer correspond to the same search region after the representation is updated. We therefore standardise the raw VAE latent coordinates and map them to a unit box, within which BO is performed and the SDR updates are applied. Further details are provided in Appendix~\ref{Appendix: Additional Implementation Details}.

\subsection{SDR in VAE Latent Spaces}\label{subsection: SDR in VAE-generated Latent Spaces}

\begin{figure}[h]
  \centering
  \includegraphics[width=\linewidth]{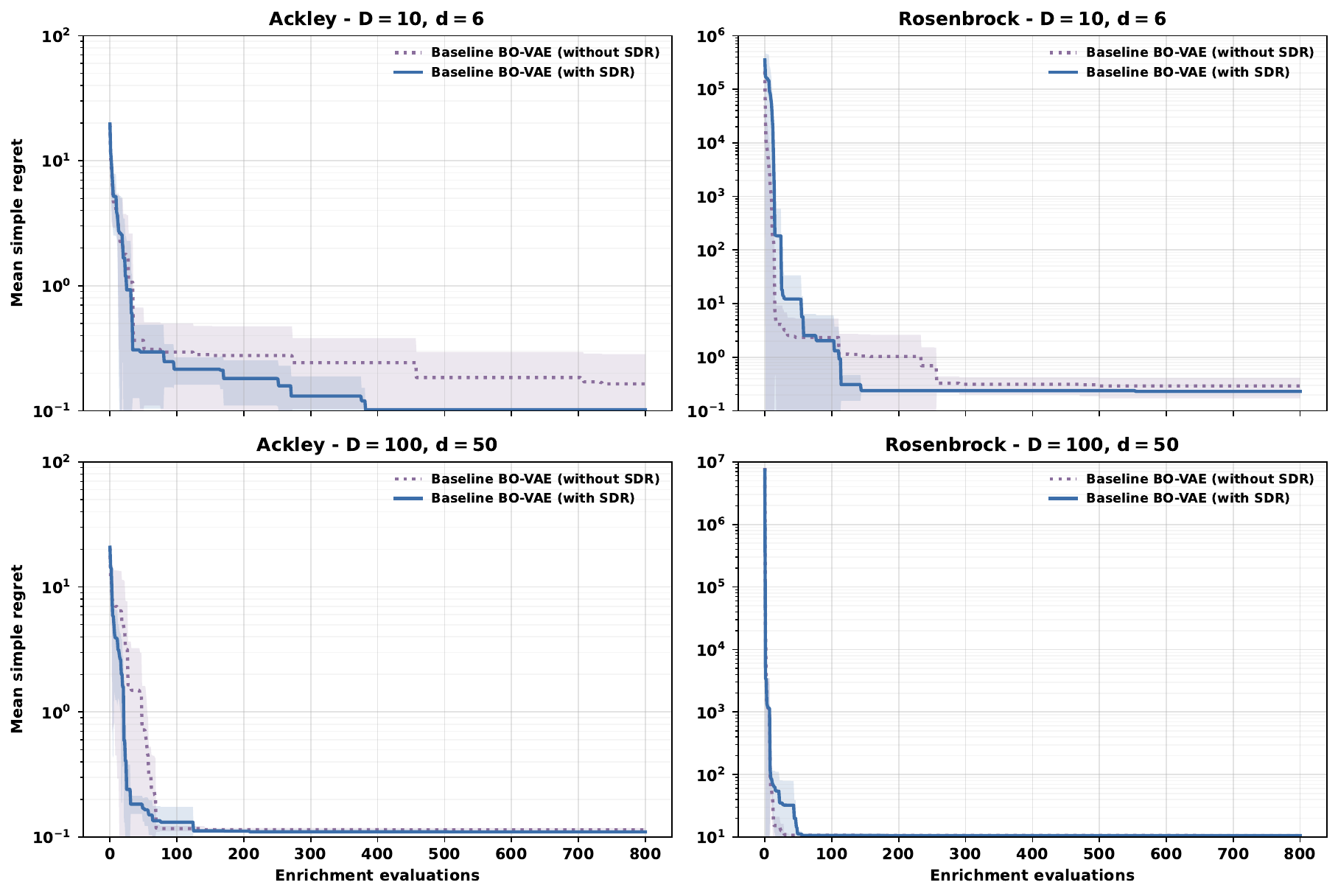}
  \caption{Simple-regret trajectories of the baseline BO--VAE  (Alg.~\ref{BOVAE-SDR}) with and without SDR on Ackley and Rosenbrock for
  $(D,d)=(10,6)$ and $(D,d)=(100,50)$. Lower values are better. Each curve shows the arithmetic mean over five distinct seeds, and the shaded region represents one sample standard deviation. The SDR region is updated every $K_{\mathrm{SDR}}=10$ enrichment evaluations in both configurations.}
  \label{fig:sdr-vae}
\end{figure}

To examine the effect of SDR in VAE-induced latent spaces, we consider the configurations $(D,d)=(10,6)$ and $(D,d)=(100,50)$. The baseline BO--VAE
algorithm (Alg.~\ref{BOVAE-SDR}) is evaluated with and without SDR on Ackley and Rosenbrock. As shown in Fig.~\ref{fig:sdr-vae}, SDR provides the clearest benefit for $(D,d)=(10,6)$. It improves the final performance on Ackley and accelerates the intermediate-stage convergence on Rosenbrock, although the difference on Rosenbrock narrows toward the end of the optimisation.

For $(D,d)=(100,50)$, the curves with and without SDR largely overlap on both problems. This does not imply that the SDR region fails to contract. Instead, the contraction does not translate into a material improvement in optimisation performance. Overall, SDR is more effective in the lower-dimensional latent space, where identifying and exploiting a promising local region is more tractable. In the higher-dimensional latent space, the quality of the incumbent and the difficulty of concentrating the search within a useful neighbourhood limit the practical benefit of contraction.

\FloatBarrier

\subsection{Latent-Dimension Effects and VAE-Induced Representation Gaps}
\label{subsec:latent-dimension-representation-gap}

We investigate how the latent dimension $d$ affects optimisation performance and whether the resulting performance is limited by the latent-space search or by the representational capacity of the VAE. We conduct controlled latent-dimension scans at fixed ambient dimensions $D=10$ and $D=100$, using the three BO-VAE algorithms and the corresponding VAEs summarised in Table~\ref{tab:CN-reference VAE family}. Specifically, we consider $d\in\{1,2,3,4,5,6,8,10\}$ for $D=10$ and $d\in\{2,10,25,50,75\}$ for $D=100$, while holding the remaining experimental protocol fixed within each ambient dimension.

In addition to comparing the achieved simple regrets, we quantify the objective-relevant limitations imposed by each VAE representation. Let
$f:\mathcal X\rightarrow\mathbb R$ denote the original minimisation objective
with known optimum value $f^\star$. 

For an adopted VAE state $s$, let $r_{\theta_s}$ denote its decoder, let $z_s:\mathcal U_d\rightarrow\mathbb R^d$ denote the stored transformation from the global unit latent box $\mathcal U_d=[0,1]^d$, and let $T:\mathcal B_{\mathrm{VAE}}\rightarrow\mathcal X$ denote the map from the usual VAE training domain to the canonical objective domain. For the test function $f:\mathcal X\rightarrow\mathbb R$ with known optimum value $f^\star$. The complete ambient-domain mapping induced by state $s$ is therefore
\begin{equation}
  \mathbf{x}_s(\mathbf{u})
  =
  T\!\left(
    \Pi_{\mathcal B_{\mathrm{VAE}}}
    \bigl(r_{\theta_s}(z_s(\mathbf{u}))\bigr)
  \right),
  \qquad \mathbf{u} \in\mathcal U_d,
\end{equation}
and its decoder-supported set is
\begin{equation}
  \mathcal S_s
  =
  \left\{\mathbf{x}_s(\mathbf{u}): \mathbf{u}\in\mathcal U_d\right\}.
\end{equation}

We define the \emph{VAE-induced representation gap}, equivalently the \emph{decoder-support floor}, of the retraining stage $s$ as
\begin{equation}
  \Delta_{\mathrm{rep}}(s)
  =
  \inf_{\mathbf{x}\in\mathcal S_s}
  \bigl[f(\mathbf{x})-f^\star\bigr]_+
  =
  \inf_{\mathbf{u} \in\mathcal U_d}
  \bigl[f(\mathbf{x}_s(\mathbf{u}))-f^\star\bigr]_+,
  \qquad
  [a]_+=\max(a,0).
  \label{eq:vae-representation-floor}
\end{equation}
This quantity is an intrinsic, objective-specific property of the adopted decoder and its latent domain. While state $s$ remains fixed, even an ideal latent-space optimiser cannot generate a decoded point whose objective gap is smaller than $\Delta_{\mathrm{rep}}(s)$.

The exact infimum in Eq.~\eqref{eq:vae-representation-floor} is generally unavailable. We therefore report a numerical estimate $\widehat{\Delta}_{\mathrm{rep}}(s)$ obtained by searching the complete global latent box and evaluating the resulting decoded points. Since this procedure returns the best feasible point found by a finite numerical search,
\begin{equation}
  \Delta_{\mathrm{rep}}(s)
  \leq
  \widehat{\Delta}_{\mathrm{rep}}(s),
\end{equation}
up to numerical evaluation error. Thus, $\widehat{\Delta}_{\mathrm{rep}}(s)$ is a feasible upper estimate of the representation floor, rather than a certified lower bound. The known optimum is used only for this retrospective diagnostic and for computing simple regret.

For run $j$, let $\mathcal A_j$ denote the set of all VAE states adopted during the run, including the initial state, and let $s_{j,\mathrm{final}}$ denote the final adopted state. We define the best-ever and final estimated representation floors as
\begin{align}
  \widehat{\Delta}_{j,\mathrm{best}}
  &=
  \min_{s\in\mathcal A_j}
  \widehat{\Delta}_{\mathrm{rep}}(s),\\
  \widehat{\Delta}_{j,\mathrm{final}}
  &=
  \widehat{\Delta}_{\mathrm{rep}}
  \bigl(s_{j,\mathrm{final}}\bigr).
\end{align}
For a fixed VAE, these two quantities coincide. For methods with retraining, their difference measures whether the objective-relevant support of the final adopted VAE is better or worse than that of the best representation encountered during the run.

Let $\mathcal Q_j$ contain all initial and enrichment points evaluated in run $j$. Its achieved simple regret is
\begin{equation}
  r_j
  =
  \left[
    \min_{\mathbf{x}\in\mathcal Q_j} f(\mathbf{x})-f^\star
  \right]_+.
\end{equation}
For a displayed set of runs $J$, the plotted quantities are the arithmetic
means
\begin{equation}
  \overline r
  =
  \frac{1}{|J|}\sum_{j\in J}r_j,
  \qquad
  \overline{\Delta}_{\mathrm{best}}
  =
  \frac{1}{|J|}
  \sum_{j\in J}\widehat{\Delta}_{j,\mathrm{best}},
  \qquad
  \overline{\Delta}_{\mathrm{final}}
  =
  \frac{1}{|J|}
  \sum_{j\in J}\widehat{\Delta}_{j,\mathrm{final}}.
  \label{eq:mean-representation-floors}
\end{equation}

These quantities distinguish representational limitations from optimisation limitations. An achieved regret close to the best-ever representation floor indicates that BO largely exploits the best decoder support made available during the run. A large separation indicates that the representation contains better objective values that remain unused because of finite-budget latent search, GP modelling, acquisition optimisation, or SDR. If the final floor is higher than the best-ever floor, a later accepted retraining step has degraded objective-relevant decoder support. Nevertheless, an incumbent can remain better than the final floor because it may have been generated by an earlier VAE state or belong to the initial design. These comparisons are therefore diagnostic rather than an exact additive decomposition of optimisation error. The representation floor should also be distinguished from average reconstruction error, which does not directly measure objective-specific reachability.

We conduct this study using Ackley, Rosenbrock, and Rastrigin. For each test function and each pair of D and d, ten distinct seeds are used. Results are reported as the \emph{mean objective gap to $f^{\star}$}. Rows correspond to Ackley, Rosenbrock, and Rastrigin, while columns correspond to Baseline BO--VAE (Alg. \ref{BOVAE-SDR}), periodically retrained BO--VAE (Alg. \ref{BOVAE-retraining}), and periodically retrained BO--VAE with DML (Alg. \ref{BOVAE-with-DML}). Since the purpose of this subsection is to study the latent dimension effects rather than retraining-frequency selection, we fix the retraining period at $q\coloneq100$ for both retraining variants. The same value is used for every objective, ambient dimension, latent dimension, and random seed. With the enrichment budget \(B\coloneq800\), this setting provides multiple representation updates during each run while allowing approximately \(100\) additional objective-labelled observations to accumulate between successive updates. It also keeps the computational cost of the full latent-dimension scan manageable.

The value \(q=100\) is used solely as a controlled experimental setting to prevent retraining frequency from becoming an additional confounding factor as \(d\) varies. It was not obtained by optimising \(q\), and the comparisons in this subsection are conditional on this fixed schedule. We therefore make no claim that \(q=100\) is generally optimal, or even preferable to other fixed retraining periods. The main algorithm-comparison experiments in Sections~\ref{sec:bovae-benchmark-performance}--\ref{sec:bovae-egorse-nonlinear} instead use the adaptive retraining policy introduced in Section~\ref{sec:bovae-benchmark-performance}, under which the interval \(q_t\) changes in response to observed retraining outcomes. 

\begin{figure}[h]
    \centering
    \includegraphics[width=\linewidth]{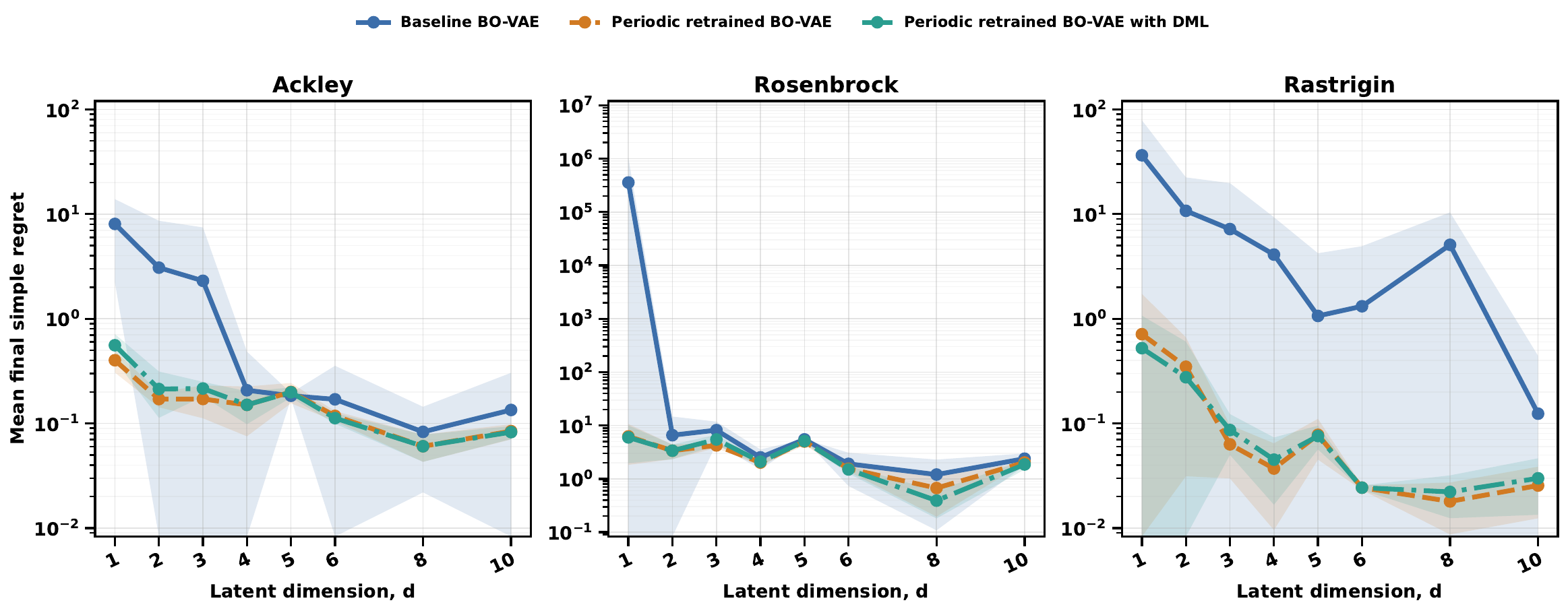}
    \caption{Effect of the latent dimension $d$ on BO--VAE performance for $D=10$ on Ackley, Rosenbrock, and Rastrigin. The latent dimensions $d\in\{1,2,3,4,5,6,8,10\}$ are considered. Each marker denotes the arithmetic mean of the simple regret, and the shaded region represents one sample standard deviation. Lower values are better, and the y--axis is shown on a logarithmic scale. The retraining period is fixed at $q\coloneq 100.$}
    \label{fig: D10 latent dim study simple regret}
\end{figure}

\begin{figure}[h]
    \centering
    \includegraphics[width=\linewidth]{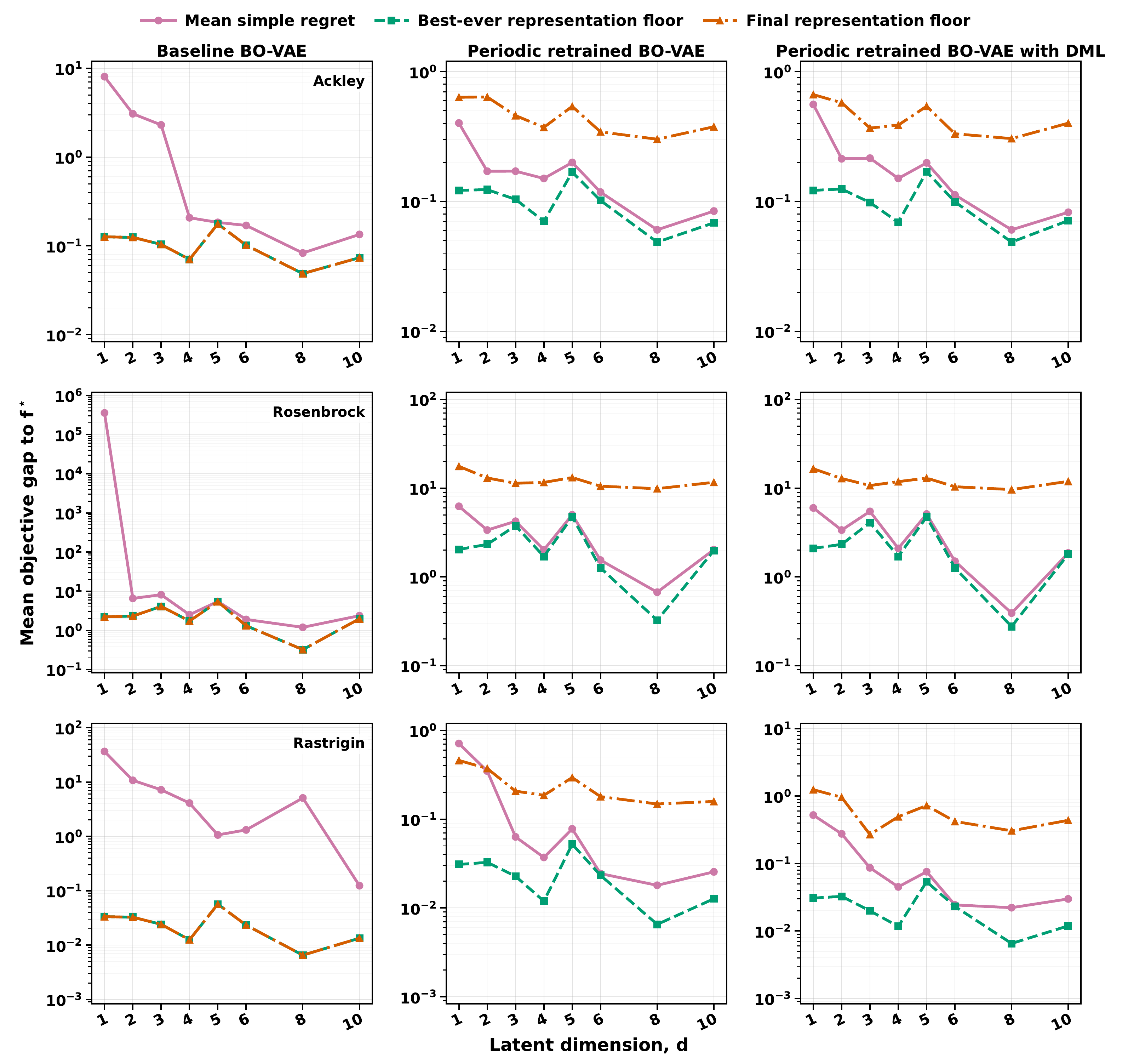}
    \caption{Achieved simple regret and estimated VAE-induced representation floors across latent dimensions for $D=10$. The pink curve shows the mean achieved simple regret $\overline r$. The green and orange curves show the mean best-ever and final representation-floor estimates, $\overline{\Delta}_{\mathrm{best}}$ and $\overline{\Delta}_{\mathrm{final}}$, respectively. They are numerical estimates for theoretical representation floors defined in~\eqref{eq:vae-representation-floor}. For the baseline BO-VAE, the best-ever and final floors coincide. A large separation between achieved regret and the best-ever floor indicates reachable performance that was not exploited by the finite-budget optimisation procedure, whereas a final floor above the best-ever floor indicates that a later accepted retraining state has degraded objective-relevant decoder support. All y--axes show mean objective gaps to $f^\star$ on logarithmic scales. The retraining period is fixed at $q\coloneq 100.$}
    \label{fig: D10 latent dim study vae representation floor}
\end{figure}

\paragraph{Results for \texorpdfstring{$D=10$}{}.} Figures \ref{fig: D10 latent dim study simple regret} and \ref{fig: D10 latent dim study vae representation floor} show the mean simple regret and representation-gap diagnostics, respectively. In the Ackley and Rosenbrock baseline panels, increasing $d$ removes the large small-$d$ search gap and the regret approaches the representation floor.  The retraining panels show substantially lower best-ever floors than final floors at several dimensions, so the strongest incumbent can be associated with an earlier decoder rather than the terminal one.  In the Rastrigin baseline panel, the smallest support floor occurs near $d=8$, whereas the smallest displayed regret occurs at $d=10$. A similar pattern appears for the two BO--VAE retraining variants, where the best-ever representation floor happens at $d = 8$ with the smallest mean simple regret at $d = 6$. This indicates the existence of a search-versus-support trade-off rather than evidence that the $d=10$ decoder has the better floor.


\FloatBarrier

\begin{figure}[h]
    \centering
    \includegraphics[width=\linewidth]{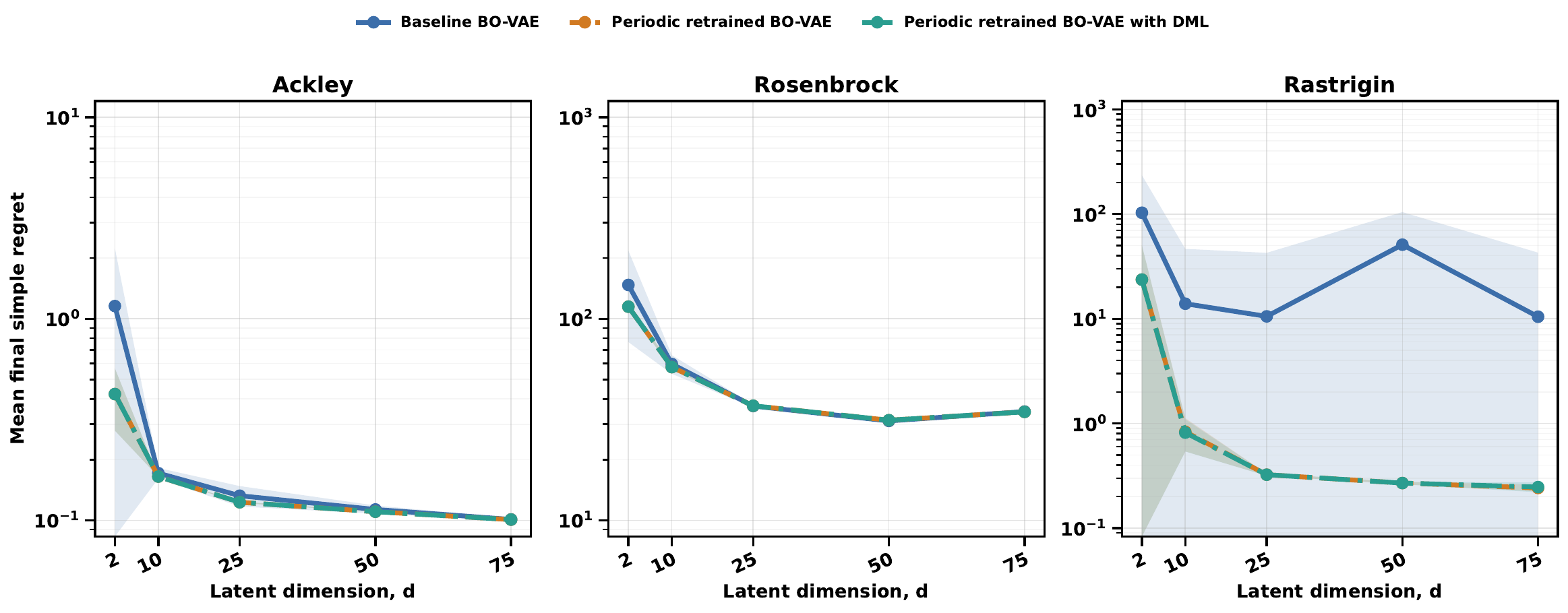}
    \caption{Effect of the latent dimension $d$ on BO--VAE performance for $D=100$ on Ackley, Rosenbrock, and Rastrigin. The latent dimensions $d\in\{2,10,25,50,75\}$ are considered. Each marker denotes the arithmetic mean of the simple regret, and the shaded region represents one sample standard deviation. Lower values are better, and the y--axis is shown on a logarithmic scale. The retraining period is fixed at $q\coloneq 100.$}
    \label{fig: D100 latent dim study simple regret}
\end{figure}

\begin{figure}[h]
    \centering
    \includegraphics[width=\linewidth]{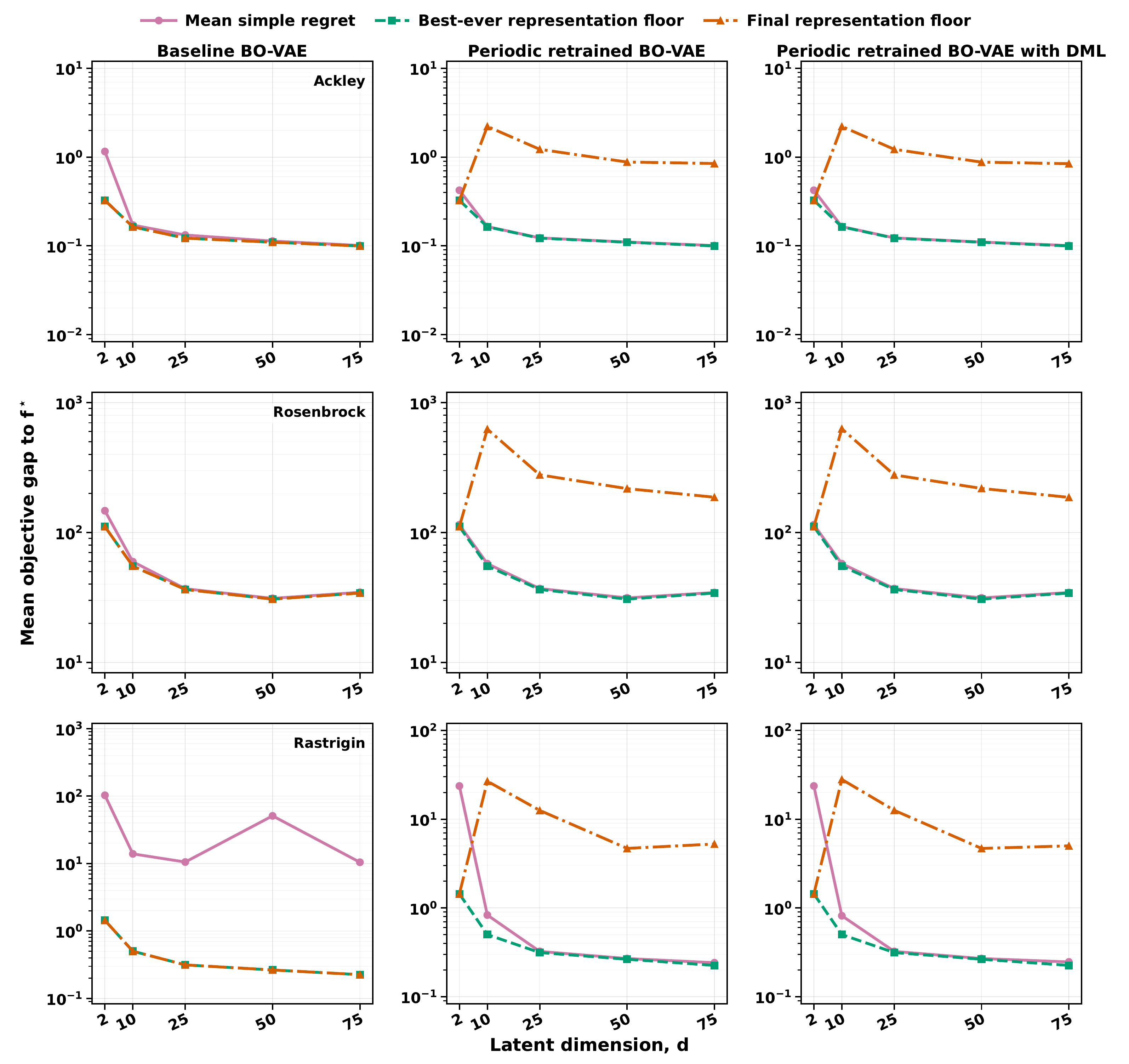}
    \caption{Achieved simple regret and estimated VAE-induced representation floors across latent dimensions for $D=100$. The pink curve shows the mean achieved simple regret $\overline r$; the green and orange curves show the mean best-ever and final representation-floor estimates, $\overline{\Delta}_{\mathrm{best}}$ and $\overline{\Delta}_{\mathrm{final}}$, respectively. They are numerical estimates for theoretical representation floors defined in~\eqref{eq:vae-representation-floor}. For the baseline BO-VAE, the best-ever and final floors coincide. A large separation between achieved regret and the best-ever floor indicates reachable performance that was not exploited by the finite-budget optimisation procedure, whereas a final floor above the best-ever floor indicates that a later accepted retraining state has degraded objective-relevant decoder support. All y--axes show mean objective gaps to $f^\star$ on logarithmic scales. The retraining period is fixed at $q\coloneq 100.$}
    \label{fig: D100 latent dim study vae representation floor}
\end{figure}

\paragraph{Results for \texorpdfstring{$D=100$}{}.} Figures \ref{fig: D100 latent dim study simple regret} and \ref{fig: D100 latent dim study vae representation floor} present the corresponding results for $D = 100$. The baseline panels show close agreement between simple regret and the fixed representation floor at the strongest dimensions, approximately $d=75$ for Ackley, $d=50$ for Rosenbrock, and $d=75$ for Rastrigin.  In both retraining columns, the best-ever floor and achieved regret generally improve toward these larger dimensions, while the final floor can remain orders of magnitude higher, especially after high-dimensional retraining. This implies that an earlier adopted VAE state supplied the useful support, BO found and retained an incumbent there, and subsequent retraining degraded the final decoder without erasing that incumbent.

\FloatBarrier

Across these figures, the two retraining variants demonstrate that the quality of the final VAE representation and the quality of the retained incumbent need not evolve in parallel. At several latent dimensions, the final representation floor is higher than the best-ever floor, while the algorithm retains an incumbent generated under an earlier, more favourable representation. This mismatch can arise because the current retraining losses, although useful for reconstruction and metric learning, do not directly control the objective-specific decoder-support floor. An important direction for future work is therefore to develop stable, objective-aware retraining losses and acceptance criteria that improve searchability while preserving access to high-quality regions discovered by earlier VAE states.

We note that this scan serves as a diagnostic study over a prespecified set of candidate latent dimensions. It demonstrates that the choice of $d$ materially affects both decoder support and latent optimisation difficulty, but it does not provide an estimator or automatic selection rule for \(d\). Consequently, the problem of selecting \(d\) when the effective dimension is unknown remains open.

\subsection{Performance of BO--VAE Algorithms on Benchmark Problems}
\label{sec:bovae-benchmark-performance}

The controlled study in Section~\ref{subsec:latent-dimension-representation-gap} shows that a VAE update accepted under the fixed \(q=100\) schedule can sometimes worsen the final objective-relevant representation floor, even though an incumbent obtained under an earlier representation is retained. This observation does not establish that \(q=100\) causes the degradation, nor does it compare \(q=100\) with alternative fixed periods. Rather, it identifies a limitation of adopting every proposed VAE update according to a predetermined schedule without evaluating the quality of the resulting representation.

The retraining period $q$ controls the trade-off between representation adaptability and computational overhead. A small value of $q$ allows the VAE to respond rapidly to newly observed data, but increases retraining cost and may adapt too aggressively to the limited objective-labelled dataset available early in BO. A large value reduces this cost and provides greater representation stability, but may leave an unsuitable latent representation unchanged for too long. These competing effects motivate an adaptive proposal-and-acceptance policy, where retraining is proposed at a variable interval, and the resulting VAE is adopted only if it passes certain quality checks, such as a held-out reconstruction test.

Let \(q_t\) denote the number of enrichment evaluations until the next retraining proposal. The adaptive policy is initialised at \(q_0\coloneq50\) and restricts the interval to
\[
    q_t\in\{50,100,200\}.
\]
A rejected proposal or a substantial observed improvement in the incumbent resets the interval to $50$, whereas consecutive proposals providing little measured improvement can increase it to \(100\) and subsequently to \(200\). Consequently, retraining remains responsive when the optimisation trajectory or the proposed representation is changing, but becomes less frequent when successive updates provide little measured benefit. These values are operational lower and upper safeguards for the present experiments, and they were not obtained from a general optimisation of the retraining schedule. Accordingly, we evaluate the behaviour of this specified adaptive policy but do not claim that its initial value, bounds, or update rule are generally optimal.

At each retraining opportunity, the available observations are divided into a retraining set $\mathcal T_t$ and a validation set $\mathcal V_t$, with $\mathcal T_t\cap\mathcal V_t=\varnothing$. The validation set is held fixed throughout the corresponding retraining proposal and is not used to update the VAE parameters. Let $s_t^{-}$ denote the currently adopted VAE state and $s_t^{+}$ the proposed state obtained by retraining on $\mathcal T_t$. Rather than using the retraining loss to decide whether the proposal should be adopted, we compare the deterministic reconstruction errors of the two states on the same held-out validation set.

For a VAE state $s$, we define the validation reconstruction error as
\begin{equation}
  \mathcal E_{\mathrm{rec}}^{\mathrm{val}}
  \bigl(s;\mathcal V_t\bigr)
  =
  \frac{1}{|\mathcal V_t|D}
  \sum_{\mathbf{x}_i\in\mathcal V_t}
  \left\|
    \mathbf{x}_i-
    r_{\theta_s}\!\left(\mathbf{m}_{\phi_s}(\mathbf{x}_i)\right)
  \right\|_2^2,
  \label{eq:validation-reconstruction-error}
\end{equation}
where $\mathbf{m}_{\phi_s}(\mathbf{x}_i)$ is the posterior-mean encoding and $r_{\theta_s}$ is the corresponding decoder for the VAE state $s$. The proposed state passes the reconstruction gate if
\begin{equation}
  \mathcal E_{\mathrm{rec}}^{\mathrm{val}}
  \bigl(s_t^{+};\mathcal V_t\bigr)
  \leq
  (1+\tau_{\mathrm{rec}})
  \mathcal E_{\mathrm{rec}}^{\mathrm{val}}
  \bigl(s_t^{-};\mathcal V_t\bigr)
  +\varepsilon,
  \qquad
  \tau_{\mathrm{rec}}=0.10,
  \quad
  \varepsilon=10^{-12}.
  \label{eq:validation-reconstruction-acceptance}
\end{equation}
Thus, the candidate VAE may increase the held-out reconstruction MSE by at most $10\%$ relative to the currently adopted state. If this condition fails, the proposal is rejected, and the previous VAE state is retained. Evaluating both states on the same held-out data prevents the acceptance decision from being based on the loss used to fit the proposed VAE and reduces the risk of adopting an update that has overfitted the limited observations available at that stage of BO.

\begin{figure}[h]
    \centering
    \includegraphics[width=\linewidth]{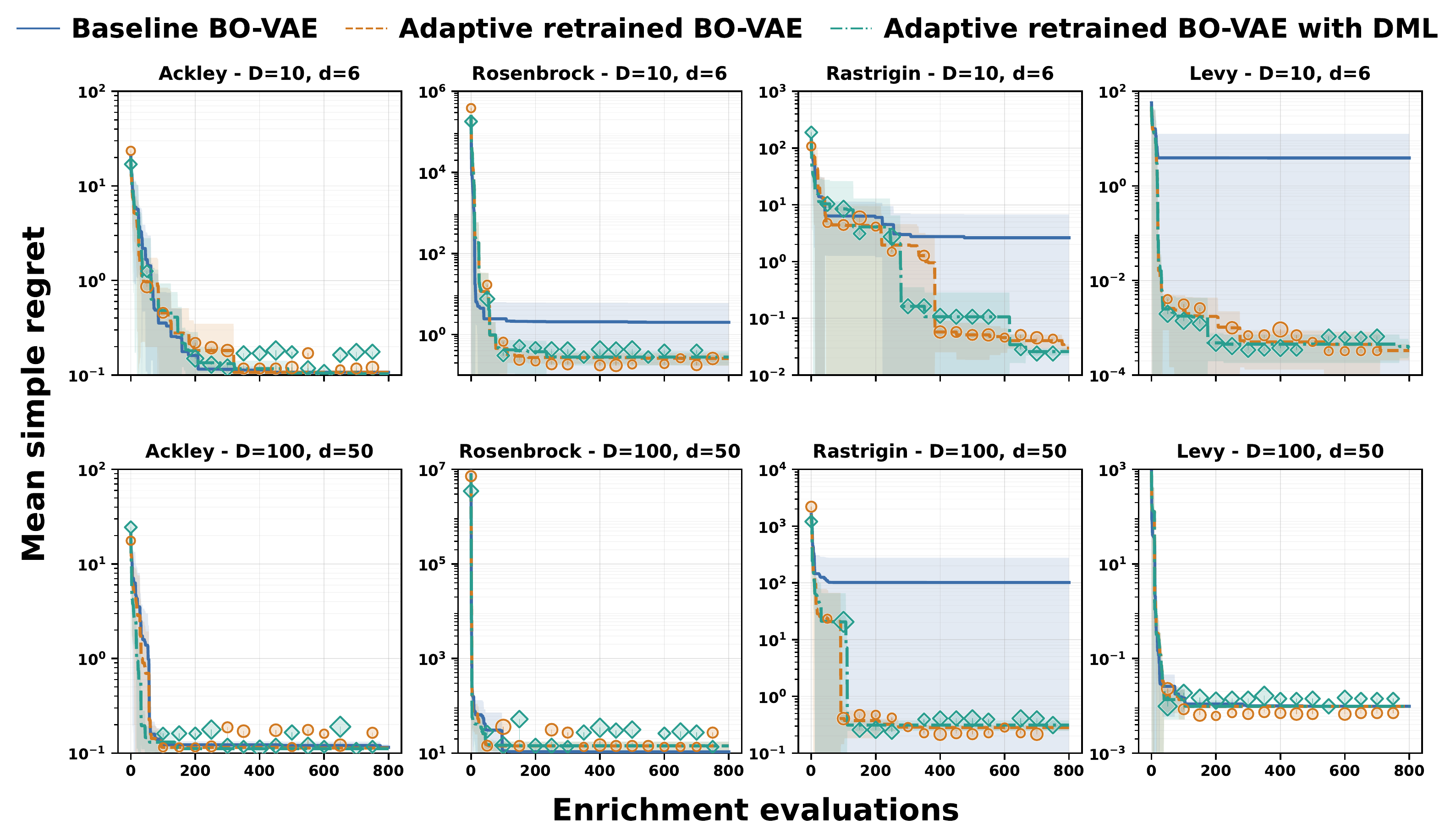}
    \caption{Convergence of the three BO--VAE variants on Ackley, Rosenbrock, Rastrigin, and L\'evy. The upper row shows $(D,d)=(10,6)$ and the lower row shows $(D,d)=(100,50)$. Each run uses $D$ initial observations followed by $800$ enrichment evaluations. The curves show the arithmetic mean of the simple regret, and the shaded regions represent one sample standard deviation over five distinct seeds. Lower values are better, and the y--axis is shown on a logarithmic scale. Markers indicate accepted VAE updates. Their sizes represent the number of displayed runs accepting an update at that iteration. For visualisation purposes only, we vertically separate the markers.}
    \label{fig: Full rank simple regret}
\end{figure}

\begin{figure}[h]
    \centering
    \includegraphics[width=\linewidth]{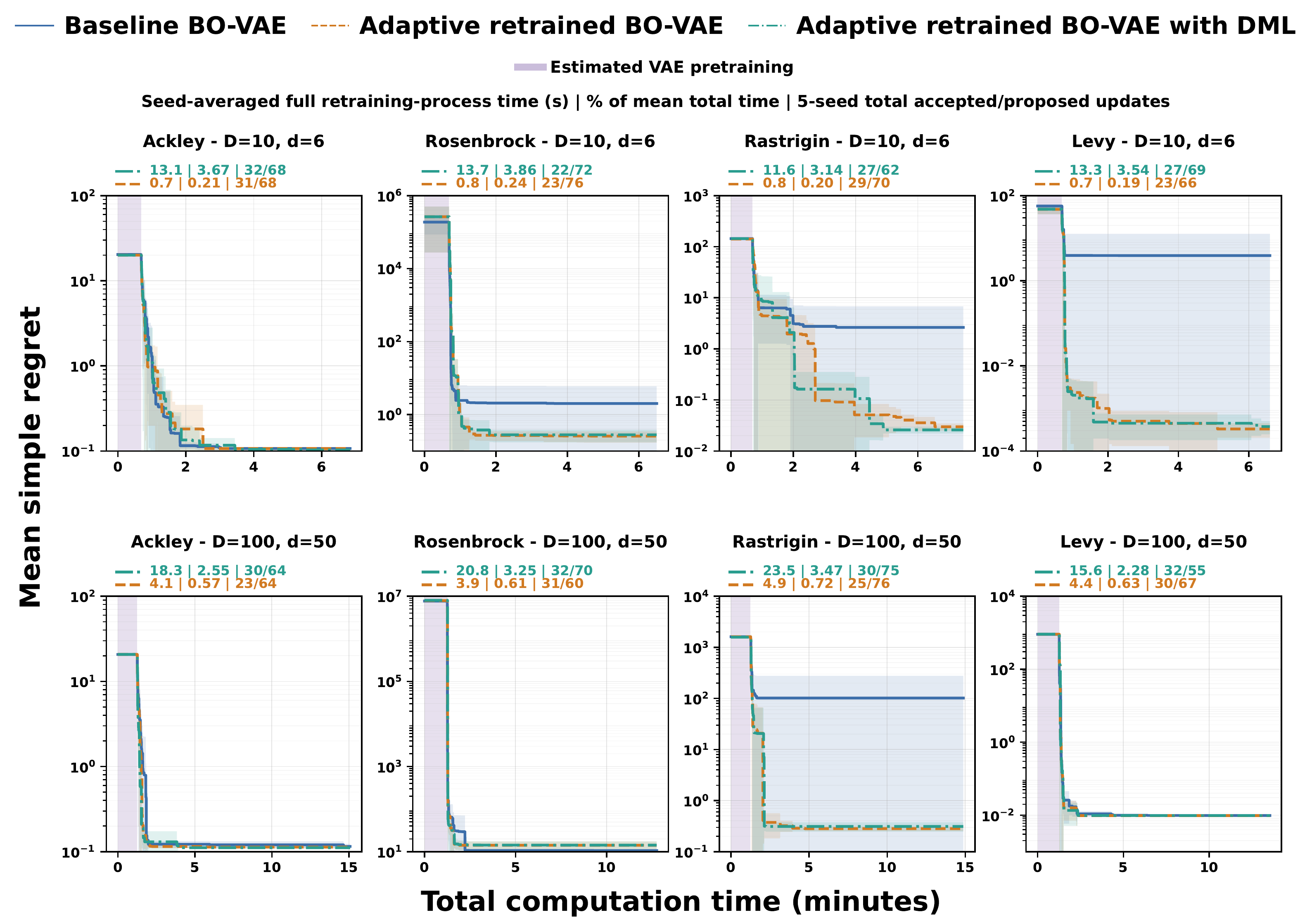}
    \caption{Simple regret of the three BO--VAE variants as a function of total accounted computation time on Ackley, Rosenbrock, Rastrigin, and L\'evy. The upper row shows $(D,d)=(10,6)$ and the lower row shows $(D,d)=(100,50)$. The curves show the arithmetic mean of the simple regret, and the shaded regions represent one sample standard deviation over the five distinct seeds. Lower values are better, and the y--axis is shown on a logarithmic scale. Total time includes the estimated cost of VAE pretraining, BO initialisation and enrichment, and every accepted or rejected retraining proposal. The grey region at the beginning of each panel denotes the estimated VAE pretraining time. Full retraining-process time is calculated in seconds by summing all retraining-process intervals within each seed, and then averaging the totals across the five seeds. It includes VAE retraining, validation, and acceptance assessment. The reported percentage is the mean full retraining-process time as a percentage of the total accounted computation time. The accompanying fraction gives the number of accepted VAE updates out of the total number of proposed updates in these 5-seed runs.}
    \label{fig: Full rank computation time}
\end{figure}

Using this adaptive policy, we evaluate Baseline BO--VAE, adaptively retrained BO--VAE, and adaptively retrained BO--VAE with DML on Ackley, Rosenbrock, Rastrigin, and L\'evy. The configurations $(D,d)=(10,6)$ and $(D,d)=(100,50)$ are adopted based on the preceding latent-dimension study. We report simple regret against both the number of enrichment evaluations and the total accounted computation time. The former measures
objective-evaluation efficiency, whereas the latter determines whether any optimisation improvement remains after accounting for VAE pretraining and the cost of all accepted and rejected retraining proposals. The results are shown in Figures \ref{fig: Full rank simple regret} and \ref{fig: Full rank computation time}. Figure~\ref{fig: Full rank simple regret} shows that the benefit of adaptive VAE retraining is strongly problem-dependent. For $D=10$, both retraining variants substantially reduce the simple regret on Rosenbrock, Rastrigin, and L\'evy, but provide little improvement on Ackley. For $D=100$, the largest benefit occurs on Rastrigin. The algorithms perform similarly on Ackley and L\'evy, while the fixed-VAE baseline performs better on Rosenbrock. Moreover, adding DML does not consistently improve upon reconstruction-based retraining alone. The update markers also show that accepting a new VAE does not necessarily produce an immediate reduction in simple regret.

Figure~\ref{fig: Full rank computation time} shows that the main iteration-based improvements, particularly for Rastrigin and $D=10$ Rosenbrock, remain visible after accounting for VAE pretraining and retraining costs. In settings where the final regrets are similar, however, adaptive retraining provides little benefit relative to its additional computational cost. Overall, adaptive retraining is most useful when the fixed latent representation is a substantial optimisation bottleneck, but its advantages are neither universal nor consistently strengthened by DML.

\FloatBarrier

\subsection{Comparisons Between BO--VAE and BO--SDR Algorithms}\label{subsection: Comparisons Between BO--VAE and BO--SDR Algorithms}

\begin{figure}[h]
    \centering
    \includegraphics[width=1\linewidth]{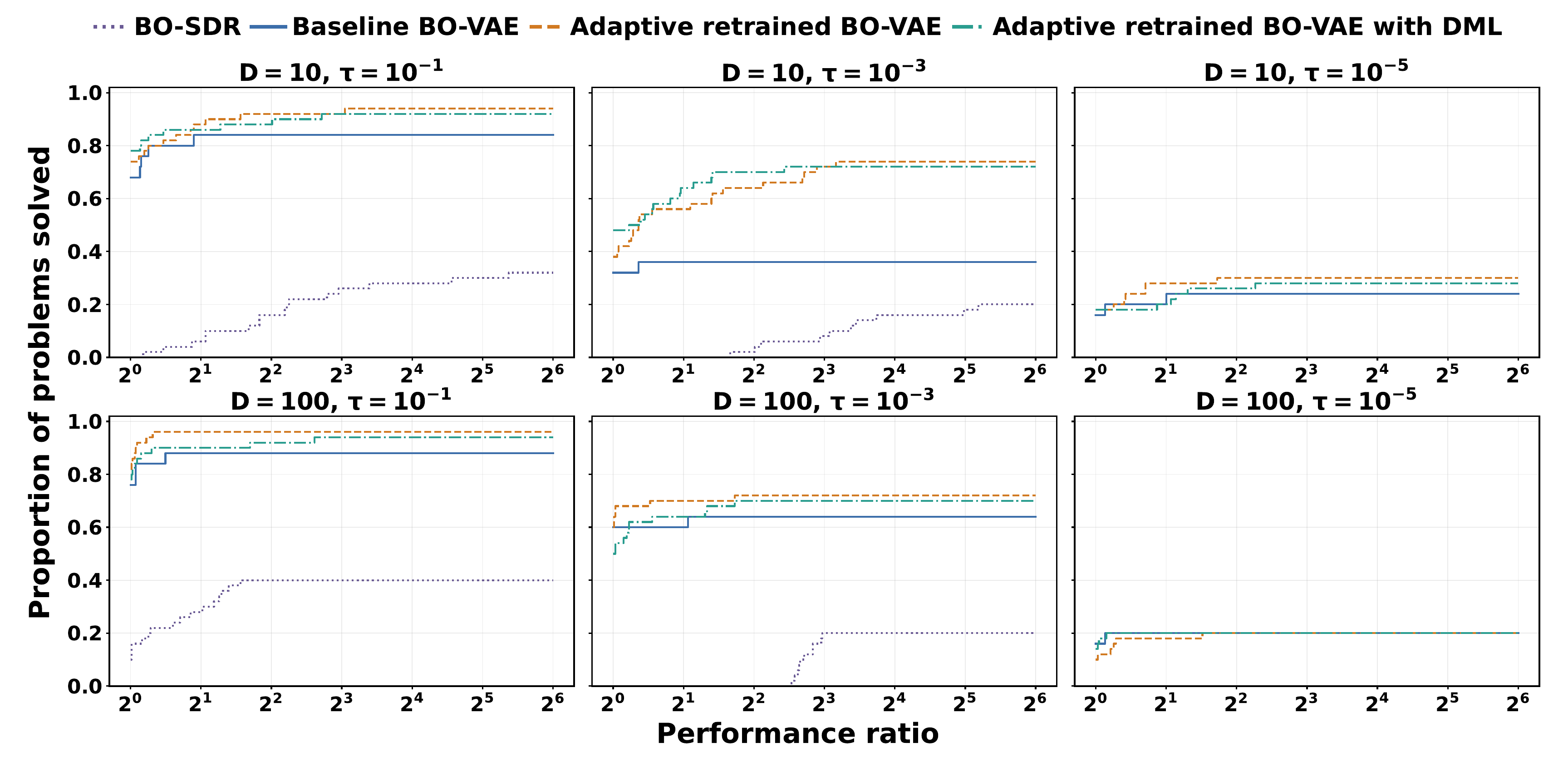}
    \caption{Performance profiles for BO--SDR and the three BO--VAE variants on the test problems in Table~\ref{tab:full-rank-test-problems} at $D=10$ and $D=100$, using tolerances $\tau\in\{10^{-1},10^{-3},10^{-5}\}$. The y--axis is the proportion of instances solved within performance ratio $\alpha$. The performance ratio is defined as the evaluation budget divided by the smallest evaluation budget required by any solver on the same problem, dimension, and seed. Thus, $\alpha=1$ corresponds to the best evaluation count.}
    \label{fig: Performance Profiles on full rank functions.}
\end{figure}

\begin{figure}[h]
    \centering
    \includegraphics[width=1.0\linewidth]{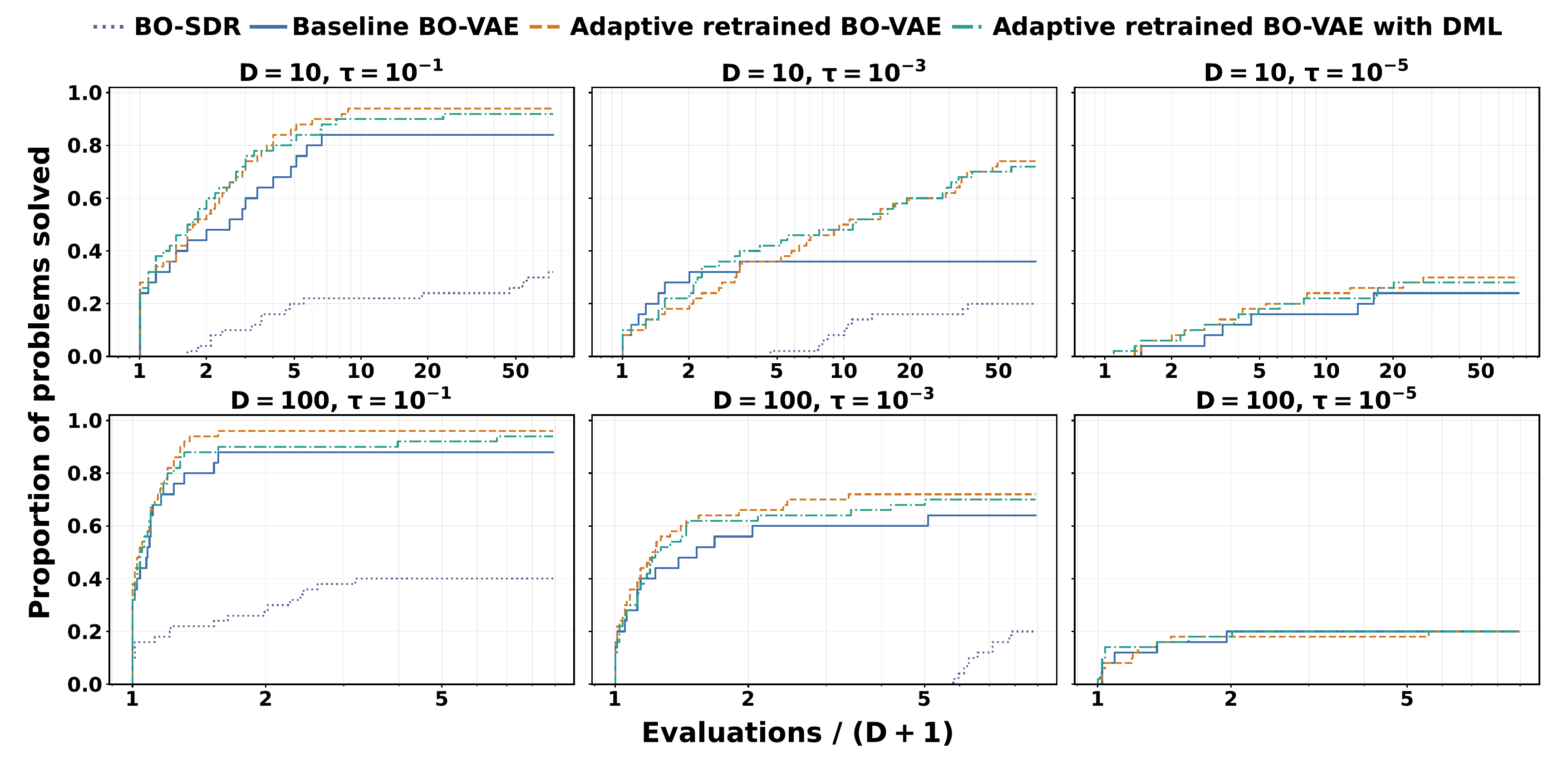}
    \caption{Data profiles for BO--SDR and the three BO--VAE variants on the test problems in Table~\ref{tab:full-rank-test-problems} at $D=10$ and $D=100$, using tolerances $\tau\in\{10^{-1},10^{-3},10^{-5}\}$. The y--axis is the proportion of problems solved within normalized evaluation budget $N/(D+1)$, where $N$ is the total evaluation count, including the initial design.}
    \label{fig: Data Profiles on full rank functions.}
\end{figure}

The preceding results describe performance on individual test problems. To assess robustness and evaluation efficiency across the complete benchmark suite, we additionally construct performance and data profiles as well as the final average percentage of problems solved accross the test function in Table~\ref{tab:full-rank-test-problems}, using ten seeds for $D=10$ and $D=100$. The BO--VAE methods use $d=6$ and $d=50$, respectively. The results are shown in Table~\ref{tab:full-rank-percentage-problems-solved}, Figure~\ref{fig: Performance Profiles on full rank functions.}, and Figure~\ref{fig: Data Profiles on full rank functions.}.

The results show that the BO-VAE variants generally outperform BO--SDR. The retrained BO--VAE variants generally attain the fastest coverage at the loose and intermediate tolerances, while solving most of the problems. The improvement is particularly evident for $D=10$ at $\tau=10^{-3}$, where both retraining variants solve substantially more instances than the fixed-VAE baseline. At $D=100$, their advantage at this tolerance is more modest. BO--SDR remains below the latent-space methods across the loose and intermediate tolerances.

At the strict tolerance $\tau=10^{-5}$, all curves plateau at relatively low levels, and the three BO--VAE variants become much less distinguishable, especially for $D=100$. This indicates that retraining improves moderate-accuracy performance but does not overcome the remaining representation and optimisation limitations at high accuracy. DML is competitive but provides no consistent advantage over retraining without DML.

\begin{table}[htbp]
  \centering
  \caption{Average percentage of problems solved at three tolerances.}
  \label{tab:full-rank-percentage-problems-solved}

  \setlength{\tabcolsep}{3pt}
  \renewcommand{\arraystretch}{1.1}

  \begin{tabularx}{\linewidth}{
    @{}
    >{\raggedright\arraybackslash}X
    *{6}{r}
    @{}
  }
    \toprule
    & \multicolumn{2}{c}{$\tau=10^{-1}$}
    & \multicolumn{2}{c}{$\tau=10^{-3}$}
    & \multicolumn{2}{c}{$\tau=10^{-5}$} \\
    \cmidrule(lr){2-3}
    \cmidrule(lr){4-5}
    \cmidrule(lr){6-7}
    Algorithm
    & $D=10$ & $D=100$
    & $D=10$ & $D=100$
    & $D=10$ & $D=100$ \\
    \midrule
    BO--SDR
    & 32\% & 40\%
    & 20\% & 20\%
    & 0\% & 0\% \\
    Baseline BO--VAE
    & 84\% & 88\%
    & 36\% & 64\%
    & 24\% & 20\% \\
    Adaptive retrained BO--VAE
    & 94\% & 96\%
    & 74\% & 72\%
    & 30\% & 20\% \\
    Adaptive retrained BO--VAE with DML
    & 92\% & 94\%
    & 72\% & 70\%
    & 28\% & 20\% \\
    \bottomrule
  \end{tabularx}
\end{table}

\FloatBarrier

\subsection{Comparison with EGORSE on Problems with Nonlinear Effective Mappings}
\label{sec:bovae-egorse-nonlinear}

EGORSE \cite{egorse} is designed to exploit low-dimensional linear structure and therefore provides a natural active-subspace baseline for high-dimensional optimisation. Such linear representations can be highly effective when the objective is well approximated by $F(\mathbf{x})=h(A\mathbf{x}),$ where $A$ is some (random) linear mapping. In practical problems, however, the effective coordinates may depend nonlinearly on the ambient variables, in which case a single linear projection may not adequately represent the objective-relevant geometry. We therefore compare the linear dimensionality reduction ability of EGORSE with the ability of the BO--VAE algorithms to learn nonlinear representations.

The nonlinear test problems are constructed from the base functions listed in Table~\ref{tab:egorse-base-functions}. Let $h:\mathbb{R}^{d_e}\rightarrow \mathbb{R}$ denote a base function and let an orthogonal matrix be partitioned
as $Q=[A^\top, C^\top]^\top$, where $A\in\mathbb{R}^{d_e\times D}$. Instead of the conventional linear construction $F_{\mathrm{lin}}(\mathbf{x}) =h(A\mathbf{x})$, we define
\[
  F_\alpha(\mathbf{x})
  =
  \widetilde h\!\left(\Phi_\alpha(\mathbf{x})\right),
  \qquad
  \Phi_\alpha(\mathbf{x})
  =
  A\mathbf{x}+\alpha G(C\mathbf{x}),
\]
where $\widetilde h$ includes the transformation to the domain of the corresponding base function and $G$ is a smooth, bounded, and centred nonlinear mapping satisfying $G(\mathbf{0})=\mathbf{0}$. When $\alpha=0$, this construction reduces to the conventional linear low-rank problem. For $\alpha>0$, the preimages of the effective coordinates become curved manifolds rather than affine subspaces. We use $\alpha=1.0$ and the ambient domain
$\mathcal{X}=[-1,1]^D$. The construction of $G$, preservation of the known minima, and geometric validation of the resulting problems are deferred to Appendix~\ref{app:curved-preimage-construction}.

We use $D\in\{10,100\}$ in the experiments. Branin has effective dimension $d_e=2$, whereas Ackley, Rosenbrock, and Rastrigin have $d_e=4$. For BO--VAE, the latent dimension is set to $d= d_e+1$. The VAE architectures are reported in Table~\ref{tab:curved-preimage-vae-family}. Each BO--VAE variant begins from the same compatible objective-blind checkpoint for a given $(D,d)$ configuration. The retraining variants propose a VAE update according to the same adaptive retraining policy in Section~\ref{sec:bovae-benchmark-performance}, while the baseline retains the initial VAE throughout. We compare these methods with EGORSE operating in dimensions $d_e$ and $2d_e$. Every method uses $D$ initial observations followed by $800$ enrichment evaluations.

\begin{table}[h]
\caption{VAEs used in the comparison against EGORSE. Bracketed lists give layer widths from the encoder input to the latent layer and from the latent layer to the decoder output.} 
\label{tab:curved-preimage-vae-family}
\centering
\setlength{\tabcolsep}{3pt}
\renewcommand{\arraystretch}{1.1}

\begin{tabularx}{\linewidth}{@{}l c c Y Y c c@{}}
\toprule
VAE ID & $D$ & $d$ & Encoder & Decoder
& \shortstack{Hidden\\activation}
& \shortstack{Decoder\\output} \\
\midrule

VAE--E10
& 10
& $d_e+1$
& $[10,\allowbreak 64,\allowbreak 32,\allowbreak d]$
& $[d,\allowbreak 32,\allowbreak 64,\allowbreak 10]$
& SiLU
& Scaled $\tanh$ \\

VAE--E100
& 100
& $d_e+1$
& $[100,\allowbreak 256,\allowbreak 128,\allowbreak d]$
& $[d,\allowbreak 128,\allowbreak 256,\allowbreak 100]$
& SiLU
& Scaled $\tanh$ \\

\bottomrule
\end{tabularx}
\end{table}

We note that the experiments should be interpreted as a comparison between learned nonlinear representations and a linear active-subspace baseline, examining whether VAE retraining can adapt the latent representation to curved effective structure. They are not intended to establish that BO--VAE generally dominates EGORSE, particularly on problems with genuinely linear structure.

\begin{figure}[h]
    \centering
    \includegraphics[width=\linewidth]{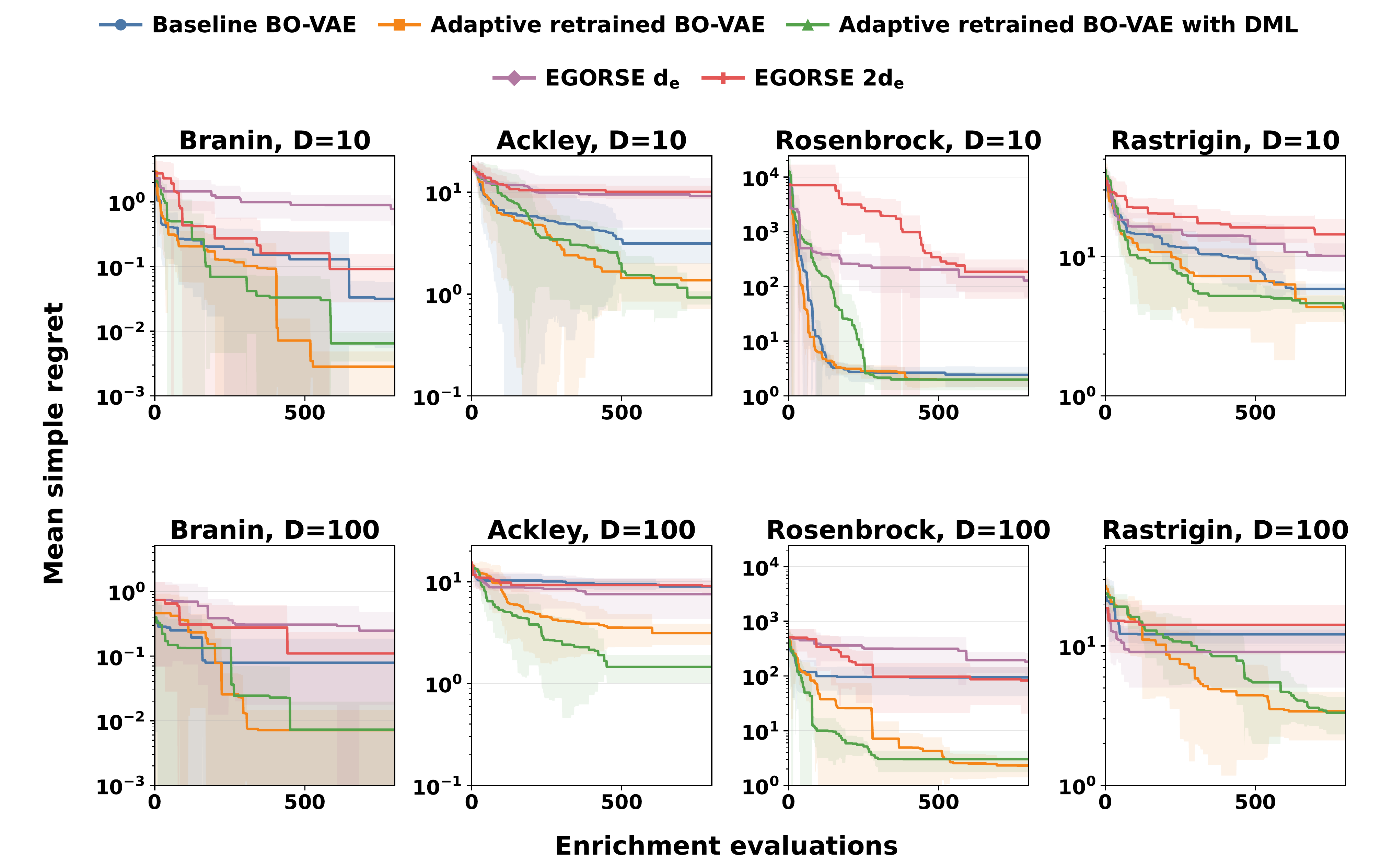}
    \caption{Simple regret plots for the three BO--VAE variants and the EGORSE solver on the test problems with nonlinear effective mappings. The upper and lower rows correspond to $D=10$ and $D=100$, respectively. Solid curves show the arithmetic mean, and shaded regions show one sample standard deviation on a logarithmic scale.}
    \label{fig:curved-preimage-anytime}
\end{figure}

Figure~\ref{fig:curved-preimage-anytime} presents the simple regret results. It shows that at least one retrained BO--VAE variant attains the lowest displayed final mean simple regret in all eight problem--dimension settings. The separation is especially clear for the $D=100$ Ackley, Rosenbrock, and Rastrigin problems. In several panels, the fixed-VAE and EGORSE curves plateau while the retraining variants continue to improve, which is consistent with retraining adapting the latent representation to the nonlinear effective structure. On $D=100$ Rastrigin, for example, EGORSE is competitive during the early iterations, but the retrained BO--VAE variants achieve lower regret later in the run.

The contribution of DML remains problem-dependent. It gives the lowest displayed final regret on $D=10$ Ackley, $D=100$ Rosenbrock, and both Rastrigin problems, whereas retraining without DML performs better on both Branin problems, $D=100$ Ackley, and $D=10$ Rosenbrock. Similarly, increasing the EGORSE search dimension from $d_e$ to $2d_e$ does not provide a consistent benefit, indicating that enlarging a linear subspace does not necessarily resolve the nonlinear representation mismatch. These results support the use of learned nonlinear representations in this problem class, but do not imply general superiority over EGORSE on linearly structured problems.

\FloatBarrier

\begin{figure}[h]
    \centering
    \includegraphics[width=\linewidth]{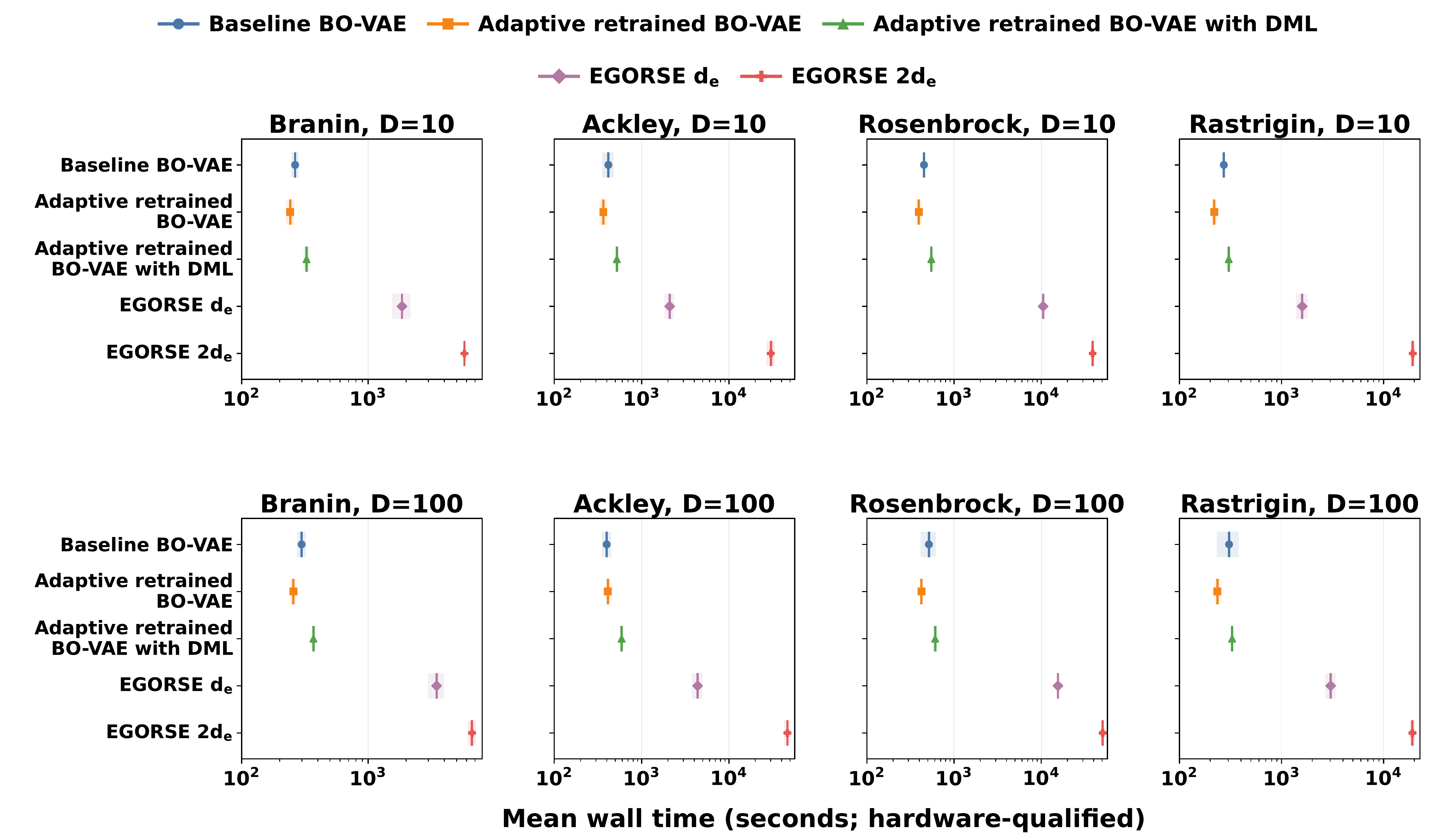}
    \caption{Mean optimisation-stage wall-clock time for the EGORSE comparison experiments. Markers show the mean and shaded bands show one sample standard deviation on a logarithmic scale, using the same displayed runs as Fig.~\ref{fig:curved-preimage-anytime}. The costs of VAE pretraining and the periodic VAE-retraining proposals are included. The BO--VAE methods were executed using one NVIDIA RTX~5080 CUDA worker, while the EGORSE methods were executed on CPU with one BLAS thread per worker and at most six concurrent workers, since it is only available on CPU. The timings are therefore provided for reference only and should not be interpreted as a hardware-normalised or algorithm-only speed comparison.}
    \label{fig:curved-preimage-runtime}
\end{figure}

Figure~\ref{fig:curved-preimage-runtime} reports the corresponding computation times. Within the BO--VAE family, retraining without DML has the lowest displayed mean runtime, while the DML variant incurs additional computation. Within the EGORSE family, operating in $2d_e$ dimensions is consistently more expensive than operating in $d_e$ dimensions. We explicitly warn that the methods were executed on different hardware and through different software backends.

\FloatBarrier





\section{Discussion and Conclusions}\label{Conclusion}

We investigated dimensionality reduction as a route to scalable BO and examined whether learned nonlinear representations can make BO more effective in high-dimensional ambient spaces. We developed three BO--VAE variants combining LSBO with SDR, periodic or adaptive VAE retraining, and deep metric learning. The main conclusion is that BO--VAE performance is governed by two distinct requirements, namely, the decoder must provide access to high-quality ambient solutions, and BO must be able to locate the corresponding latent points within its finite evaluation budget.

The fixed-representation analysis makes this finding explicit. Under the conditions of Theorem~\ref{thm:fixed-representation-bovae}, the expected ambient simple regret is bounded by a latent optimisation term of order $O(n^{-1/d})$ and the VAE-induced representation gap $\Delta_{\mathrm{rep}}$. The computational domain gap $\Delta_{\mathcal{X}}$ can also be included when the search domain does not contain an ambient global minimiser. Thus, increasing the BO budget can eliminate the latent optimisation error, but cannot remove a representation floor created by a fixed decoder. Under the stated regularity and local-mass assumptions, the representation gap can additionally be controlled by a data-dependent high-probability reconstruction certificate, highlighting that a small average reconstruction loss alone is insufficient. 

This theorem is deliberately a minimal guarantee for a controlled fixed-representation subcase. It assumes a fixed VAE and latent region, deterministic mean decoding, noiseless evaluations, a fixed GP prior, and exact EI optimisation. It is not a convergence theorem for the complete algorithms with SDR, VAE retraining, or DML. Those mechanisms change the search region, latent coordinates, decoder-supported set, or latent objective during the run and therefore require a stagewise nonstationary analysis together with additional reasonable assumptions controlling the stability and improvement of successive representations. The result should consequently be interpreted as a rigorous baseline that identifies the error terms which any successful adaptive representation method must control. 

The numerical representation-floor diagnostics support this interpretation, although they are not formal certificates. They reveal a search--support
trade-off: increasing \(d\) can improve decoder support but also makes latent BO more difficult. Moreover, a retrained VAE can have a worse final floor than the best representation encountered earlier, while the corresponding incumbent is retained. Representation quality and incumbent quality should therefore be monitored separately.

SDR is most beneficial in lower-dimensional latent spaces, while adaptive retraining and DML remain problem-dependent. An accepted VAE update does not necessarily improve optimisation, and DML is not uniformly superior to reconstruction-based retraining, motivating more objective-aware acceptance criteria.

Within the controlled benchmark suite, BO--VAE generally achieves greater coverage than BO--SDR, particularly at the loose and intermediate tolerances, although the advantage narrows at the strictest tolerance. On the nonlinear curved-preimage problems, increasing the EGORSE embedding dimension from \(d_e\) to \(2d_e\) does not consistently overcome the linear-representation mismatch. This supports nonlinear representation learning in the problem class considered, but does not imply general superiority over EGORSE on linearly structured objectives.

However, several limitations remain. First, the latent-dimension scans diagnose the consequences of different choices of $d$, but they do not provide an automatic rule for selecting $d$ when the effective dimension is unknown. The practical dimension-selection problem therefore remains unresolved. Second, the empirical study uses controlled analytical benchmarks with known global optima. While this is important for computing simple regret and estimating representation floors, it does not constitute validation on genuinely expensive real-world black-box problems. Third, the adopted retraining intervals, acceptance thresholds, and SDR settings are experimental design choices rather than generally optimal configurations.

The current theoretical analysis concerns the fixed-representation setting. A natural next step is to extend the analysis to the retraining setting. Another important future direction is to develop objective-aware retraining policies such as augmented VAE losses, retraining intervals, and acceptance criteria that explicitly control decoder support near promising solutions. Such mechanisms should preserve access to the best representation encountered so far without knowing the global optimum, rather than accepting updates only because they improve reconstruction or latent metric quality. Automatic or adaptive selection of the latent dimension, retraining interval, and SDR policy is another important direction. Finally, experiments on genuinely expensive scientific or engineering black-box problems are needed to determine whether the advantages observed on the controlled nonlinear benchmarks persist under real evaluation costs, constraints, noise, and unknown effective dimension.

\newpage

\begin{appendices}

\section{Test Sets}\label{appendix: test sets}
\subsection{Solver Comparison Test Function Set}
\begin{table}[h]
\caption{Benchmark test problems for BO solvers. Problems marked with $\ast$ have variable dimension; chosen settings follow \cite{Cartis_2021,cartis2021globaloptimizationusingrandom}.}
\label{tab:solver-test}
\centering
\begin{tabular}{@{}l l c l c@{}}
\toprule
\# & Function & Dim.\ $d$ & Domain & {Global min.} \\
\midrule
1 & Beale \cite{Ernesto2005} & 2 & $\mathbf{x}\in[-4.5,4.5]^2$   & 0 \\
2 & Rosenbrock$^{\ast}$ \cite{Surjanovic2013} & 3 & $\mathbf{x}\in[-5,10]^3$     & 0 \\
3 & Hartmann-3 \cite{Ernesto2005} & 3 & $\mathbf{x}\in[0,1]^3$        & -3.86278 \\
4 & Hartmann-6 \cite{Ernesto2005} & 6 & $\mathbf{x}\in[0,1]^6$        & -3.32237 \\
5 & Shekel-5 \cite{Surjanovic2013} & 4 & $\mathbf{x}\in[0,10]^4$       & -10.1532 \\
6 & Rastrigin$^{\ast}$ \cite{Surjanovic2013} & 5 & $\mathbf{x}\in[-5.12,5.12]^5$ & 0 \\
\botrule
\end{tabular}
\end{table}
\FloatBarrier
\subsection{High-dimensional Full-rank Test Set}
Table~\ref{tab:full-rank-test-problems} lists the scalable full-rank benchmark families used in the high-dimensional experiments. The symbol \(D\) denotes the ambient dimension of the instantiated problem. These benchmark functions can be defined for variable \(D\). In our experiments, \(D\) is chosen to match the input dimension of the VAE architecture used in the corresponding run. The tested values are \(D=10\) and \(D=100\), as specified in the table and in Table~\ref{tab:CN-reference VAE family}.

\begin{table}[h]
\centering
\caption{Benchmark high-dimensional full-rank test problems. The listed functions
are scalable in the ambient dimension \(D\).}
\label{tab:full-rank-test-problems}
\begin{tabular}{llllll}
\toprule
\# & Function & Scalable dim. & Domain & Global min. \\
\midrule
1 & Ackley & any \(D\geq 1\)  & \(\mathbf{x}\in[-30,30]^D\) & $0$ \\
2 & L\'evy & any \(D\geq 1\) & \(\mathbf{x}\in[-10,10]^D\) & $0$ \\
3 & Rosenbrock & any \(D\geq 2\)  & \(\mathbf{x}\in[-5,10]^D\) & $0$ \\
4 & Styblinski--Tang & any \(D\geq 1\)  & \(\mathbf{x}\in[-5,5]^D\) & $-39.16599\times D$\\
5 & Rastrigin & any \(D\geq 1\) & \(\mathbf{x}\in[-5.12,5.12]^D\) & $0$ \\
\bottomrule
\end{tabular}
\end{table}

\subsection{High-dimensional Low-rank Test Set}\label{Appendix High-dimensional Low-rank Test Set}

The base functions for the test problems with nonlinear effective mappings are listed in Table~\ref{tab:egorse-base-functions}.

\begin{table}[htbp]
  \centering
  \caption{Base functions used to construct the test problems with nonlinear effective mappings for comparison against EGORSE.}
  \label{tab:egorse-base-functions}
  \begin{tabular}{lccc}
    \toprule
    Function & Eff.\ dim.\ $d_e$ & Domain & Global min. \\
    \midrule
    Branin--Hoo \cite{Ernesto2005} & $2$ & $[-5,10]\times[0,15]$
      & $0.397887$ \\
    Ackley \cite{Ernesto2005} & $4$ & $[-30,30]^4$
      & $0$ \\
    Rosenbrock \cite{Surjanovic2013} & $4$ & $[-5, 10]^4$
      & $0$ \\
    Rastrigin \cite{Surjanovic2013} & $4$ & $[-5.12,5.12]^4$
      & $0$ \\
    \bottomrule
  \end{tabular}
\end{table}

\section{Construction and Validation of the Curved-Preimage Problems}
\label{app:curved-preimage-construction}

\paragraph{Nonlinear effective mapping}

Let $\mathcal{X}=[-1,1]^D$ denote the ambient search domain and let $d_e<D$ be the effective dimension of a base function $h$. We construct an orthogonal matrix
\[
Q=
\begin{bmatrix}
A\\
C
\end{bmatrix}
\in\mathbb{R}^{D\times D},
\]
where $A\in\mathbb{R}^{d_e\times D}$ and $C\in\mathbb{R}^{(D-d_e)\times D}$. Hence,
\[
AA^\top=I_{d_e},\qquad
CC^\top=I_{D-d_e},\qquad
AC^\top=0,
\]
and
\[
A^\top A+C^\top C=I_D.
\]
For any $\mathbf{x}\in\mathcal{X}$, define $\mathbf{a}=A\mathbf{x}$ and $\mathbf{c}=C\mathbf{x}$. The conventional
linear construction would evaluate the base function using only $\mathbf{a}$. To introduce a nonlinear effective structure, we instead define
\[
\Phi_{\alpha}(\mathbf{x})
=
A\mathbf{x}+\alpha G(C\mathbf{x}),
\]
where $\alpha\geq 0$ controls the strength of the nonlinearity.

The nonlinear function component $G:\mathbb{R}^{D-d_e}\rightarrow\mathbb{R}^{d_e}$ is constructed as
\[
G(\mathbf{c})
=
\begin{bmatrix}
0 &
\bar g_2(\mathbf{c}) &
\cdots &
\bar g_{d_e}(\mathbf{c})
\end{bmatrix}^{\top},
\]
with
\[
g_i(\mathbf{c})
=
\tanh\left(
\frac{1}{\sqrt{D-d_e}}
\sum_{k=1}^{D-d_e}
W_{ik}\sin(\pi c_k+\theta_{ik})
\right),
\qquad i=2,\ldots,d_e,
\]
and
\[
\bar g_i(\mathbf{c})
=
\frac{g_i(\mathbf{c})-g_i(\mathbf{0})}{2}.
\]
Here, $W_{ik}\in\{-1,+1\}$ are fixed Rademacher signs and $\theta_{ik}\sim\mathcal{U}[0,2\pi)$ are fixed phase shifts. The centring ensures that
\[
G(\mathbf{0})=\mathbf{0},
\]
while the factor $1/2$ limits the magnitude of the nonlinear perturbation. The matrices and other parameters are generated once using chosen seeds and then held fixed throughout all runs.

The resulting ambient-space problem is
\[
F_{\alpha}(\mathbf{x})
=
h\!\left(T_h\!\left(\Phi_{\alpha}(\mathbf{x})\right)\right),
\qquad \mathbf{x}\in[-1,1]^D,
\]
where $T_h$ maps the effective coordinates to the conventional domain of the corresponding base function in Table~\ref{tab:egorse-base-functions}. 

For a fixed effective coordinate $\mathbf{z}$, the associated preimage satisfies
\[
\mathbf{a}+\alpha G(\mathbf{c})=\mathbf{z}.
\]
Using $\mathbf{x}=A^\top\mathbf{a}+C^\top\mathbf{c}$, it can therefore be parametrised by
\[
\mathbf{x}(\mathbf{c})
=
A^\top\!\left\{\mathbf{z}-\alpha G(\mathbf{c})\right\}
+
C^\top\mathbf{c}.
\]
When $\alpha=0$, this expression reduces to the affine preimage $A^\top\mathbf{z}+C^\top\mathbf{c}$ associated with a conventional linear effective mapping. For $\alpha>0$, the dependence of the active coordinates on $\mathbf{c}$ generally produces a curved $(D-d_e)$-dimensional preimage within the ambient box.

\paragraph{Preservation of the known global minima}

The centring condition $G(\mathbf{0})=\mathbf{0}$ provides an explicit feasible ambient-space witness for every known minimiser of the base function. Let $\mathbf{z}^{\star}$ be an effective coordinate satisfying
\[
h\!\left(T_h(\mathbf{z}^{\star})\right)=h^\star,
\]
where $h^\star$ is the known global minimum of $h$. Define
\[
\mathbf{x}^{\star}=A^\top\mathbf{z}^{\star}.
\]
Orthogonality gives
\[
A\mathbf{x}^{\star}=\mathbf{z}^{\star},
\qquad
C\mathbf{x}^{\star}
=
CA^\top\mathbf{z}^{\star}
=
\mathbf{0}.
\]
It follows that, for any value of $\alpha$,
\[
\Phi_{\alpha}(\mathbf{x}^{\star})
=
\mathbf{z}^{\star}+\alpha G(\mathbf{0})
=
\mathbf{z}^{\star},
\]
and hence
\[
F_{\alpha}(\mathbf{x}^{\star})
=
h\!\left(T_h(\mathbf{z}^{\star})\right)
=
h^\star.
\]
Moreover,
\[
\|\mathbf{x}^{\star}\|_{\infty}
\leq
\|\mathbf{x}^{\star}\|_2
=
\|\mathbf{z}^{\star}\|_2.
\]
The effective minimisers used in the experiments satisfy $\|\mathbf{z}^{\star}\|_2\leq 0.5$, so that $\mathbf{x}^{\star}\in[-1,1]^D$. Because the lifted objective is an evaluation of the original base function, it cannot attain a value below $h^\star$. Thus, the construction preserves the known global minimum exactly.

We first specify the transformations from the base coordinates to the benchmark coordinates. Throughout, $\tanh$ acts componentwise, and $\overline{S}$ denotes the Euclidean closure of a set $S$. Let
\[
  \widehat{\mathbf z}
  =\left(\frac12,0,0,0\right)^\top,
\]
and let $\mathbf 1\in\mathbb R^4$ denote the all-ones vector.

For Ackley and Rastrigin, respectively, define
\[
  T_{\rm A}(\mathbf z)
  =30\tanh(\mathbf z-\widehat{\mathbf z}),
  \qquad
  T_{\rm Ra}(\mathbf z)
  =5.12\tanh(\mathbf z-\widehat{\mathbf z}),
  \qquad \mathbf z\in\mathbb R^4.
\]
Since $\tanh$ maps $\mathbb R$ bijectively onto $(-1,1)$, these transformations satisfy
\[
\begin{aligned}
  T_{\rm A}(\mathbb R^4)&=(-30,30)^4,
  &
  \overline{T_{\rm A}(\mathbb R^4)}&=[-30,30]^4,\\
  T_{\rm Ra}(\mathbb R^4)&=(-5.12,5.12)^4,
  &
  \overline{T_{\rm Ra}(\mathbb R^4)}&=[-5.12,5.12]^4.
\end{aligned}
\]

For Rosenbrock, let
\[
  \rho_R
  =\operatorname{artanh}\!\left(-\frac15\right),
  \qquad
  \mathbf b
  =\widehat{\mathbf z}-\rho_R\mathbf 1,
\]
and define
\[
  T_{\rm Ro}(\mathbf z)_j
  =-5+\frac{15}{2}\bigl[1+\tanh(z_j-b_j)\bigr],
  \qquad j=1,\ldots,4,\quad \mathbf z\in\mathbb R^4.
\]
Then
\[
  T_{\rm Ro}(\mathbb R^4)=(-5,10)^4,
  \qquad
  \overline{T_{\rm Ro}(\mathbb R^4)}=[-5,10]^4.
\]
The shifts are chosen so that
\[
  T_{\rm A}(\widehat{\mathbf z})
  =T_{\rm Ra}(\widehat{\mathbf z})
  =\mathbf 0,
  \qquad
  T_{\rm Ro}(\widehat{\mathbf z})=\mathbf 1.
\]
Thus, the transformed Ackley, Rastrigin, and Rosenbrock objectives share the global minimiser $\widehat{\mathbf z}$ in base coordinates.

For Branin--Hoo, define, for $\mathbf z\in\mathbb R^2$,
\[
\begin{aligned}
  T_{\rm B}(\mathbf z)_1
  &=-5+\frac{15}{2}\bigl[1+\tanh(4z_1)\bigr],\\
  T_{\rm B}(\mathbf z)_2
  &=\frac{15}{2}\bigl[1+\tanh(4z_2)\bigr].
\end{aligned}
\]
Its image and closure are
\[
  T_{\rm B}(\mathbb R^2)=(-5,10)\times(0,15),
  \qquad
  \overline{T_{\rm B}(\mathbb R^2)}
  =[-5,10]\times[0,15].
\]

Each transformation therefore maps the base-coordinate space onto the interior of the corresponding closed benchmark domain. Boundary points belong to the closure of the image but are not attained at finite base coordinates. The bound $\lVert\mathbf z^\star\rVert_2\leq 0.5$ is a direct numerical consequence of these specified transformations. For Branin--Hoo, the inverse used in the table is
\[
 z_1^\star=\frac14\operatorname{artanh}\!\left(
   2\frac{u_1^\star+5}{15}-1\right),
 \qquad
 z_2^\star=\frac14\operatorname{artanh}\!\left(
   2\frac{u_2^\star}{15}-1\right).
\]

\begin{table}[h]
\centering
\tiny
\setlength{\tabcolsep}{3pt}
\renewcommand{\arraystretch}{1.1}

\begin{tabularx}{\linewidth}{@{}l l l l @{}}
\toprule
Base function & Conventional minimiser $\mathbf u^\star$
 & Effective minimiser $\mathbf z^\star=T_h^{-1}(\mathbf u^\star)$
 & $\lVert\mathbf z^\star\rVert_2$ \\
\midrule

Branin--Hoo
 & $(-\pi,12.275)$
 & $(-0.244507794584447,\phantom{-}0.188137030695672)$
 & $0.308511918621525$ \\
Branin--Hoo
 & $(\pi,2.275)$
 & $(\phantom{-}0.021438821128525,-0.215198563754209)$
 & $0.216263831680786$ \\
Branin--Hoo
 & $(3\pi,2.475)$
 & $(\phantom{-}0.402743322143486,-0.202685781368784)$
 & $0.450869947435229$ \\
Ackley
 & $(0,0,0,0)$
 & $(0.5,0,0,0)$
 & $0.5$ \\
Rosenbrock
 & $(1,1,1,1)$
 & $(0.5,0,0,0)$
 & $0.5$ \\
Rastrigin
 & $(0,0,0,0)$
 & $(0.5,0,0,0)$
 & $0.5$ \\

\bottomrule
\end{tabularx}
\end{table}

The Ackley and Rastrigin entries follow immediately from their centred transformations.  For Rosenbrock, $z_j^\star-b_j=\rho_R$ and $\tanh(\rho_R)=-1/5$, so $T_{\rm Ro}(\mathbf z^\star)=(1,1,1,1)$.

\section{Proofs for the fixed-representation analysis}
\label{app:fixed-representation-proofs}

This appendix proves Theorem~\ref{thm:fixed-representation-bovae}.  The certificate sample size $N$ and the BO budget $n$ index different experiments and different probability spaces.  Throughout this section, the decoder and encoder are fixed, and probabilities not explicitly marked as EI probabilities refer to the independent sample $\mathcal{D}_{\mathrm{cert}}\sim\mu^{\otimes N}$.

\subsection{From empirical to uniform reconstruction}

Define the population reconstruction function
\begin{equation}
  h(\mathbf{x}):=
  \mathbb{E}_{\mathbf{Z}\sim q_{\phi^\star}(\cdot\mid\mathbf{x})}
  \bigl[\|\mathbf{x}-r_{\theta^\star}(\mathbf{Z})\|\bigr],
  \qquad \mathbf{x}\in\mathcal{M}.
  \label{eq:population-reconstruction-function}
\end{equation}

\begin{lemma}[PAC--Bayes population certificate]
\label{lem:pac-bayes-population-certificate}
Under Assumption~\ref{ass:fixed-representation-certificate}, with probability
at least $1-\eta$,
\begin{equation}
  \int_{\mathcal{M}}h(\mathbf{x})\,\mu(\mathrm d\mathbf{x})
  \leq \overline\varepsilon_N^{\mathrm{PB}}.
  \label{eq:population-reconstruction-bound}
\end{equation}
\end{lemma}

\begin{proof}
For $\mathbf{x},\mathbf{y}\in\mathcal{X}$, the closed form of the Wasserstein-$2$ distance between diagonal Gaussian measures gives
\begin{align}
  W_1\!\left(q_{\phi^\star}(\cdot\mid\mathbf{x}),
             q_{\phi^\star}(\cdot\mid\mathbf{y})\right)
  &\leq
  W_2\!\left(q_{\phi^\star}(\cdot\mid\mathbf{x}),
             q_{\phi^\star}(\cdot\mid\mathbf{y})\right)
  \notag\\
  &=\left(
      \|\mathbf{m}_{\phi^\star}(\mathbf{x})
        -\mathbf{m}_{\phi^\star}(\mathbf{y})\|^2
      +\|\mathbf{s}_{\phi^\star}(\mathbf{x})
        -\mathbf{s}_{\phi^\star}(\mathbf{y})\|^2
    \right)^{1/2}
  \notag\\
  &\leq K_\phi\|\mathbf{x}-\mathbf{y}\|.
  \label{eq:encoder-wasserstein-lipschitz}
\end{align}
For fixed $\mathbf{x}$, the reverse triangle inequality and decoder Lipschitzness imply
\begin{equation}
  \left|
    \|\mathbf{x}-r_{\theta^\star}(\mathbf{z})\|
    -\|\mathbf{x}-r_{\theta^\star}(\mathbf{z}')\|
  \right|
  \leq K_\theta\|\mathbf{z}-\mathbf{z}'\|.
  \label{eq:reconstruction-loss-lipschitz}
\end{equation}
These are precisely the conditional-posterior regularity conditions used by Proposition~4.1 of \cite{mbacke2023statisticalguaranteesvariationalautoencoders}.  Applying their Theorem~4.3 with instance-space diameter $D_{\mathcal{X}}$, loss $\|\mathbf{x}-r_{\theta^\star}(\mathbf{z})\|$, prior $p=\mathcal{N}(\mathbf{0},\mathbf{I}_d)$, sample size $N$, and confidence parameter $\eta$ yields
\begin{equation}
  \int_{\mathcal{M}}h(\mathbf{x})\,\mu(\mathrm d\mathbf{x})
  \leq \varepsilon_N^{\mathrm{PB}}.
  \label{eq:mbacke-specialization}
\end{equation}
Because $\mathbf{x}$ and $r_{\theta^\star}(\mathbf{z})$ lie in $\mathcal{X}$, their distance is at most $D_{\mathcal{X}}$, so the left-hand side is also at most $D_{\mathcal{X}}$.  Taking the smaller of these two valid bounds proves \eqref{eq:population-reconstruction-bound}.
\end{proof}

The independence condition on $\mathcal{D}_{\mathrm{cert}}$ is essential in this idealised case.  The cited PAC--Bayes result is uniform over the encoder posterior for a fixed decoder, but it is not uniform over decoders trained on the certificate sample.  Requiring the entire VAE to be fixed before drawing $\mathcal{D}_{\mathrm{cert}}$ is a simple sufficient condition.

\begin{lemma}[Average-to-uniform conversion]
\label{lem:average-to-uniform}
Suppose \eqref{eq:local-mass-condition} holds.  If $u:\mathcal{M}\to[0,\infty)$ is $K_u$-Lipschitz and $\int_{\mathcal{M}}u\,\mathrm d\mu\leq\varepsilon$, then 
\begin{equation}
  \sup_{\mathbf{x}\in\mathcal{M}}u(\mathbf{x})
  \leq
  \inf_{0<\rho\leq\rho_0}
  \left[K_u\rho+\frac{\varepsilon}{c_\mu\rho^m}\right].
  \label{eq:average-to-uniform-bound}
\end{equation}
If $\varepsilon>0$, $K_u>0$, and
\begin{equation}
  \rho_\star
  :=\left(\frac{m\varepsilon}{c_\mu K_u}\right)^{1/(m+1)}
  \leq\rho_0,
  \label{eq:optimal-local-radius}
\end{equation}
then the right-hand side of \eqref{eq:average-to-uniform-bound} equals
\begin{equation}
  C_m K_u^{m/(m+1)}
  \left(\frac{\varepsilon}{c_\mu}\right)^{1/(m+1)},
  \qquad
  C_m:=(m+1)m^{-m/(m+1)}.
  \label{eq:average-to-uniform-closed-form}
\end{equation}
\end{lemma}

\begin{proof}
Fix an arbitrary $\mathbf{x}\in\mathcal{M}$ and $0<\rho\leq\rho_0$.  For every $\mathbf{y}\in B(\mathbf{x},\rho)\cap\mathcal{M}$, Lipschitz continuity gives $u(\mathbf{y})\geq u(\mathbf{x})-K_u\rho$.  Since $u$ is nonnegative,
\begin{align}
  \varepsilon
  &\geq \int_{B(\mathbf{x},\rho)\cap\mathcal{M}}
    u(\mathbf{y})\,\mu(\mathrm d\mathbf{y})
  \notag\\
  &\geq
  \mu\bigl(B(\mathbf{x},\rho)\cap\mathcal{M}\bigr)
  \bigl[u(\mathbf{x})-K_u\rho\bigr]_+
  \notag\\
  &\geq c_\mu\rho^m\bigl[u(\mathbf{x})-K_u\rho\bigr]_+.
  \label{eq:average-uniform-proof-core}
\end{align}
In either case, whether the positive part vanishes or not, $u(\mathbf{x})\leq K_u\rho+\varepsilon/(c_\mu\rho^m)$.  The bound is independent of $\mathbf{x}$. Thus, taking the supremum over $\mathbf{x}$ and then the infimum over $\rho$ proves \eqref{eq:average-to-uniform-bound}.  Differentiating the expression in brackets gives \eqref{eq:optimal-local-radius} and substitution yields \eqref{eq:average-to-uniform-closed-form}.
\end{proof}

\begin{lemma}[Uniform reconstruction within the BO region]
\label{lem:uniform-reconstruction-in-region}
On the event in Lemma~\ref{lem:pac-bayes-population-certificate},
\begin{equation}
  \sup_{\mathbf{x}\in\mathcal{M}}
  \mathbb{E}\!\left[
    \|\mathbf{x}-r_{\theta^\star}(\mathbf{Z})\|
    \mathrel{\big|} \mathbf{Z}\in\mathcal{R}_0,\ \mathbf{X}=\mathbf{x}
  \right]
  \leq \varepsilon_N^{\mathrm{cert}}.
  \label{eq:conditional-uniform-reconstruction}
\end{equation}
\end{lemma}

\begin{proof}
We first show that $h$ is Lipschitz.  For $\mathbf{x},\mathbf{y}\in\mathcal{M}\subseteq\mathcal{X}$, add and subtract the integral of the $\mathbf{x}$-loss under the $\mathbf{y}$-encoder law:
\begin{align*}
\bigl|h(\mathbf{x})-h(\mathbf{y})\bigr|
&=
\Bigg|
\int
  \bigl\|\mathbf{x}-r_{\theta^\star}(\mathbf{z}_x)\bigr\|
  q_{\phi^\star}
    (\mathrm{d}\mathbf{z}_x\mid\mathbf{x})
-
\int
  \bigl\|\mathbf{x}-r_{\theta^\star}(\mathbf{z}_y)\bigr\|
  q_{\phi^\star}
    (\mathrm{d}\mathbf{z}_y\mid\mathbf{y})
\\
&\qquad
+
\int
  \left(
    \bigl\|\mathbf{x}-r_{\theta^\star}(\mathbf{z}_y)\bigr\|
    -
    \bigl\|\mathbf{y}-r_{\theta^\star}(\mathbf{z}_y)\bigr\|
  \right)
  q_{\phi^\star}
    (\mathrm{d}\mathbf{z}_y\mid\mathbf{y})
\Bigg|
\\
&\leq
\left|
\int
  \bigl\|\mathbf{x}-r_{\theta^\star}(\mathbf{z}_x)\bigr\|
  q_{\phi^\star}
    (\mathrm{d}\mathbf{z}_x\mid\mathbf{x})
-
\int
  \bigl\|\mathbf{x}-r_{\theta^\star}(\mathbf{z}_y)\bigr\|
  q_{\phi^\star}
    (\mathrm{d}\mathbf{z}_y\mid\mathbf{y})
\right|
\\
&\qquad
+
\left|
\int
  \left(
    \bigl\|\mathbf{x}-r_{\theta^\star}(\mathbf{z}_y)\bigr\|
    -
    \bigl\|\mathbf{y}-r_{\theta^\star}(\mathbf{z}_y)\bigr\|
  \right)
  q_{\phi^\star}
    (\mathrm{d}\mathbf{z}_y\mid\mathbf{y})
\right|.
\end{align*}

For the first term, define
\[
  \ell_{\mathbf{x}}(\mathbf{z})
  \coloneqq
  \bigl\|\mathbf{x}-r_{\theta^\star}(\mathbf{z})\bigr\|.
\]
We first show that $\ell_{\mathbf{x}}$ is $K_\theta$--Lipschitz. For any $\mathbf{z}_x,\mathbf{z}_y\in\mathcal{Z}$,
\begin{align*}
\bigl|
  \ell_{\mathbf{x}}(\mathbf{z}_x)
  -
  \ell_{\mathbf{x}}(\mathbf{z}_y)
\bigr|
&=
\left|
  \bigl\|\mathbf{x}-r_{\theta^\star}(\mathbf{z}_x)\bigr\|
  -
  \bigl\|\mathbf{x}-r_{\theta^\star}(\mathbf{z}_y)\bigr\|
\right|
\\
&\leq
\bigl\|
  r_{\theta^\star}(\mathbf{z}_x)
  -
  r_{\theta^\star}(\mathbf{z}_y)
\bigr\|
\\
&\leq
K_\theta\bigl\|\mathbf{z}_x-\mathbf{z}_y\bigr\|.
\end{align*}

The first inequality follows from the reverse triangle inequality, and the second follows from Assumption~\ref{ass:fixed-representation-certificate}. Hence, when $K_\theta>0$, the function
\[
  \psi_{\mathbf{x}}(\mathbf{z})
  \coloneqq
  \frac{\ell_{\mathbf{x}}(\mathbf{z})}{K_\theta}
\]
is $1$--Lipschitz. If $K_\theta=0$, the following bound is immediate. The Kantorovich--Rubinstein inequality therefore gives
\begin{align*}
&
\left|
\int
  \ell_{\mathbf{x}}(\mathbf{z}_x)
  q_{\phi^\star}
    (\mathrm{d}\mathbf{z}_x\mid\mathbf{x})
-
\int
  \ell_{\mathbf{x}}(\mathbf{z}_y)
  q_{\phi^\star}
    (\mathrm{d}\mathbf{z}_y\mid\mathbf{y})
\right|
\\
&\qquad\leq
K_\theta
W_1\!\left(
  q_{\phi^\star}(\cdot\mid\mathbf{x}),
  q_{\phi^\star}(\cdot\mid\mathbf{y})
\right)
\\
&\qquad\leq
K_\theta K_\varphi
\bigl\|\mathbf{x}-\mathbf{y}\bigr\|,
\end{align*}
where the second inequality follows from the encoder regularity assumption.

For the second term, the reverse triangle inequality yields
\begin{align*}
&
\left|
\int
  \left(
    \bigl\|\mathbf{x}-r_{\theta^\star}(\mathbf{z}_y)\bigr\|
    -
    \bigl\|\mathbf{y}-r_{\theta^\star}(\mathbf{z}_y)\bigr\|
  \right)
  q_{\phi^\star}
    (\mathrm{d}\mathbf{z}_y\mid\mathbf{y})
\right|
\\
&\qquad\leq
\int
\left|
  \bigl\|\mathbf{x}-r_{\theta^\star}(\mathbf{z}_y)\bigr\|
  -
  \bigl\|\mathbf{y}-r_{\theta^\star}(\mathbf{z}_y)\bigr\|
\right|
q_{\phi^\star}
  (\mathrm{d}\mathbf{z}_y\mid\mathbf{y})
\\
&\qquad\leq
\int
  \bigl\|\mathbf{x}-\mathbf{y}\bigr\|
  q_{\phi^\star}
    (\mathrm{d}\mathbf{z}_y\mid\mathbf{y})
=
\bigl\|\mathbf{x}-\mathbf{y}\bigr\|.
\end{align*}

Combining the two estimates gives
\begin{equation*}
  \bigl|h(\mathbf{x})-h(\mathbf{y})\bigr|
  \leq
  \bigl(1+K_\theta K_\varphi\bigr)
  \bigl\|\mathbf{x}-\mathbf{y}\bigr\|.
  \label{eq:h-lipschitz}
\end{equation*}

Lemma~\ref{lem:average-to-uniform}, applied with $K_u=1+K_\phi K_\theta$ and $\varepsilon=\overline\varepsilon_N^{\mathrm{PB}}$, therefore yields
\begin{equation}
  \sup_{\mathbf{x}\in\mathcal{M}}h(\mathbf{x})
  \leq
  \inf_{0<\rho\leq\rho_0}
  \left[
    (1+K_\phi K_\theta)\rho
    +\frac{\overline\varepsilon_N^{\mathrm{PB}}}
           {c_\mu\rho^m}
  \right].
  \label{eq:h-uniform-bound}
\end{equation}
For a Borel set $\mathcal{A}\subseteq\mathcal{R}_0$, define the conditioned encoder law
\begin{equation}
  q_{\phi^\star}^{\mathcal{R}_0}(\mathcal{A}\mid\mathbf{x})
  :=\frac{q_{\phi^\star}(\mathcal{A}\mid\mathbf{x})}
          {q_{\phi^\star}(\mathcal{R}_0\mid\mathbf{x})}.
  \label{eq:conditioned-encoder}
\end{equation}
Nonnegativity of the reconstruction loss and \eqref{eq:latent-region-mass} imply
\begin{align}
  \mathbb{E}_{\mathbf{Z}\sim
    q_{\phi^\star}^{\mathcal{R}_0}(\cdot\mid\mathbf{x})}
    \bigl[\|\mathbf{x}-r_{\theta^\star}(\mathbf{Z})\|\bigr]
  &=
  \frac{
    \mathbb{E}_{\mathbf{Z}\sim q_{\phi^\star}(\cdot\mid\mathbf{x})}
      \bigl[\|\mathbf{x}-r_{\theta^\star}(\mathbf{Z})\|
        \mathbf{1}_{\{\mathbf{Z}\in\mathcal{R}_0\}}\bigr]
  }{q_{\phi^\star}(\mathcal{R}_0\mid\mathbf{x})}
  \notag\\
  &\leq \frac{h(\mathbf{x})}{\beta_0}.
  \label{eq:conditioning-penalty}
\end{align}
This expectation is also at most $D_{\mathcal{X}}$.  Combining these observations with \eqref{eq:h-uniform-bound} gives exactly \eqref{eq:uniform-reconstruction-certificate} and proves the claim.
\end{proof}

\subsection{Representation gap}

\begin{proposition}[Certified representation gap]
\label{prop:certified-representation-gap}
On the event in Lemma~\ref{lem:pac-bayes-population-certificate},
\begin{equation}
  0\leq\Delta_{\mathrm{rep}}
  \leq L_f\bigl(\delta_{\mathcal{M}}
                 +\varepsilon_N^{\mathrm{cert}}\bigr).
  \label{eq:certified-representation-gap}
\end{equation}
\end{proposition}

\begin{proof}
The lower bound follows from
$r_{\theta^\star}(\mathcal{R}_0)\subseteq\mathcal{X}$.  Compactness of $\mathcal{M}$ and $\mathcal{S}_{\mathcal{X}}^\star$ gives $\mathbf{x}_{\mathcal{M}}\in\mathcal{M}$ and $\mathbf{x}^\star\in\mathcal{S}_{\mathcal{X}}^\star$ satisfying $\|\mathbf{x}_{\mathcal{M}}-\mathbf{x}^\star\|=\delta_{\mathcal{M}}$.  An infimum is no larger than an average, so Lemma~\ref{lem:uniform-reconstruction-in-region} and Lipschitz continuity of $f$ give
\begin{align}
  \Delta_{\mathrm{rep}}
  &=g^\star-f_{\mathcal{X}}^\star
  \notag\\
  &\leq
  \int_{\mathcal{R}_0}
    \left[f\bigl(r_{\theta^\star}(\mathbf{z})\bigr)
      -f(\mathbf{x}^\star)\right]
    q_{\phi^\star}^{\mathcal{R}_0}
      (\mathrm d\mathbf{z}\mid\mathbf{x}_{\mathcal{M}})
  \notag\\
  &\leq L_f
  \int_{\mathcal{R}_0}
    \|r_{\theta^\star}(\mathbf{z})-\mathbf{x}^\star\|
    q_{\phi^\star}^{\mathcal{R}_0}
      (\mathrm d\mathbf{z}\mid\mathbf{x}_{\mathcal{M}})
  \notag\\
  &\leq L_f\left(
    \int_{\mathcal{R}_0}
      \|r_{\theta^\star}(\mathbf{z})-\mathbf{x}_{\mathcal{M}}\|
      q_{\phi^\star}^{\mathcal{R}_0}
        (\mathrm d\mathbf{z}\mid\mathbf{x}_{\mathcal{M}})
    +\|\mathbf{x}_{\mathcal{M}}-\mathbf{x}^\star\|
  \right)
  \notag\\
  &\leq L_f\bigl(\varepsilon_N^{\mathrm{cert}}
                  +\delta_{\mathcal{M}}\bigr).
  \label{eq:representation-gap-proof}
\end{align}
\end{proof}

This step explains why average reconstruction error is not, by itself, an objective-specific reachability guarantee.  The lower-mass condition converts the population average into control at the unknown point of $\mathcal{M}$ nearest the optimizer, while $\delta_{\mathcal{M}}$ records whether the population support approaches the optimizer at all.

\subsection{Latent EI Rate and Proof of the Main Theorem}

\begin{proof}[Proof of Theorem~\ref{thm:fixed-representation-bovae}]
Under Assumption~\ref{ass:fixed-prior-ei}, the Mat\'ern-$5/2$ kernel has smoothness parameter $\nu=5/2>1$, and $g$ lies in its RKHS ball of radius
$B_g$.  Theorem~2 of \cite{bull2011ego} therefore gives constants $C_{\mathrm{EI}}<\infty$ and $n_0\in\mathbb{N}$ such that
\begin{equation}
  \mathbb{E}_{\mathrm{EI}}\!\left[
    \min_{1\leq i\leq n}g(\mathbf{z}_i)-g^\star
  \right]
  \leq C_{\mathrm{EI}}n^{-1/d},
  \qquad n\geq n_0.
  \label{eq:bull-ei-specialization}
\end{equation}
The constant can depend on the fixed domain, kernel, GP scale and length-scales, and $B_g$, but not on $n$.

On the certificate event of probability at least $1-\eta$, Proposition~\ref{prop:certified-representation-gap} bounds $\Delta_{\mathrm{rep}}$.  Taking $\mathbb{E}_{\mathrm{EI}}$ in the exact decomposition \eqref{eq:fixed-regret-decomposition}, applying \eqref{eq:bull-ei-specialization}, and then applying \eqref{eq:certified-representation-gap} proves both inequalities in \eqref{eq:fixed-representation-main-bound}.
\end{proof}

\begin{proof}[Proof of Corollary~\ref{cor:fixed-representation-floor}]
The latent error in \eqref{eq:fixed-regret-decomposition} is nonnegative and, by \eqref{eq:bull-ei-specialization}, its expectation converges to zero. The remaining two terms are deterministic with respect to the EI randomness, which proves \eqref{eq:fixed-representation-asymptote}.
\end{proof}

\begin{proof}[Proof of Corollary~\ref{cor:fixed-exact-representation}]
Set $\Delta_{\mathcal{X}}=\Delta_{\mathrm{rep}}=0$ in the first inequality of \eqref{eq:fixed-representation-main-bound}.  The sufficient certificate condition follows from \eqref{eq:certified-representation-gap}.
\end{proof}

\section{Standard BO Algorithm with SDR}\label{appendix: BO-SDR}

Here we present Bayesian Optimisation with SDR in the ambient space, followed by a brief discussion about SDR. A related description appears in the appendix of our work \cite{long2024dimensionalityreductiontechniquesglobal}. The notation \(\mathcal D_k\) denotes the labelled BO dataset after \(k\) evaluations, and \(\mathcal R_k\subseteq\mathcal X\) denotes the current region of interest. The acquisition function is maximised over \(\mathcal R_k\), not over the full initial search box \(\mathcal X\).

\begin{algorithm}[h]
\caption{Bayesian Optimisation with Sequential Domain Reduction}
\label{BOSDR}
\begin{algorithmic}[1]
\Require Initial labelled dataset $\mathcal D_0\coloneq\{(\mathbf x_i,y_i)\}_{i=1}^{N_0}$, where $y_i\coloneq f(\mathbf x_i)$; budget $B$; initial search box $\mathcal X\subseteq\mathbb R^D$; acquisition rule $u(\cdot\mid\cdot)$ \textrm{(EI)}; SDR parameters $(\gamma_o,\gamma_p,\eta_{\mathrm{SDR}})$; minimum RoI side length $t_{\min}$; SDR update period $K_{\mathrm{SDR}}$.
\Ensure Best objective value $f_{\min}$ found.

\State Initialise the region of interest: $\mathcal R_0 \leftarrow \mathcal X$.

\For{$k=0,1,\ldots,B-1$}
    \State Fit GP surrogate $h_k:\mathcal X\to\mathbb R$ on $\mathcal D_k$.

    \State $u_k(\mathbf x)\coloneq u(\mathbf x\mid\mathcal D_k)$. \Comment{The current acquisition function.}

    \State $\hat{\mathbf x}_{k+1} \leftarrow \operatorname*{arg\,max}_{\mathbf x\in\mathcal R_k} u_k(\mathbf x)$. \Comment{Select the next point within $\mathcal{R}_k.$}

    \State $y_{k+1}\leftarrow f(\hat{\mathbf x}_{k+1})$.

    \State $\mathcal D_{k+1} \leftarrow \mathcal D_k\cup\{(\hat{\mathbf x}_{k+1},y_{k+1})\}$. 

    \If{$(k+1)\bmod K_{\mathrm{SDR}}=0$} \Comment{Update ROI;$ \mathcal{R}_k \subseteq \mathcal R_0.$}
        \State $\mathcal R_{k+1} \leftarrow \mathrm{SDR\_Update}\bigl(\mathcal R_k,\mathcal D_{k+1};\gamma_o,\gamma_p,\eta_{\mathrm{SDR}},t_{\min},\mathcal X\bigr)$.
    \Else
        \State $\mathcal R_{k+1}\leftarrow\mathcal R_k$.
    \EndIf
\EndFor

\State \Return $f_{\min} \leftarrow \min\{\,y:(\mathbf x,y)\in\mathcal D_B\,\}$.
\end{algorithmic}
\end{algorithm}

\paragraph{Sequential domain reduction in brief.}
To formally introduce SDR \cite{SDR}, let the initial search box be
\[
    \mathcal X
    =
    \prod_{i=1}^D [a_i,b_i],
\]
and let the current region of interest be the axis-aligned box
\[
    \mathcal R_k
    =
    \prod_{i=1}^D [\ell_{k,i},u_{k,i}]
    \subseteq \mathcal X,
    \qquad
    r_{k,i}\coloneq u_{k,i}-\ell_{k,i}.
\]
The incumbent after \(k\) evaluations is denoted by
\[
    \mathbf x_k^{\mathrm{inc}}
    \in
    \operatorname*{arg\,min}
    \{\,y:(\mathbf x,y)\in\mathcal D_k\,\}.
\]
If there are several incumbents, any one of them may be selected.

At an SDR update, the region is re-centred around the current incumbent and its side lengths are modified coordinatewise. For \(k\geq 1\), define the scaled incumbent movement in coordinate \(i\) by
\[
    s_{k,i}
    \coloneq
    \frac{2\bigl(x_{k,i}^{\mathrm{inc}}-x_{k-1,i}^{\mathrm{inc}}\bigr)}
         {r_{k-1,i}} .
\]
For \(k\geq 2\), define the oscillation indicator
\[
    c_{k,i}\coloneq s_{k,i}s_{k-1,i},
    \qquad
    \widehat c_{k,i}
    \coloneq
    \operatorname{sgn}(c_{k,i})\sqrt{|c_{k,i}|}.
\]
For the first SDR update, where \(s_{k-1,i}\) is not yet available, we set
\(\widehat c_{k,i}=0\).

The coordinatewise panning--damping parameter is
\[
    \gamma_{k,i}
    \coloneq
    \frac12
    \left[
        \gamma_p(1+\widehat c_{k,i})
        +
        \gamma_o(1-\widehat c_{k,i})
    \right],
\]
where \(\gamma_p\approx 1\) corresponds to panning and \(\gamma_o\in[0.5,0.7]\) gives stronger shrinkage when the incumbent movement oscillates. The side-length contraction factor is
\[
    \lambda_{k,i}
    \coloneq
    \eta_{\mathrm{SDR}}
    +
    |s_{k,i}|
    \bigl(\gamma_{k,i}-\eta_{\mathrm{SDR}}\bigr),
    \qquad
    \eta_{\mathrm{SDR}}\in[0.5,1).
\]
The proposed new side length is clipped to remain feasible and not fall below the minimum RoI size:
\[
    \widetilde r_{k,i}
    \coloneq
    \min\{\,b_i-a_i,\,
    \max\{t_{\min},\,\lambda_{k,i}r_{k-1,i}\}\,\}.
\]
Before trimming to \(\mathcal X\), the candidate bounds are
\[
    \widetilde \ell_{k,i}
    \coloneq
    x_{k,i}^{\mathrm{inc}}-\frac12\widetilde r_{k,i},
    \qquad
    \widetilde u_{k,i}
    \coloneq
    x_{k,i}^{\mathrm{inc}}+\frac12\widetilde r_{k,i}.
\]
The bounds are then shifted, if necessary, to keep the RoI inside the initial search box:
\[
    \ell_{k,i}
    \coloneq
    \min\{
        \max\{\widetilde \ell_{k,i},a_i\},
        b_i-\widetilde r_{k,i}
    \},
    \qquad
    u_{k,i}
    \coloneq
    \ell_{k,i}+\widetilde r_{k,i}.
\]
Thus,
\[
    \mathcal R_k
    =
    \prod_{i=1}^D [\ell_{k,i},u_{k,i}]
    \subseteq \mathcal X.
\]

The construction guarantees that every region remains feasible,
\[
    \mathcal{R}_k
    =\prod_{i=1}^{D}[\ell_{k,i},u_{k,i}]
    \subseteq\mathcal{X},
\]
and that its side lengths satisfy the prescribed minimum- and maximum-size constraints. It does not, however, impose the set-inclusion condition
$\mathcal{R}_k\subseteq\mathcal{R}_{k-1}$. Recentring the box around a new incumbent can violate this inclusion even when every side length decreases. Furthermore, under sustained incumbent movement, the scaling factor $\lambda_{k,i}$ may permit temporary expansion in a coordinate. Consequently, the implemented SDR procedure should be understood as an incumbent-centred panning--zooming method with adaptive side-length control, rather than as a strictly nested sequence of regions.

For the retraining BO--VAE algorithms, this interpretation is stagewise. After a VAE update, the latent coordinate system changes and the region is reinitialised as $\mathcal{R}_{(\ell;0)}=\mathcal{R}_0$. Hence, no nesting relation is imposed, or generally meaningful, between regions belonging to different VAE states.

\section{Additional Implementation Details}\label{Appendix: Additional Implementation Details}

In this section, we describe how the empirical encoded distribution is calibrated to the initial SDR region $[0, 1]^d$, enabling the standard SDR update to be applied in latent coordinates. 

Let $\mathcal{A}  =\{\mathbf{x}_i\}_{i=1}^{N}$ denote the dataset associated with the pretrained VAE. Let $q_\phi(\cdot|\mathbf{x})$ denote the corresponding encoder. The points in set $\mathcal{A}$ are mapped to the VAE latent space using the posterior mean of the encoder,
\[
\mathbf{z}_i
    \coloneq \mathbb E_{q_{\phi}(\cdot|\mathbf{x}_i)}[\mathbf{z}]
   \in\mathbb{R}^{d}.
\]
We define the empirical latent mean
\[
\mathbf{c}
   \coloneq \frac{1}{N}
    \sum_{i=1}^{N}\mathbf{z}_i,
\]
and the regularised covariance matrix
\[
\Sigma_{\varepsilon}
   \coloneq \frac{1}{N-1}
    \sum_{i=1}^{N}
    (\mathbf{z}_i-\mathbf{c})
    (\mathbf{z}_i-\mathbf{c})^{\top}
    +\varepsilon I_d,
    \qquad \varepsilon=10^{-6}.
\]
Let
\[
\Sigma_{\varepsilon}=LL^{\top}
\]
be its Cholesky factorisation. Then, we have the corresponding normalised latent coordinates
\[
\mathbf{w}_i
    \coloneq L^{-1}(\mathbf{z}_i-\mathbf{c}).
\]

We consider the empirical $99.5\%$ quantile of the normalised latent data to define the radius of the calibrated latent box as
\[
r
   \coloneq \max\left\{
      1,\,
      Q_{0.995}
      \left(
        \left\{
          \max_{1\leq j\leq d}
          |\mathbf{w}_{i}^{(j)}|
        \right\}_{i=1}^{N}
      \right)
    \right\},
\]
where $\mathbf{w}_{i}^{(j)}$ indicates the jth entry of $\mathbf{w}_i$. Hence, the calibrated latent search box is
\[
\mathcal{W}_0=[-r,r]^d.
\]
A point \(\mathbf{v}\in\mathcal{W}_0\) is mapped to the normalized SDR coordinates by
\[
\mathbf{u}
   \coloneq \frac{\mathbf{v}+r\mathbf{1}}{2r},
\]
where $\mathbf{1}$ denotes the vector of ones. We therefore initialise the SDR region as the unit hypercube $R_0=[0,1]^d$. The inverse transformation from an SDR coordinate \(u\in R_0\) to the raw VAE latent coordinate is
\begin{equation}\label{eq: normal latent to raw latent}
    \mathbf{v}=2r \mathbf{u} - r \mathbf{1}, \qquad \mathbf{z}=\mathbf{c}+L\mathbf{v}.
\end{equation}
The corresponding point in the original input space is obtained from the VAE decoder as usual.

At BO iteration \(k\),  we optimise the acquisition function to obtain the next iterate $u_{k+1} \in [0, 1]^d$. The point is then converted to the raw VAE latent coordinate using \eqref{eq: normal latent to raw latent} and evaluated through the decoder. The SDR region is updated periodically according to \(K_{\mathrm{SDR}}\), with every updated region satisfying $R_{k+1}\subseteq[0,1]^d.$

For the fixed-VAE algorithm, the latent calibration parameters \((\mathbf{c},L,r)\) remain unchanged throughout the optimization. For the VAE
retraining algorithms, an accepted retraining step changes the encoder and decoder. The dataset $\mathcal{D}_Z^{(l)}$ is therefore re-encoded and the calibration parameters \((\mathbf{c},L,r)\) are recomputed. The complete observation history is then re-encoded using the updated VAE, and the SDR region is reinitialised as $[0, 1]^d$ at the beginning of the BO loop. If a proposed VAE retraining step is rejected, the VAE parameters, latent
calibration, and current SDR region remain unchanged.

\section{Further Experimental Details}


\subsection{Experimental Setup for Comparing State-of-the-art BO Solvers}\label{Appendix: BO Solver Exp Setup Details}
The details are listed in Table~\ref{tab: BO solver table exp details}.
\begin{table}[h]
\caption{Details of the selected state-of-the-art BO solvers. $n_p$ denotes the problem dimension.}
\label{tab: BO solver table exp details}
\centering
\begin{tabular}{@{}l l l c c@{}}
\toprule
BO solver & Initial sampling strategy & Kernel & Acquisition function & Acquisition optimiser \\
\midrule
\textit{GPyOpt} \textit{v1.2.1}   & $2n_p$ uniform random points & Matérn-$\tfrac{5}{2}$ & EI & L-BFGS-B \\
\textit{BayesOpt} \textit{v1.5.1}  & $2n_p$ uniform random points & Matérn-$\tfrac{5}{2}$ & EI & L-BFGS-B \\
\textit{BoTorch} \textit{v0.11.1}  & $2n_p$ uniform random points & Matérn-$\tfrac{5}{2}$ & EI & L-BFGS-B \\
\botrule
\end{tabular}
\end{table}

\subsection{Experimental Configurations for Sections \ref{subsection: SDR in VAE-generated Latent Spaces}\mbox{-}\ref{subsection: Comparisons Between BO--VAE and BO--SDR Algorithms}}\label{Appendix Numerical Experiments Details}
The training details of the VAEs in Table~\ref{tab:CN-reference VAE family} are shown in Table~\ref{tab:cn-reference-vae-pretraining}. 

\begin{table}[h]
\caption{Common pretraining configuration for the VAEs in Table~\ref{tab:CN-reference VAE family}. Here, $M$ is the total number of objective-free pretraining samples. The tuple $(\beta_i,\beta_f,\beta_s,\beta_a)$ gives the initial and final KL weights, the annealing interval in epochs, and the increment applied at each interval, respectively.}
\label{tab:cn-reference-vae-pretraining}
\centering
\small
\setlength{\tabcolsep}{3pt}
\renewcommand{\arraystretch}{1.1}

\begin{tabularx}{\linewidth}{@{}l c c c c c c c @{}}
\toprule
VAE ID
& $M$
& \shortstack{Train/\\validation}
& Epochs
& Optimiser
& Learning rate
& \shortstack{Batch\\size}
& $(\beta_i,\beta_f,\beta_s,\beta_a)$ \\
\midrule

VAE--D10
& $50\,000$
& $45\,000/5\,000$
& 180
& Adam
& $10^{-3}$
& 512
& $(0,1,10,0.1)$ \\

VAE--D100
& $50\,000$
& $45\,000/5\,000$
& 180
& Adam
& $10^{-3}$
& 512
& $(0,1,10,0.1)$ \\

\bottomrule
\end{tabularx}
\end{table}

We use the Adam optimiser with an $L_2$ weight-decay coefficient of $10^{-5}$. This optimiser regularisation is distinct from the coefficient $\beta$ multiplying the KL-divergence term in the VAE objective. The pretraining data are sampled from a zero-mean truncated multivariate normal distribution on $[-5,5]^D$. Pretraining is objective-blind. The VAE is trained by minimising the reconstruction loss plus the $\beta$-weighted KL divergence. The value of $\beta$ is increased from $0$ to $1$ in increments of $0.1$ every ten epochs. To prevent checkpoint selection from favouring an early model trained with a weak KL penalty, the retained checkpoint is the model with the lowest validation loss among the epochs for which $\beta=1$.

The retraining details are given in Table~\ref{tab:cn-reference-vae-retraining}. The first retraining proposal is made when the initial labelled dataset is available. The initial gap between subsequent proposals is $50$ enrichment evaluations. If two accepted updates each produce less than a $1\%$ gain in the validation quality measure, the proposal gap is doubled, up to a maximum of $200$ evaluations. Conversely, a normalised incumbent improvement of at least $5\%$ resets the gap to $50$. A rejected proposal is rolled back, so that neither the VAE nor the latent coordinate system is changed. 

\begin{table}[h]
\caption{Online VAE retraining configurations used in the adaptive-retraining experiments. Each proposal is warm-started from the current accepted VAE and uses the complete accumulated labelled dataset. The KL weight is fixed at $\beta=1$ during retraining.}
\label{tab:cn-reference-vae-retraining}
\centering
\tiny
\setlength{\tabcolsep}{3pt}
\renewcommand{\arraystretch}{1.12}

\begin{tabularx}{\linewidth}{@{}l Y c c c c Y@{}}
\toprule
Method
& Proposal schedule
& \shortstack{Epochs per\\proposal}
& Optimiser
& Learning rate
& Weight decay
& DML parameters \\
\midrule

Baseline BO--VAE
& --
& --
& --
& --
& --
& -- \\

Retrained BO--VAE
& Acceptance-feedback gaps in $[50,200]$
& 5
& Adam
& $10^{-4}$
& $0$
& -- \\

Retrained BO--VAE with DML
& Acceptance-feedback gaps in $[50,200]$
& 5
& Adam
& $10^{-4}$
& $0$
& $(\lambda_{\mathrm{DML}},\eta, \nu)
   =(3,0.1,0.2)$ \\

\bottomrule
\end{tabularx}
\end{table}

For Section \ref{subsec:latent-dimension-representation-gap}, the pretraining settings remain those in Table~\ref{tab:cn-reference-vae-pretraining}, whereas the retraining settings are given separately in Table~\ref{tab:latent-dimension-retraining}.

\begin{table}[h]
\caption{Online retraining settings used specifically in the controlled latent-dimension study. These settings are separate from the adaptive policy in Table~\ref{tab:cn-reference-vae-retraining}.}
\label{tab:latent-dimension-retraining}
\centering
\tiny
\setlength{\tabcolsep}{3pt}
\renewcommand{\arraystretch}{1.1}

\begin{tabularx}{\linewidth}{@{}l l c c c c Y@{}}
\toprule
Setting & Method & $q$ & Epochs & Learning rate & Weight decay & DML settings \\
\midrule

$D=10,100$
& Retrained BO--VAE
& 100 & 2 & $3\times10^{-4}$ & $0$
& -- \\

\shortstack{$D=10$: Ackley,\\Rosenbrock}
& Retrained BO--VAE with DML
& 100 & 2 & $3\times10^{-4}$ & $0$
& $(\lambda_{\mathrm{DML}},\eta, \nu)
  =(0.1,0.02,0.2)$ \\

$D=10$: Rastrigin
& Retrained BO--VAE with DML
& 100 & 5 & $10^{-4}$ & $10^{-5}$
& $(0.3,0.01,0.1)$ \\

$D=100$
& Retrained BO--VAE with DML
& 100 & 2 & $3\times10^{-4}$ & $0$
& $(0.1,0.02,0.2)$ \\

\bottomrule
\end{tabularx}
\end{table}

We highlight three implementation details.

\begin{enumerate}
    \item The objective-free dataset $\mathcal D_{\mathbb U}$ is used only for VAE pretraining and latent-coordinate calibration. Its correlated truncated normal distribution is intended to provide structured ambient data from which the VAE can learn a compressed representation. Because no objective values are used, this stage does not consume the black-box evaluation budget.

    \item The labelled dataset $\mathcal D_{\mathbb L}^{(k)}$ contains only points at which the objective has actually been evaluated. These points are encoded using posterior means,
    \[
        \mathbf z_i
        =
        \mathbb E_{q_{\phi}(\mathbf z\mid\mathbf x_i)}[\mathbf z],
    \]
    to form the latent GP dataset
    \[
        \mathcal D_{\mathbb Z}^{(k)}
        =
        \{(\mathbf z_i,y_i):
        (\mathbf x_i,y_i)\in\mathcal D_{\mathbb L}^{(k)}\}.
    \]
    The latent-dimension study used $100$ initial observations and $800$ enrichment evaluations, whereas the subsequent full-rank, SDR-ablation, and profile experiments used $D$ initial observations and $800$ enrichment evaluations.

    \item Each retraining proposal uses the complete accumulated labelled dataset and is warm-started from the most recently accepted VAE. A proposed update is retained only if the observation-based quality and representation safeguards are satisfied. Following an accepted update, all previous observations are re-encoded, the latent transformation is recalibrated, the GP is rebuilt, and the SDR region is reinitialised in the updated coordinate system. If the proposal is rejected, the previous VAE and SDR state are retained.
\end{enumerate}

The SDR update hyperparameters are listed in Table~\ref{tab:sdr-experimental-settings}.

\begin{table}[h]
\caption{SDR and latent-coordinate settings used in the numerical experiments. The SDR regions are maintained in the standardised unit coordinates $\mathbf u$.}
\label{tab:sdr-experimental-settings}
\centering
\small
\setlength{\tabcolsep}{4pt}
\renewcommand{\arraystretch}{1.1}

\begin{tabularx}{\linewidth}{@{}l Y@{}}
\toprule
Parameter & Setting \\
\midrule

Latent reference set
& $45\,000$ objective-free VAE pretraining points \\

Covariance regularisation
& $\epsilon=10^{-6}$ \\

Coverage quantile
& $0.995$ \\

Initial SDR region
& $\mathcal R_0=[0,1]^d$ in standardised coordinates \\

Earliest SDR update
& After $10$ enrichment iterations \\

Latent-dimension study
& $K_{\mathrm{SDR}}=10$ \\

SDR-effect study
& $K_{\mathrm{SDR}}=10$ whenever SDR is enabled \\

Full-rank and profile experiments
& $K_{\mathrm{SDR}}=10$ for all BO--VAE variants \\

Accepted VAE update
& Refit the standardisation transformation, re-encode the complete labelled history, and reset $\mathcal R_k$ to $\mathcal R_0$ \\

Rejected VAE update
& Retain the current transformation and SDR region \\

\bottomrule
\end{tabularx}
\end{table}

\subsection{Experiment Configurations for Section \ref{sec:bovae-egorse-nonlinear}} \label{Appendix Exp set up for comparison against EGORSE}

Table~\ref{tab:curved-preimage-vae-pretraining} contains the pretraining details. Pretraining does not use any objective values or the nonlinear effective maps. The checkpoint with the smallest validation loss after $\beta$ has reached $\beta_f$ is selected. The pretraining data consist of $50\,000$ uniform samples from $[-1,1]^{10}$ with $45\,000/5\,000$ training/validation split.
\begin{table}[h]
\caption{Pretraining settings for the VAEs in Table~\ref{tab:curved-preimage-vae-family}. Here, $\beta_i$ and $\beta_f$ are the initial and final KL weights, and $\beta_a$ is added every $\beta_s$ epochs. The Adam weight decay is an optimiser regularisation and is distinct from the KL weight $\beta$.}
\label{tab:curved-preimage-vae-pretraining}
\centering
\small
\setlength{\tabcolsep}{4pt}
\renewcommand{\arraystretch}{1.1}

\begin{tabularx}{\linewidth}{@{}l c c c c c c@{}}
\toprule
VAE & Epochs & Optimiser & LR & Batch
& Weight decay & $(\beta_i,\beta_f,\beta_s,\beta_a)$\\
\midrule

VAE--E10
& 180
& Adam
& $10^{-3}$
& 512
& $10^{-5}$
& $(0,10^{-3},10,10^{-4})$ \\

VAE--E100
& 180
& Adam
& $10^{-3}$
& 512
& $10^{-5}$
& $(0,10^{-1},10,10^{-2})$\\

\bottomrule
\end{tabularx}
\end{table}

For the retraining setting applied to the VAEs in Table~\ref{tab:curved-preimage-vae-family}, we refer the reader to Table~\ref{tab:cn-reference-vae-retraining}. Finally, the hyperparameters for EGORSE are in Table~\ref{tab:reconstructed-egorse-settings}.

\begin{table}[h]
\caption{Settings of the EGORSE baseline. The original proprietary SEGOMOE/SNOPT implementation was unavailable. Therefore, these settings describe the reconstructed constrained-BO backend used to produce the reported results.}
\label{tab:reconstructed-egorse-settings}
\centering
\small
\renewcommand{\arraystretch}{1.1}

\begin{tabularx}{\linewidth}{@{}l Y@{}}
\toprule
Component & Setting \\
\midrule
Displayed variants
& EGORSE $d_e$ and EGORSE $2d_e$, with optimiser dimension
  $m=d_e$ and $m=2d_e$, respectively \\

Embedding sequence
& PLS followed by a Gaussian random embedding \\

Subspace block length
& $20m$ enrichment evaluations per embedding \\

Maximum outer cycles
& $10$, subject to the overall $800$-evaluation budget \\

Initial design
& The same $D$ uniformly sampled ambient points used by BO--VAE \\

Reduced search box
& $\mathcal{B}_j=[-\sum_k|A_{jk}|,\sum_k|A_{jk}|]$ \\

Backward map
& Constrained quadratic programme using CVXOPT, with SciPy SLSQP
  as a numerical fallback; clipped pseudoinverse projection is used
  when the equality-constrained map is infeasible \\

Objective and constraint surrogates
& SMT KRG with squared-exponential correlation and constant trend \\

Initial KRG length-scale
& $\theta_0=10^{-3}$ \\

KRG multistarts
& $1$ \\

Nugget sequence
& $2.22\times10^{-14},10^{-13},10^{-12},10^{-11},
  10^{-10},10^{-8}$ \\

Acquisition
& Constrained expected improvement using the surrogate mean constraint \\

Initial acquisition screen
& $1024$ random candidates and the centre of the reduced box \\

Global acquisition optimiser
& NLOPT ISRES, maximum $256$ evaluations \\

Local refinement
& SciPy SLSQP from $4$ starts, maximum $80$ iterations,
  $\mathrm{ftol}=10^{-9}$ \\

Backward-map tolerance
& $10^{-7}$ \\

Duplicate tolerance
& $10^{-8}$ \\

Arithmetic and device
& CPU execution with one BLAS thread per worker \\
\bottomrule
\end{tabularx}
\end{table}

We highlight the following implementation details.

\begin{enumerate}
    \item The VAEs used in this comparison are pretrained without objective values or knowledge of the nonlinear effective map. Unlike the experiments in Sections \ref{subsec:latent-dimension-representation-gap}--\ref{sec:bovae-benchmark-performance}, the pretraining data are sampled uniformly from $[-1,1]^D$. The scaled-$\tanh$ decoder matches this bounded ambient domain.

    \item For every $(D,d)$ pair, the three BO--VAE variants start from the same pretrained checkpoint. Consequently, differences
    among the variants arise from their online representation updates and the DML term rather than from different initial VAEs.

    \item The VAE pretraining set is not used as GP training data. It is used only for VAE pretraining and latent calibration. The GP is fitted to all objective-labelled observations collected by BO.

    \item An accepted VAE update changes the meaning of the latent coordinates. The complete observation history is therefore re-encoded, the GP is rebuilt, and SDR is reset to the full calibrated unit box. A rejected proposal is rolled back without modifying the existing latent calibration or SDR region.
\end{enumerate}

\FloatBarrier

\subsection{Solver Results}\label{appendix: solver results}
\begin{figure}[h]
    \centering
    \includegraphics[width=\linewidth]{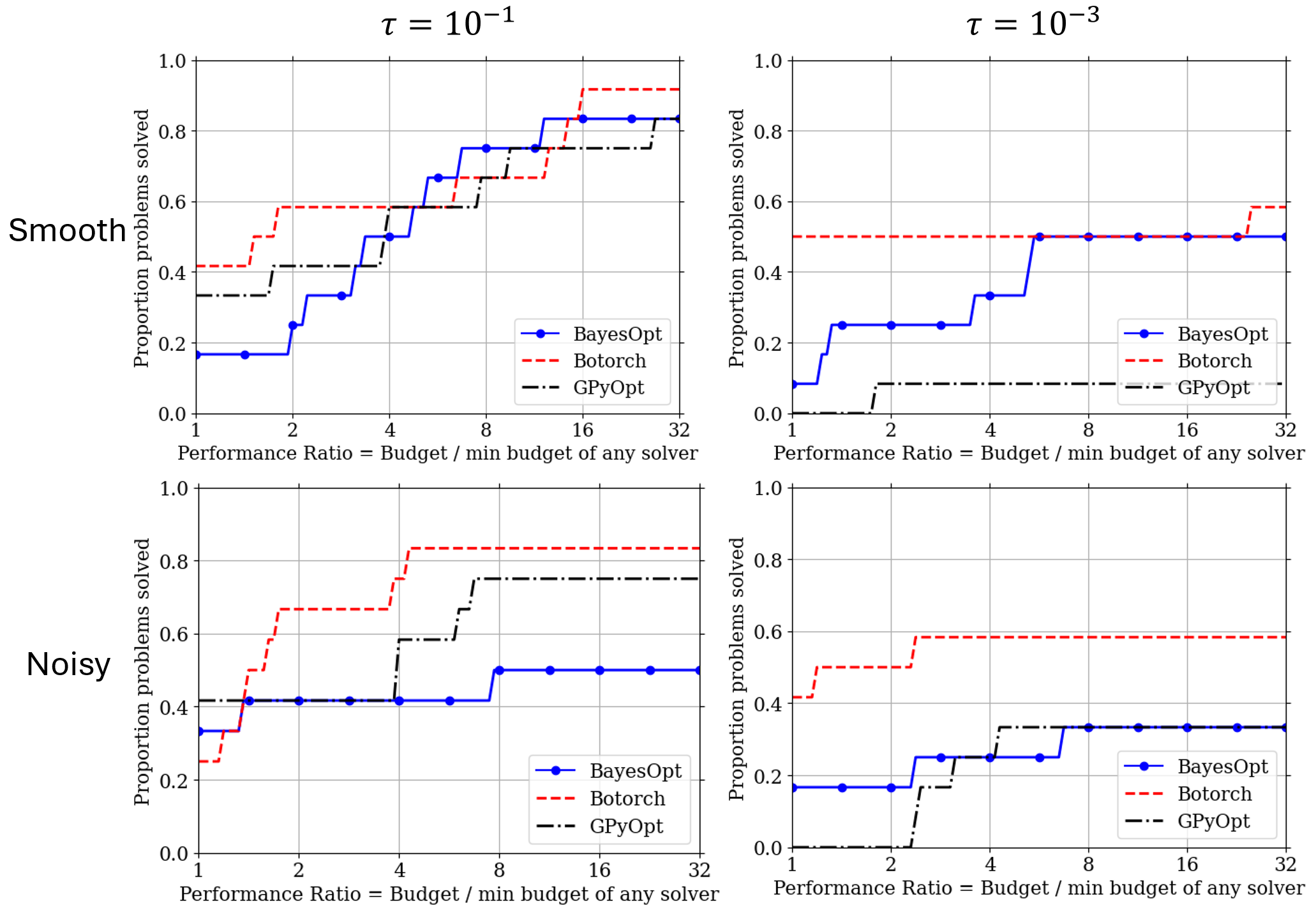}
    \caption{Performance profiles of \textit{GPyOpt, BayesOpt,} and \textit{BoTorch} on the benchmark test problems listed in Table \ref{tab:solver-test}.}
    \label{fig: BOS comparison}
\end{figure}

We compare \textit{GPyOpt}, \textit{BoTorch}, and \textit{BayesOpt} for noisy and smooth $f$. When the function evaluation is noisy, we mean that additive Gaussian noise is added to $f$, i.e., 
\begin{equation}\label{noisy problem}
    \Tilde{f}(\mathbf{x}) = f(\mathbf{x}) + \sigma \epsilon, 
\end{equation}
where $\epsilon \sim \mathcal{N}(0, 1)$ i.i.d. for each $\mathbf{x}$ and each solver. Here, we set $\sigma$ to be $10^{-2}.$  The results as a function of the number of objective evaluations are presented in Table \ref{tab:solver data-profiles} and Figure \ref{fig: BOS comparison}. The selected accuracy levels \( \tau \) are set to \( 10^{-1} \) and \( 10^{-3} \). Each problem in Table \ref{tab:solver-test} is configured with \( N_g = 50 \) and repeated once. Therefore, the effective number of test instances for this comparison experiment is $12$. The experimental settings of the BO solvers are given in Table~\ref{tab: BO solver table exp details}.

\begin{table}[h]
\caption{Average percentage of problems solved from Table~\ref{tab:solver-test} at two tolerances $\tau$ (higher is better).}
\label{tab:solver data-profiles}
\centering
\begin{tabular}{@{}ccccccc@{}}
\toprule
& \multicolumn{3}{c}{$\tau=10^{-1}$ } & \multicolumn{3}{c}{$\tau=10^{-3}$ } \\
\cmidrule(lr){2-4}\cmidrule(lr){5-7}
Regime & {BoTorch} & {BayesOpt} & {GPyOpt} & {BoTorch} & {BayesOpt} & {GPyOpt} \\
\midrule
Smooth & 100   & 83.33 & 83.33 & 58.33 & 66.67 & 8.33  \\
Noisy  & 83.33 & 75    & 75    & 50    & 58.33 & 33.33 \\
\botrule
\end{tabular}
\end{table}

As illustrated in the figure, in the noiseless setting, in the low-accuracy case \( \tau = 10^{-1} \), \textit{BoTorch} and \textit{BayesOpt} outperform \textit{GPyOpt}, with \textit{BoTorch} solving a higher percentage of problems and demonstrating better performance compared to \textit{BayesOpt}. For \( \tau = 10^{-3} \), although \textit{BayesOpt} ultimately solves more problems, \textit{BoTorch} solves more problems at smaller evaluation budgets and exhibits slightly better performance, as shown in the top-right plot. In the noisy setting, \textit{BoTorch} clearly outperforms the other two solvers, solving more problems with limited budgets.

Based on this comparison, \textit{BoTorch} and \textit{BayesOpt} are the strongest candidates among the three solvers considered. \textit{BoTorch} performs better in the noisy setting, offers greater flexibility for custom datasets, and supports GPU execution, making it the more suitable choice for our implementation.

\end{appendices}

\bibliography{sn-bibliography}

\end{document}